\documentclass[nosumlimits,twoside,reqno]{amsart}
\usepackage{amssymb}
\usepackage{amsmath}
\usepackage[bookmarksnumbered,colorlinks,plainpages,hyperindex]{hyperref}
\usepackage{amsfonts}
\usepackage{color}
\usepackage{float}
\usepackage{graphicx}
\usepackage{xy}
\xyoption{all}
\newtheorem{theorem}{\sc Theorem}[section]
\newtheorem{proposition}[theorem]{\sc Proposition}
\newtheorem{lemma}[theorem]{\sc Lemma}
\newtheorem{corollary}[theorem]{\sc Corollary}
\theoremstyle{definition}
\newtheorem{definition}[theorem]{\sc Definition}

\newtheorem{example}[theorem]{\sc Example}

\theoremstyle{remark}
\newtheorem{remark}[theorem]{\sc Remark}
\newtheorem{remarks}[theorem]{\sc Remarks}

\newtheoremstyle{mystyle}{}{}{}{}{}{.}{ }
  {\thmname{#1}\thmnote{\scshape #3}}
\theoremstyle{mystyle}

\newcommand{\bs}[1]{{\boldsymbol{#1}}}\newcommand{\cs}[1]{{\mathcal{#1}}}
\newcommand{\mon}[1]{(\boldsymbol{#1},\boldsymbol{m_{#1}},\boldsymbol{u_{#1}})}

\def\I{{\boldsymbol{I}}}
\def\O{{\boldsymbol{O}}}

\def\wotR{\widehat{\otimes}_R\,}

\def\calM{\frak{M}}

\def\calM{\mathcal{M}}

\def\ot{\otimes}

\def\stac#1{\raise-.2cm\hbox{$\,\stackrel{\displaystyle\otimes}{\scriptscriptstyle{#1}}\,$}}
\def\cten#1{\raise-.2cm\hbox{$\stackrel{\displaystyle\widehat{\otimes}}
{\scriptscriptstyle{#1}}$}}
\numberwithin{equation}{section}

\def\I{{\boldsymbol{T_r}}}
\def\O{{\boldsymbol{T_l}}}
\def\calS{\mathcal{S}}
\def\calC{\mathcal{C}}

\def\rmr{R\text{-}\mathrm{Mod}\text{-}R}

\def\calD{\mathcal{S}}
\def\calW{\mathcal{W}}
\begin{document}

\title{(Co)cyclic (co)homology of bialgebroids: An approach via (co)monads}
\author{Gabriella B\"ohm}
\author{Drago\c s \c Stefan}
\address{Research Institute for Particle and Nuclear Physics, Budapest,
H-1525 Budapest 114, P.O.B.49, Hungary.}
\address{University of Bucharest, Faculty of Mathematics and Informatics,
Bucharest, 14 Academiei Street, Ro-010014, Romania.}
\date{July 2007}

\begin{abstract}
For a (co)monad $\boldsymbol{T_l}$ on a category $\mathcal{M}$, an object $X$
in $\mathcal{M}$, and a functor $\bs{\Pi}:\mathcal{M}\to \mathcal{C}$, there
is a (co)simplex $Z^\ast:=\bs{\Pi} \boldsymbol{T_l}^{\ast +1} X$ in
$\mathcal{C}$. The
aim of this paper is to find criteria for para-(co)cyclicity of $Z^\ast$. Our
construction is built on a distributive law of $\boldsymbol{T_l}$ with a
second (co)monad $\boldsymbol{T_r}$ on $\mathcal{M}$, a natural transformation
$\boldsymbol{i}:\bs{\Pi} \O \to \bs{\Pi} \I$, and a morphism
$w:\boldsymbol{T_r}X \to \boldsymbol{T_l}X$ in $\mathcal{M}$. The
(symmetrical) relations $\boldsymbol{i}$ and $w$ need to satisfy are
categorical versions of Kaygun's axioms of a transposition map.
Motivation comes from the observation that a (co)ring $T$ over an algebra $R$
determines a distributive law of two (co)monads $\boldsymbol{T_l}=T \otimes_R
(-)$ and $\boldsymbol{T_r}=(-)\otimes_R T$ on the category of $R$-bimodules.
The functor $\bs{\Pi}$ can be chosen such that $Z^n= T\widehat{\otimes}_R\dots
\widehat{\otimes}_R T \widehat{\otimes}_R X$ is the cyclic $R$-module tensor
product. A natural transformation $\boldsymbol{i}:T \widehat{\otimes}_R (-)
\to (-) \widehat{\otimes}_R T$ is given by the flip map and a morphism $w: X
\otimes_R T \to T\otimes_R  X$ is constructed whenever $T$ is a
(co)module algebra or coring of an $R$-bialgebroid. The notion of a stable
anti Yetter-Drinfel'd module over certain bialgebroids, so called
$\times_R$-Hopf algebras, is introduced. In the particular example when $T$ is
a module coring of a $\times_R $-Hopf algebra $\mathcal{B}$ and $X$ is a
stable anti Yetter-Drinfel'd $\mathcal{B}$-module, the para-cyclic object
$Z_\ast$ is shown to project to a cyclic structure on $T^{\otimes_R\,
  \ast+1} \otimes_{\mathcal{B}} X$. For a $\mathcal{B}$-Galois extension
$S \subseteq T$, a stable anti Yetter-Drinfel'd
$\mathcal{B}$-module $T_S$ is constructed, such that the cyclic
objects $\mathcal{B}^{\otimes_R\, \ast+1} \otimes_{\mathcal{B}}
T_S$ and $T^{\widehat{\otimes}_S\, \ast+1}$ are isomorphic. This
extends a theorem by Jara and \c Stefan for Hopf Galois
extensions. As an application, we compute Hochschild and cyclic
homologies of a groupoid with coefficients in a stable anti
Yetter-Drinfel'd module, by tracing it back to the group case. In
particular, we obtain explicit expressions for (coinciding
relative and ordinary) Hochschild and cyclic homologies of a
groupoid. Latter extends results of Burghelea on cyclic homology
of groups.
\end{abstract}
\maketitle
\tableofcontents

\section*{Introduction}

Cyclic cohomology of Hopf algebras is originated from the work
\cite{ConMos:HopfCyc} of Connes and Moscovici on the index theory of
transversally elliptic operators. Their local index formula in
\cite{ConMos:loc_index} gives a generalization of the Chern character to
non-commutative geometry. In order to give a geometrical interpretation of the
non-commutative Chern character in terms of non-commutative foliations, in
\cite{ConMos:HopfCyc} a cocyclic structure was constructed on a cosimplex
$Z^{CM}_n = H^{\otimes n}$, associated to the coalgebra underlying a Hopf
algebra $H$ over a field $\mathbb{K}$. The cocyclic operator was given in
terms of a so called modular pair in involution. 

In the subsequent years the Connes-Moscovici cocyclic module was placed in a 
broader and broader context. In \cite{KhaRan:InvInvCycHom} (see also
\cite{HaKhRaSo2}) to any (co)module algebra ${T}$ of a Hopf algebra $H$, and
any $H$-(co)module $X$, there was associated a para-cyclic module with
components {${T}^{\otimes\ast+1}\otimes X$}. Dually, for any (co)module
coalgebra ${T}$ of a Hopf algebra $H$, and any $H$-(co)module $X$, there is a
para-cocyclic module with components {${T}^{\otimes\ast+1}\otimes X$}. The
Connes-Moscovici cosimplex $Z^{CM}_\ast$ turns out to be isomorphic to a
quotient of the para-cocyclic module associated to the regular module
coalgebra ${T}:= H$ and an $H$-comodule defined on $\mathbb{K}$.    
For bialgebras, the Connes-Moscovici construction was generalized  in
\cite{Kay}. 

In the papers \cite{HaKhRaSo1} and \cite{JaSt:CycHom}, a modular pair in
involution was proven to be equivalent to a stable anti Yetter-Drinfel'd
module structure on the ground field $\mathbb{K}$. In \cite{HaKhRaSo2}, the
para-cocyclic module ${T}^{\otimes\ast+1}\otimes X$, associated to an
$H$-module coalgebra ${T}$ and a stable anti Yetter-Drinfel'd $H$-module $X$,
was shown to project to a cocyclic object whose components are the $H$-module
tensor products ${T}^{\otimes\ast+1}\otimes_H X$. The way in which the
para-cocyclic object $H^{\otimes  \ast+1}$ projects to the Connes-Moscovici
cosimplex $Z^{CM}_\ast$, is an example of this scenario. 
Dually, the para-cyclic module, associated to an $H$-comodule algebra, was
proven to have a cyclic submodule. 

In the spirit of \cite{JaSt:CycHom}, one can follow a dual approach. That is,
para-cocyclic modules can be constructed for (co)module algebras of Hopf 
algebras, and para-cyclic modules for (co)module coalgebras, in both cases
with coefficients in $H$-(co)modules. Taking coefficients in a stable anti
Yetter Drinfel'd module, it was shown in \cite{JaSt:CycHom} that
in this case the para-cyclic object associated to a module coalgebra possesses
a cyclic quotient. 
In \cite{KR} an isomorphism was proven between the cyclic quotient of the
para-cyclic object in \cite{JaSt:CycHom} of $H$ as an $H$-module coalgebra on 
one hand, and the cyclic subobject of the para-cyclic object in \cite{HaKhRaSo2}
of $H$ as an $H$-comodule algebra on the other.

Constructions in Section \ref{sec:bgd} of the current paper follow
the root in \cite{JaSt:CycHom}. Since this framework is dual to
that suggested in \cite{HaKhRaSo2} (cf. also \cite{Kay:UniHCyc}), some might
like to call it a {\em dual} Hopf (co)cyclic theory. However, we do not use
this somewhat involved terminology in the paper, but remind the reader to the
difference between the two possible dual approaches.  

In \cite{Kay:UniHCyc} Kaygun proposed a unifying approach to the
para-(co)cyclic objects corresponding to a (co)module (co)algebra
of a Hopf algebra. Starting with a (co)algebra ${T}$ and an
object $X$ in a symmetric monoidal category $\mathcal{S}$, he
introduced the notion of a transposition map. It is a morphism
$w:X \otimes {T}\to {T}\otimes X$ in $\mathcal{S}$, satisfying
conditions reminiscent to half of the axioms of a distributive law
in \cite{Be}. Any transposition map $w$ was shown to determine a
para-(co)cyclic structure on the (co)simplex {
$T^{\otimes\ast+1}\otimes X$} in $\mathcal{S}$. In particular,
canonical transposition maps were constructed for (co)module
(co)algebras ${T}$ and (co)modules $X$ of a bialgebra.

Connes and Moscovici's index theory of transversally elliptic
operators lead beyond cyclic homology of Hopf algebras. In dealing
with the general, non-flat case, in \cite{ConMos:DiffCyc} certain
bialgebroids (in fact $\times_{R}$-Hopf algebras) arose
naturally. Bialgebroids can be thought of as a generalization of
bialgebras to a non-commutative base algebra ${R}$, while
$\times_{R}$-Hopf algebras generalize Hopf algebras. There are
a few papers in the literature, e.g. \cite{KhaRan:Para_Hopf} and
\cite{Rang:CycCor}, attempting to extend Hopf cyclic theory to
non-commutative base algebras. However, an understanding of the
subject, comparable to that in the classical case of a commutative
base ring (or field), is missing yet. The aim of the current paper
is to give a universal construction of para-(co)cyclic
(co)simplices, including examples coming from (co)module algebras
and (co)module corings for bialgebroids.

When replacing bialgebras over a commutative ring $\mathbb{K}$ by
bialgebroids over a non-commutative $\mathbb{K}$-algebra ${R}$, the monoidal
category of $\mathbb{K}$-modules becomes replaced by the monoidal category of
${R}$-bimodules. Indeed, (co)module algebras of an ${R}$-bialgebroid are in
particular algebras, and (co)module corings are coalgebras, in the category of
$R$-bimodules. The main difference is that the category of
$\mathbb{K}$-modules is symmetric. In contrast, the category of
${R}$-bimodules is not even braided in general. Hence Kaygun's elegant theory
\cite{Kay:UniHCyc}, formulated in a symmetric monoidal category $\mathcal{S}$,
is not applicable.

Our key observation is that the role, the symmetry plays in Kaygun's work, is
  that it defines a compatible natural transformation $\bs{i}$ between
  the two (co)monads ${T}\otimes (-)$ and $(-)\otimes {T}$ on
  the symmetric monoidal category $\mathcal {S}$, induced by a (co)algebra
  ${T}$ in $\mathcal{S}$. Note
  that these (co)monads on $\mathcal {S}$ are connected by a trivial
  distributive law.
Guided by this observation, in Section \ref{sec:adm.sep} we start with a
  distributive law of two (co)monads $\O$ and $\I$ on any category
  $\mathcal{M}$. In addition, we allow for the presence of a functor
  $\bs{\Pi}:\mathcal{M} \to \mathcal {C}$ (it is the identity functor on
  $\mathcal{S}$ in \cite{Kay:UniHCyc}). Then, for any object $X$ in
  $\mathcal{M}$, there is a (co)simplex in $\mathcal{C}$, given at degree $n$
  by $\bs{\Pi} \O^{n+1} X$. In Sections \ref{sec:main} and \ref{sec:dual} we
  show that it is para-(co)cyclic provided that
  there exist a natural transformation $\bs{i}: \bs{\Pi} \O\to \bs{\Pi} \I$
  and a morphism $\I X \to \O X$ in $\mathcal{M}$, satisfying symmetrical
  conditions generalizing Kaygun's axioms of a transposition map. Examples of
  this situation are collected in Section \ref{sec:ex_adm_sep}.
  Among other (classical) examples, we show that \v Skoda's functorial
  construction in \cite{Sko:CycCom} of a para-cyclic object
  in the category of endofunctors,
  Majid and Akrami's para-cyclic modules associated to a ribbon algebra in
  \cite{AM}, and Rangipour's cyclic module in \cite{Rang:CycCor} determined
  by a coring, fit our framework. It is
  discussed in Sections \ref{sec:(co)mod.alg} and \ref{sec:(co)mod_coring} how
  the general results in Sections \ref{sec:main} and \ref{sec:dual} cover the
  particular cases when the two (co)monads
  $\O={T}\otimes_{R} (-)$ and $\I=(-)\otimes_{R} {T}$ are induced by a
  (co)module algebra or
  (co)module coring ${T}$ of an ${R}$-bialgebroid $\mathcal{B}$, the
  functor $\bs{\Pi}$ is defined via the coequalizer of the ${R}$-actions in a
  bimodule, and $X$ is a $\mathcal{B}$-(co)module. The components of the
  resulting para-(co)cyclic module are cyclic ${R}$-module tensor products
  ${T}\widehat{\otimes}_{R} \dots \widehat{\otimes}_{R} {T}
  \widehat{\otimes}_{R} X$. {In this way we obtain examples which extend
    both some para-(co)cyclic objects in \cite{Kay:UniHCyc} and
    \cite{HaKhRaSo2} and their cyclic duals.}

{By the above procedure, we associate a para-cyclic object in particular to
a module coring ${C}$ and a comodule $X$ of an ${R}$-bialgebroid
$\mathcal{B}$. Following \cite{JaSt:CycHom},}
in Section \ref{sec:st.a.Yet.Dri} we look
for situations in which it projects to the $\mathcal{B}$-module
tensor product $\big( {C} \otimes_{R} \dots \otimes_R
{C}\big)\otimes_{\mathcal{B}} X$. Restricting at this point
our study to $\times_{R}$-Hopf algebras $\mathcal{B}$, we
define stable anti Yetter-Drinfel'd modules for $\mathcal{B}$. In
parallel to the case of Hopf algebras \cite[Theorem
  4.13]{JaSt:CycHom}, \cite[Theorem 2.1]{HaKhRaSo2}, we prove
cyclicity of the simplex $\big( {C} \otimes_{R} \dots
\otimes_{R} {C}\big)\otimes_{\mathcal{B}} X$, whenever $X$
is a  stable anti Yetter-Drinfel'd module.

A simplest example of a cyclic simplex is associated to an algebra
extension
  ${S}\subseteq {T}$. Its components are given by the $n+1$ fold cyclic
  tensor product ${T}^{\widehat{\otimes}_{S}\, n+1}$, face and
  degeneracy maps
  are determined by the algebra structure of ${T}$ and the cyclic map is
  given by the cyclic permutation of the tensor factors.
In Section \ref{sec:Galois}, for a Galois extension ${S}\subseteq {T}$
  by a $\times_{R}$-Hopf
algebra $\mathcal{B}$, we construct a stable anti Yetter-Drinfel'd
module ${T}_{S}:= {T}/\{\ s\cdot t - t\cdot s \ |\
s\in {S},\ t\in {T}\ \}$. We prove that the cyclic
simplices ${T}^{\widehat{\otimes}_{S}\,
  n+1}$ and $\mathcal{B}^{\otimes_{R}\, n+1}\otimes_{\mathcal{B}} {
  T}_{S}$ are isomorphic. This extends \cite[Theorem 3.7]{JaSt:CycHom}.

A most fundamental class of examples of bialgebroids (in fact $\times_R$ Hopf
algebras) is given by {algebras (over fields), generated by a groupoid of
  finitely many objects}. As an application of
our abstract theory, we compute explicitly the relative Hochschild and cyclic
homologies of {such} a groupoid, with coefficients in a stable anti
Yetter-Drinfel'd module. By our results, any Galois extension by the groupoid
provides us with a stable anti Yetter-Drinfel'd module. In particular, the
groupoid algebra ${\mathcal B}$ is a Galois extension of its base algebra
$R$. Applying the isomorphism of the simplices ${{\mathcal
    B}}^{\widehat{\otimes}_{R}\, n+1}$ and ${\mathcal B}^{\otimes_{R}\,
  n+1}\otimes_{\mathcal{B}} {{\mathcal B}}_{R}$, we obtain the $R$-relative
cyclic homology of ${\mathcal B}$. Since $R$ is a separable algebra, it is
equal to ordinary cyclic homology of ${\mathcal B}$, hence our results extend
those by Burghelea on the cyclic homology of groups \cite{Burg}. Similar
formulae were obtained by Crainic for cyclic homology of {\em \'etale}
groupoids \cite{Cra}. {Observe that any groupoid (with arbitrary set of
objects) can be obtained as a direct limit of groupoids with finite sets of
objects. Certainly, the algebra generated by a groupoid with infinitely many
objects is no longer unital. However, one can still consider its cyclic
homology, as a homology of Connes' complex, associated to a presimplicial
object. Since the homology functor commutes with direct limits, we can extend
our formula of cyclic homology to arbitrary groupoids.}

Throughout the paper $\mathbb{K}$ denotes a commutative and associative unital
ring. The term \emph{$\mathbb{K}$-algebra} means an associative and unital
algebra over $\mathbb{K}$.

\section{The (co)cyclic object associated to a transposition map}
\label{sec:adm.sep}

In this first section we establish a general categorical framework -- in terms
of {\em admissible septuples} and their {\em transposition maps} -- to produce
para-cocyclic, and dually, para-cyclic objects.

\subsection{Notation and conventions}

In the $2$-category \textrm{CAT } we denote horizontal composition (of
functors) by juxtaposition, while $\circ $ is used for
vertical composition (of natural transformations). That is, for
two functors $\bs{F}:{\mathcal C}\to {\mathcal C}'$,
$\bs{G}:{\mathcal
    C}'\to {\mathcal C}''$ and an object $X$ in ${\mathcal C}$, instead
of $\bs{G}(\bs{F}(X))$ we write $\bs{GF}X$. For two natural
transformations $\bs{\mu}:\bs{F}\rightarrow \bs{F}^{\,\prime }$
and $ \bs{\nu}:\bs{G}\rightarrow \bs{G}^{\,\prime }$ we write
$\bs{G}^{\,\prime }\bs{\mu} X\circ \bs{\nu F}X:\bs{GF}X\rightarrow
\bs{G}^{\,\prime }\bs{F}^{\,\prime }X$ instead of
$\bs{G}^{\,\prime }(\bs{\mu}_ X)\circ \bs{\nu}_{\bs{F}(X)}.$ In
equalities of natural transformations we shall omit the object $X$ in our
formulae.

Inspired by the diagrammatic computation in a 2-category (in
particular \textrm{CAT}), we shall use a graphical representation
of morphisms in a category. For functors
$\bs{F}_1,\dots,\bs{F}_n,\bs{G}_1,\dots,\bs{G}_m$, which can be
composed to $\bs{F}_1\bs{F}_2\dots \bs{F}_n:{\mathcal D}_1\to
{\mathcal C}$ and $\bs{G}_1\bs{G}_2\dots \bs{G}_m:{\mathcal
D}_2\to {\mathcal C}$, and objects $X$ in ${\mathcal D}_1$ and $Y$
in ${\mathcal D}_2$, a morphism $ f:\bs{F}_1\bs{F}_2\dots \bs{F}_n
X\to \bs{G}_1\bs{G}_2\dots \bs{G}_mY$ will be represented
vertically, with the domain up, as in Figure~\ref{Fig:morphisms in
C}(a). Furthermore, for a functor
$\bs{T}:\mathcal{C}\to\mathcal{C}'$, the morphism $\bs{T}f$ will
be drawn as in (b). Keeping the notation from the first paragraph
of this section, the picture representing $\bs{\mu G}X$ is shown
in diagram (c). The composition $g\circ f$ of the morphisms
$f:X\to Y$ and $g:Y\to Z$ will be represented as in diagram (d).
For the multiplication $\bs{m_T}$ and the unit $\bs{u_T}$ of a
monad $\bs{T}$ on $\mathcal{C}$ (see Definition~\ref {de:monad}),
and an object $X$ in $\mathcal{C}$, to draw $\bs{m_T}X$ and
$\bs{u_T} X$ we shall use the diagrams (e) and (f), while for a
distributive law $\bs{t}:\bs{RT}\to \bs{TR}$ (see
Definition~\ref{de:distributivity law}) $\bs{t} X$ will be drawn
as in the picture (g). If $\bs{t}$ is invertible, the
representation of $\bs{t}^{-1}X$ is shown in diagram (h).
\begin{figure}[h]
\begin{center}
{\includegraphics[scale=1]{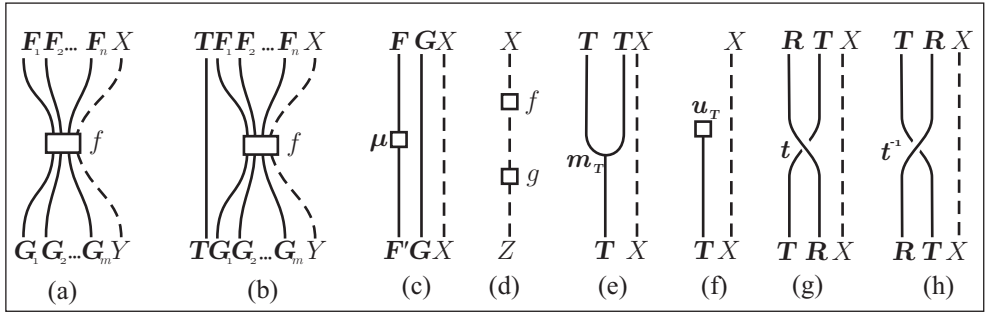}}
\end{center}
\caption{Diagrammatic representation of morphisms in a category}
\label{Fig:morphisms in C}
\end{figure}
For simplifying diagrams containing only natural transformations, we
shall always omit the last string that corresponds to an object in the
category. That is, we work with diagrams in \textrm{CAT} whenever it is
possible.

We shall use the following method to perform computations with such
diagrams. In view of associativity of composition, any diagram, representing
a well-defined composition of morphisms, can be thought of as a tower with
several layers. Any part of the diagram, corresponding to a layer, can be
substituted with any other equivalent representation of it. Usually,
equivalent representations are obtained from formulas that define the
notions that we deal with, or equations that have been previously proved.

\subsection{Monads and distributive laws}
Monads represent the main ingredient in our approach to cyclic
(co)homology. The definition of monads traces back to Godement's
book \cite{Go}, where they are called ``standard constructions''.
In the literature they are also called ``triples'', see for
example \cite{EiM}.

\begin{definition}
\label{de:monad} A {\em monad} on a category $\mathcal{C}$ is a
triple $(\bs{T},\bs{m_T},\bs{u_T} ),$ where $
\bs{T}:\mathcal{C}\rightarrow \mathcal{C}$ is a functor,
$\bs{m_T}:\bs{T}^{2}\rightarrow \bs{T}$ and $\bs{u_T}
:\bs{Id}_{\mathcal{C}}\rightarrow \bs{T}$ are natural
transformations such that the first two diagrams in
Figure~\ref{fig:monad}, expressing associativity and unitality,
are commutative.
\begin{figure}[h]
\begin{center}
\fbox{
\xymatrix{
  \bs{T}^3 \ar[d]_{\bs{T m_T}} \ar[r]^{\bs{m_TT}} & \bs{T}^2 \ar[d]^{\bs{m_T}}
  \\
  \bs{T}^2\ar[r]^{\bs{m_T}} &  \bs{T} }\
\xymatrix{ &  \bs{T} \ar[d]^{\bs{{Id}_T}} \ar[dl]_{\bs{u_TT}}
  \ar[dr]^{\bs{Tu_T} }  & \\
 \bs{T}^2\ar[r]_{\bs{m_T}}&  \bs{T}   & \bs{T}^2 \ar[l]^{\bs{m_T}}
 }\quad
\xymatrix{
  \bs{T}^2 \ar[d]^{\bs{m_T}} \ar[rr]^{\bs{\varphi T'}\circ \bs{T\varphi}} &&
  \bs{T'}^2 \ar[d]_{\bs{m_{T'}}} \\
  \bs{T}\ar[rr]_{\bs{\varphi}} && \bs{ T'} }\
\xymatrix{
                &\bs{Id}_{\mathcal{C}} \ar[ld]_{\bs{u_T} }\ar[dr]^{\bs{u_{T'}}
  } \\
 \bs{T}  \ar[rr]_{{\varphi}} & &    \bs{ T'}       }}
\end{center}
\caption{Monads and morphisms of monads.}
\label{fig:monad}
\end{figure}
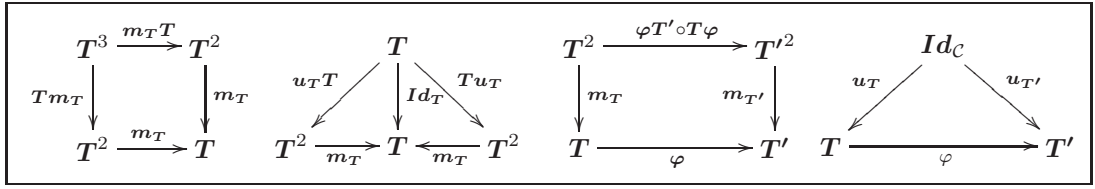
\noindent We call $\bs{m_T}$ and $\bs{u_T}$ the
\emph{multiplication} and the \emph{unit} of the monad $\bs{T}$,
respectively.

For two monads $\mon{T}$ and $\mon{T'}$ on $\mathcal{C}$, we say
that a natural transformation $\bs{\varphi} :\bs{T}\rightarrow
\bs{T'}$ is a \emph{morphism of monads} if the last two diagrams
in Figure~\ref{fig:monad} are commutative.
\end{definition}

\begin{example}\label{ex:k_alg}
Let $(\mathcal{C},\otimes,\bs{a},\bs{l},\bs{r},\mathbf{1})$ be a
monoidal category with unit object $\mathbf{1}$, associativity
constraint $\bs{a}$ and unit constraints $\bs{l},\bs{r}$. For
details about monoidal categories the reader is referred to
\cite[Chapter XI]{Ka}. An algebra in $\mathcal{C}$ is a triple
$(T,m_T,u_T)$ such that $m_T:T\otimes T\rightarrow T$ defines an
associative multiplication on $T$ with unit
$u_T:\mathbf{1}\rightarrow T$. To such an algebra one can
associate two monads $\bs{T_l}:= { T} \otimes (-)$ and
$\bs{T_r}:=  (-)\otimes{ T}$ on $\cs{C}$. The multiplication
$\bs{m_{T_l}}$ and the unit $\bs{u_{T_l}}$ of $\bs{T_l}$ are given
by
 \[\bs{m_{T_l}}X :=(m_{T}\otimes X)\circ\bs{a}_{T,T,X}^{-1}\qquad
 \text{and}\qquad
 \bs{u_{T_l}}X :=(u_{T}\otimes X)\circ\bs{l}_X^{-1},
 \]
for every $X$ in $\calC$. Analogously, for  $X$ in $\calC$,
$\bs{m_{T_r}}X$ and $\bs{u_{T_r}}X$ are defined by
 \[
 \bs{m_{T_r}}X :=(X\otimes m_{T})\circ \bs{a}_{X,T,T}\qquad
 \text{and} \qquad \bs{u_{T_r}}X :=(X\otimes
 u_{T})\circ\bs{r}_X^{-1}.
 \]
A  homomorphism ${\varphi}:{ T}\to { T'}$ of algebras in
$\mathcal{C}$ induces monad morphisms $\boldsymbol{\varphi_l} :
\bs{T_l} \to \bs{T_l^\prime}$ and $\boldsymbol{\varphi_r} :
\bs{T_r} \to \bs{T_r^\prime}$. For example,
$\boldsymbol{\varphi_l}X:=\varphi\otimes\mathrm{id}_X$, for any
object $X$ in $\calC$.

A  particular case of these constructions, which is very important
for our work, is obtained when we take $\calC$ to be the category
$\rmr$ of bimodules over an ordinary $\mathbb{K}$-algebra $R$ (i.e $R$ is
an algebra in the category of $\mathbb{K}$-modules, where $\mathbb{K}$ is a
commutative ring). The category $\rmr$ is monoidal with respect to
the $R$-module tensor product $\otimes_R$. Unit object is $R$. An algebra in
$\rmr$ is called an {\em $R$-ring}. $R$-rings $(T,m_T,\varphi)$ are in
bijective correspondence with $\mathbb{K}$-algebra maps $\varphi:R\to T$.
Indeed, for an algebra $(T,m_T,\varphi)$ in $\rmr$, composition
of the canonical epimorphism $T\otimes_{\mathbb K} T \to T\otimes_R T$ with
$m_T:T\otimes_R T\to T$ defines a $\mathbb{K}$-algebra structure on $T$ such
that $\varphi$ is a $\mathbb{K}$-algebra homomorphism.
Conversely, via a $\mathbb{K}$-algebra homomorphism $\varphi:R\to T$, $T$
becomes an $R$-bimodule. Multiplication of $T$ induces a morphism $m_T$ from
$T\otimes_R T$ to $T$, which makes $T$ an associative algebra in
$\rmr$.
The unit of $T$  is $\varphi$. (With a slight abuse of notation, we denote both
multiplication maps $T\ot_R T \to T$ and $T\ot_{\mathbb K} T\to T$ by the same
symbol $m_T$.)
Consequently, a $\mathbb{K}$-algebra homomorphism
$\varphi:R\to T$ defines two monads $T\otimes_R(-)$ and $(-)\otimes_R T$ on
$\rmr$.
\end{example}

Distributive laws were introduced by J. Beck \cite{Be}. They give
a way to compose two monads in order to obtain a monad.

\begin{definition}
\label{de:distributivity law} A \emph{distributive law} between
two monads $\mon{R}$ and $\mon{T} $ is a natural transformation
$\bs{t} :\bs{RT}\rightarrow \bs{TR}$ satisfying the four
conditions in Figure~\ref{fig:DLaw}.
\begin{figure}[h]
\begin{center}
{\includegraphics[scale=1]{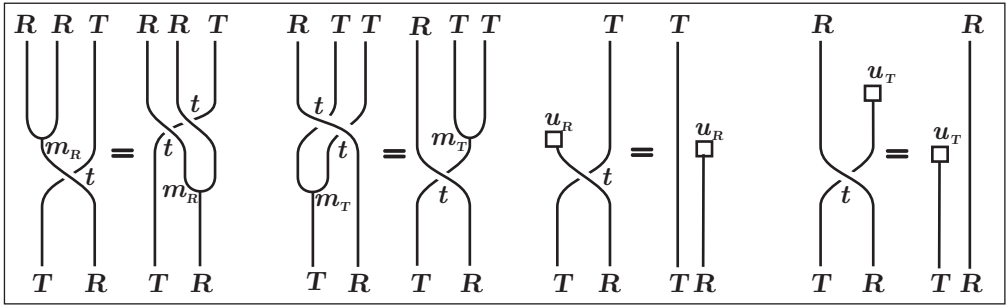}}
\end{center}
\caption{The definition of distributive laws.} \label{fig:DLaw}
\end{figure}
\end{definition}

\begin{remark}
Since we are using for the first time the diagrammatic
representation of morphisms, let us write out explicitly
the first and the third relations in Figure~\ref {fig:DLaw}. They read as
\begin{equation*}
\bs{t}\circ \bs{m_RT} =\bs{Tm_R}\circ \bs{t R}\circ
\bs{Rt},\qquad\qquad \bs{t}\circ\bs{u_RT} = \bs{Tu_R}.
\end{equation*}
\end{remark}

\begin{example}\label{ex:DL_ass.constraint}
Let $T$ be an algebra in a monoidal category as in Example \ref{ex:k_alg}.
Keeping the the notation from Example \ref{ex:k_alg}, the natural
transformation $\bs{t}:\bs{T_r}\bs{T_l}\rightarrow\bs{T_l}\bs{T_r}$, given by
\begin{equation}\label{eq:DL_ass.constraint}
 \bs{t}X:=\bs{a}_{T,X,T}
\end{equation}
for any object $X$ in $\mathcal{C}$, is a distributive law. Note
that the first equality in Figure~\ref{fig:DLaw} follows by the
Pentagon Axiom \cite[p. 282, Diagram (2.6)]{Ka}, applied to the
quadruple $(T,X,T,T)$. Similarly, by applying The Pentagon Axiom
for $(T,T,X,T)$ we deduce the second equality in the definition of
distributive laws. The fourth equality in Figure~\ref{fig:DLaw} is
a consequence of $\bs{l}_{X\otimes
T}\circ\bs{a}_{\mathbf{1},X,T}=\bs{l}_X\otimes T$, see \cite[Lemma
XI.2.2]{Ka}. The other relation in the above cited lemma can be
used to prove that the third equality in Figure~\ref{fig:DLaw}
holds too.
\end{example}

\begin{example}\label{ex:braided}
Let $\mathcal{C}$ be a braided monoidal category with
braiding $\bs{c}_{X,Y}:X\otimes Y\rightarrow Y\otimes X$. For the
definition and properties of braided monoidal categories see
\cite[Chapter XIII]{Ka}. If $R$ and $T$ are algebras in
$\mathcal{C}$ then
 \[
 \bs{t}X:=\bs{a}_{T,R,X}\circ(\bs{c}_{R,T}\otimes X)\circ \bs{a}_{R,T,X}^{-1}
 \]
defines a distributive law
$\bs{t}:\bs{R_l}\bs{T_l}\rightarrow\bs{T_l}\bs{R_l}$, where
$\bs{R_l}$ and $\bs{T_l}$ are constructed as in Example
\ref{ex:k_alg}. Obviously,
$\bs{t}^{-1}:\bs{T_l}\bs{R_l}\rightarrow\bs{R_l}\bs{T_l}$ is also
a distributive law.
\end{example}

\subsection{Admissible septuples and transposition maps. The main
  result.}\label{sec:main}
In this section we introduce admissible septuples and
transposition morphisms of them. We show that to these data one
associates functorially  para-cocyclic objects. Our aim is
twofold. On one hand, in this way we obtain a very general but at
the same time technically very simple framework.
In particular, it can be used to associate para-cocyclic objects to (co)module
algebras of bialgebroids, cf. Section \ref{sec:(co)mod.alg}.
On the other hand, the resulting setting will be easily dualized to describe
in Section \ref{sec:(co)mod_coring} the situation dual to that in Section
\ref{sec:(co)mod.alg}, i.e. the (para-)cyclic objects associated to (co)module
corings of bialgebroids.

\begin{definition}\label{cl:monad_data}
An \emph{admissible septuple}
$\calD:=(\mathcal{M},\mathcal{C},\O,\I,\bs{\Pi}, \boldsymbol{t,i})$ is
defined by the following data:
\begin{itemize}
\item Two categories $\mathcal{M}$ and $\mathcal{C}$;
 \item Two monads $\O$ and $\I$ on $\mathcal{M}$;
 \item A functor $\bs{\Pi}:\mathcal{M}\to \mathcal{C}$; \item A
distributive law $\boldsymbol{t}:\boldsymbol{T_rT_l}\to
\boldsymbol{T_lT_r}$; \item A natural transformation
$\boldsymbol{i}:\bs{\Pi} \boldsymbol{{T_l}}\to \bs{\Pi}\, \boldsymbol{T_r}$.
\end{itemize}
These data are assumed to satisfy the relations
\begin{equation}
\boldsymbol{i}\, \circ \bs{\Pi}\, \boldsymbol{u_{T_l}}\, = \bs{\Pi}\,
\boldsymbol{u_{T_r}}
\qquad \textrm{and}\qquad \boldsymbol{i} \circ
\bs{\Pi}\boldsymbol{m_{T_l}} = \bs{\Pi} \boldsymbol{m_{T_r}}
\circ\boldsymbol{i} \I{}\ \circ \bs{\Pi}{\boldsymbol{t}} \circ
\boldsymbol{i} \O{}. \label{eq:i}
\end{equation}
\end{definition}

Examples of admissible septuples will be given in Section \ref{sec:ex_adm_sep},
where also several applications of the main result of this section, Theorem
\ref{thm:coc_distr_law}, will be indicated.
\medskip

By \cite[p. 281]{We}, to every monad $\O: \mathcal{M}\rightarrow \mathcal{M}$
and object $X$ in $\calM$ one can associate a cosimplicial object
of components $\O^{n+1} \,X$
in $\mathcal{M}$. Thus in particular an admissible septuple $\calD$ in
Definition \ref{cl:monad_data} determines a cosimplicial object
in $\mathcal{M}$. It can be transported to $\mathcal{C}$ via the functor
$\bs{\Pi}:\mathcal{M}\rightarrow\mathcal{C}$ in Definition
\ref{cl:monad_data}. The resulting cosimplex in $\mathcal{C}$ will
be  denoted by $Z^{\ast }(\calD,X)$. By construction,
$Z^{n}(\calD,X)=\bs{\Pi} \O^{n+1}X$ and, for every $k\in \{0,\dots
,n\}$, the coface maps $ {d}_{k}:\bs{\Pi}\O^{n}X\rightarrow \bs{\Pi}
\O^{n+1}X$ and the codegeneracy maps ${s}_{k}:\bs{\Pi} \O
^{n+2}X\rightarrow \bs{\Pi} \O^{n+1}X$ are given by
\begin{equation}\label{eq:coface.codegeneracy}
{d}_{k}:=\bs{\Pi} \O^{k}\boldsymbol{u_{T_l}}\O^{n-k} X, \qquad\qquad
{s}_{k}:=\bs{\Pi}\O^{k}\boldsymbol{m_{T_l}}\O^{n-k}X.
\end{equation}

Our aim is to construct a category $\calW_\calS$ such that
$Z^{\ast }(\calD,-)$ can be regarded as a functor from
$\calW_\calS$ to the category of para-cocyclic objects in $\calC$.
Observe that, for an admissible septuple $\calD$ in Definition
\ref{cl:monad_data}, the distributive law $\boldsymbol{t}$ is
lifted to a natural transformation
$\boldsymbol{t}_n:\bs{\Pi}\I\O^n\rightarrow \bs{\Pi}\O^n\I$,
\begin{equation}\label{eq:t_n}
{\boldsymbol{t}_n:=\bs{\Pi} \O^{n-1}{\boldsymbol{t}}\circ
\bs{\Pi} \O^{n-2}{\boldsymbol{t}}\O \circ \dots \circ \bs{\Pi}
\bs{T_l}{\boldsymbol{t}}\bs{T_l}^{n-2}  \circ
\bs{\Pi}{\boldsymbol{t}}\bs{T_l}^{n-1}.}
\end{equation}

\begin{definition}
Let $\calD:=(\mathcal{M},\mathcal{C},\O,\I,\bs{\Pi}, \boldsymbol{t,i})$
be an admissible septuple. We say that an arrow
$w:\bs{T_r}X\to\bs{T_l}X$ in $\calM$ is a \emph{transposition
morphism} with respect to $\calD$ if
\begin{equation}
w\circ \boldsymbol{u_{T_r}} X =\boldsymbol{u_{T_l}} X\qquad
\textrm{and}\qquad w\circ \boldsymbol{m_{T_r}} X =
\boldsymbol{m_{T_l}} X \circ \O{} w \circ {\boldsymbol{t}} X \circ
\I{} w \label{eq:w}.
\end{equation}
The category of pairs $(X,w)$, with $w: \I X \to \O X$ a
transposition morphism of $\calD$, will be denoted by $\calW_\calD$.
A morphism from $(X,w)$ to $(X',w')$ is an arrow $f:X\rightarrow
X'$ in $\calM$ such that $\bs{T_l}f\circ w=w'\circ\bs{T_r}f$.
\end{definition}

Morphisms $w:\bs{T_r}X\to\bs{T_l}X$ satisfying (\ref{eq:w}), for a
distributive law $\boldsymbol{t}:\boldsymbol{T_rT_l}\to \boldsymbol{T_lT_r}$,
were termed \emph{$\boldsymbol{t}$-algebras} in \cite{Bur:D-alg}. Based on
\cite[Proposition I.1.1]{Bur:D-alg}, transposition morphisms can be
characterized as in Proposition \ref{prop:w_char} below. Recall that a module
of a monad $\bs{(T,m_T,u_T)}$ on a category $\mathcal{M}$ is a pair
$(Y,\varrho)$,
consisting of an object $Y$ and a morphism $\varrho:\bs{T} Y \to Y$ in
$\mathcal{M}$, such that $\varrho \circ \bs{T} \varrho =\varrho \circ \bs{m_T}
Y$
and $\varrho \circ \bs{u_T} Y =\mathrm{Id}_Y$ {(i.e. $\varrho$ is
  associative and unital)}. A morphism of $\bs{T}$-modules
$(Y,\varrho)\to (Y',\varrho')$ is a morphism $f:Y\to Y'$ in $\mathcal{M}$,
such that
$f\circ \varrho= \varrho'\circ \bs{T}f$.

\begin{proposition}\label{prop:w_char}
Consider an admissible septuple
$\calD:=(\mathcal{M},\mathcal{C},\O,\I,\bs{\Pi}, \boldsymbol{t,i})$. There is
a bijective correspondence between objects $(X,w)$ in the category
$\calW_\calD$ and $\I$-modules of the form $(\O X, \varrho)$, satisfying
\begin{equation}\label{eq:bim_cond}
\bs{m}_\O X \circ \O \varrho \circ \boldsymbol{t} \O X = \varrho
\circ \I \bs{m}_\O X.
\end{equation}
Moreover, a morphism $f:X \to X'$ in $\mathcal{M}$ is a morphism in
$\calW_\calD$ if and only if $\O f$ is a $\I$-module morphism.
\end{proposition}

\begin{proof}
{Similarly to the proof of \cite[Proposition I.1.1]{Bur:D-alg} one checks
  that, for an object $(X,w)$ in $\calW_\calD$, an associative and unital
  $\I$-action on $\O X$ satisfying (\ref{eq:bim_cond}) is given by
  $\varrho_w:= \bs{m}_\O X \circ \O w \circ \boldsymbol{t} X: \I\O X \to \O
  X$.}
Conversely, note that for a $\I$-module
$(\O X, \varrho)$, (\ref{eq:bim_cond}) is equivalent to $\varrho= \bs{m}_\O X
\circ \O \varrho\circ \O \I \bs{u}_\O X \circ \boldsymbol{t} X$. With this
identity at hand, the pair $(X,w_\varrho:=\varrho\circ \I \bs{u}_\O X)$ is
checked to be an object in $\calW_\calD$. A straightforward computation shows
that the two constructions are mutual inverses. A morphism $\O f$ is a
morphism of $\I$-modules $(\O X,\varrho_w)\to (\O X',\varrho_{w'})$ if
\begin{equation}\label{eq:f_bim_map}
\bs{m}_{\O} X'\circ \O w'\circ \O \I f \circ \boldsymbol{t} X =
\bs{m}_{\O} X'\circ \O \O f \circ \O w \circ \boldsymbol{t} X.
\end{equation}
If $f$ is a morphism in $\calW_\calD$ then (\ref{eq:f_bim_map}) obviously
holds. In order to prove the converse implication, compose both sides of
(\ref{eq:f_bim_map}) with $\I \bs{u}_\O X$ on the right.
\end{proof}

\begin{theorem}\label{thm:coc_distr_law}
Consider an admissible septuple $\calD$ and a transposition map
$w:\bs{T_r}X\to\bs{T_l}X$ in $\calW_\calD$. The cosimplicial
object $Z^\ast(\calD,X)$ is para-cocyclic with respect to
\begin{equation}\label{eq:para_coc_op_distr_law}
{w}_n:=\bs{\Pi} \O^n w \circ \boldsymbol{t}_n X\circ \boldsymbol{i}\O^n
X.
\end{equation}
We shall denote this para-cocyclic object by $Z^\ast(\calD,w)$.
For a morphism $f:(X,w)\to (X',w')$ in $\calW_\calD$, the morphisms $\bs{\Pi}
\O^{n+1} f:Z^n (\calD,w) \to Z^n(\calD,w')$ determine a morphism of
para-cocyclic objects.
\end{theorem}
\begin{proof} In Figure~\ref{main_thm1} we show that the
  morphism \eqref{eq:para_coc_op_distr_law} is
compatible with the coface maps, that is
\begin{equation}
{w}_{n}\circ {d}_{0}={d}_{n}\qquad \text{and\qquad } {w}_{n}\circ
{d}_{k}={d}_{k-1}\circ {w} _{n-1} \label{ec:cocyc_1}
\end{equation}
{for any $k\in\{1,\ldots, n\}$}. The proof of the first
equation is given in three steps in the left picture. To simplify
the diagrams, we draw the $n$ strings representing $\bs{{T_l}}^n$
as a black stripe. For the first equality  we used the
compatibility between $\bs{i}$ and the unit of $\bs{{T_l}}$, that
is the first equation in (\ref{eq:i}). Next we applied $n$ times
the compatibility relation between the distributive law $\bs{t}$
and the unit of $\bs{{T_r}}$, i.e. the third equality in
Figure~\ref{fig:DLaw}. The first relation in (\ref{eq:w}) implies
the third equality.
\begin{figure}[h]
\begin{center}
{\includegraphics[scale=1]{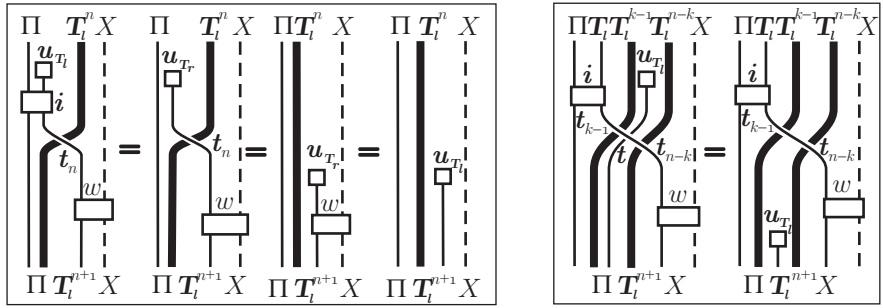}}
\end{center}
\caption{The proof of the relations (\ref{ec:cocyc_1}).}
\label{main_thm1}
\end{figure}
The second relation in (\ref{ec:cocyc_1}) follows in a similar
way, as it is shown in the right picture in
Figure~\ref{main_thm1}. Note that the leftmost black stripe
represents $\bs{T_l}^{k-1}$ and the other one represents
$\bs{T_l}^{n-k}$. Since $\bs{u_{T_l}}$ is a natural
transformation, the box representing it can be pushed down along
the string until it meets the crossing $\bs{t}$. By the fourth
identity in Figure~\ref{fig:DLaw}, one can push $\bs{u_{T_l}}$
under the string in the crossing. To conclude the proof of this
equality, we use once again that $\bs{u_{T_l}}$ is a natural
transformation to move it to the bottom of the diagram.

Next we prove that the morphism \eqref{eq:para_coc_op_distr_law}
and the codegeneracy maps are compatible too, that is
 \begin{equation}\label{ec:cocyc_2}
{w}_{n}\circ {s}_{0}
 ={s}_{n}\circ
 ({w}_{n+1})^{2}\qquad\text{and}\qquad{w}_{n}\circ
{s}_{k}={s}_{k-1}\circ {w}_{n+1}
 \end{equation}
{for any $k\in\{1,\ldots, n\}$}. The proof of the first
relation can be found in the left picture in
Figure~\ref{main_thm2}. As before, the black stripe represents
$\bs{{T_l}}^n$. The morphisms corresponding to the first two
diagrams are equal in view of the second equation in (\ref{eq:i}).
By applying $n$ times the first identity in Figure~\ref{fig:DLaw},
it follows that the second and the third diagrams represent the
same morphism. Taking into account the second relation in
(\ref{eq:w}) we got the penultimate equality, while for the last
one we used that ${\boldsymbol{i}}$ is a natural transformation.
\begin{figure}[h]
\begin{center}
{\includegraphics[scale=1]{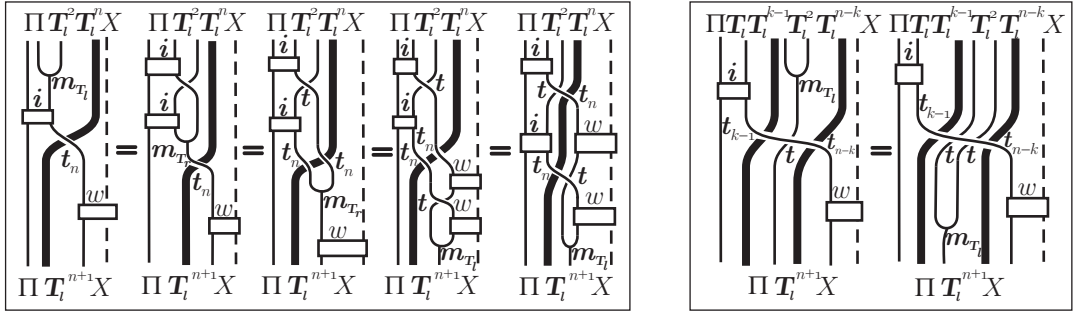}}
\end{center}
\caption{The proof of the relations (\ref{ec:cocyc_2}).}
\label{main_thm2}
\end{figure}

The other relation in (\ref{ec:cocyc_2}) immediately follows by
the second identity in Figure~\ref{fig:DLaw} and the fact
$\bs{m_{T_l}}$ is natural (see the second picture in
Figure~\ref{main_thm2}).

Since the coface and codegeneracy morphisms (\ref{eq:coface.codegeneracy}) are
defined in terms of natural transformations, the morphisms $\bs{\Pi}
\O^{n+1} f:Z^n (\calD,X) \to Z^n(\calD,X')$ determine a morphism
$Z^\ast (\calD,f):Z^\ast (\calD,X)\to Z^\ast (\calD,X')$ of
cosimplicial objects, for any morphism $f:X\to X'$ in
$\mathcal{M}$. It follows immediately from the definition of a morphism in
$\calW_\calD$ that if $f:(X,w)\to (X',w')$ is a morphism in $\calW_\calD$ then
$Z^\ast (\calD,f)$ is a morphism of para-cocyclic objects $Z^\ast (\calD,w)\to
Z^\ast (\calD,w')$.
\end{proof}

\begin{corollary}\label{cor:cocyc_quotient}
Let $\calD$ be an admissible septuple as in Definition
\ref{cl:monad_data} and let $w:\bs{T_r}X\to\bs{T_l}X$ be a
transposition morphism in $\calW_\calS$.
Consider the corresponding para-cocyclic morphism ${w}_n$ in
(\ref{eq:para_coc_op_distr_law}). If the coequalizer
\[
\xymatrix{ Z^n(\calD,w) \ar@<.5mm>[rr]^{({w}_n)^{n+1}}
\ar@<-.5mm>[rr]_{\mathrm{Id}_{Z^n(\calD,w)}}&& Z^n(\calD,w)
\ar[r]& \widehat{Z}^n(\calD,w)}
\]
exists in $\mathcal{C}$, for every non-negative integer $n$, then it defines a
cocyclic cosimplex $\widehat{Z}^\ast(\calD,w)$.
\end{corollary}
\begin{proof}
{Let $k\in\{0,\ldots, n\}$}. It follows by (\ref{ec:cocyc_1})
that ${d}_k$ satisfies ${d}_k\circ ({w}_n)^{n+1} = ({w}_{n+1})^{n+2}
\circ {d}_k$. Similarly, by (\ref{ec:cocyc_2}), the codegeneracy
morphism ${s}_k$ satisfies ${s}_k\circ ({w}_n)^{n+1} =
({w}_{n-1})^{n} \circ {s}_{k}$. Hence ${d}_k$ and ${s}_k$
determine coface morphisms $\widehat{d}_k$ and codegeneracy
morphisms $\widehat{s}_k$ on $\widehat{Z}^n(\calD,X)$. Together
with the projection $\widehat{w}_n$ of ${w}_n$ onto
$\widehat{Z}^n(\mathcal{S},X)$ they define a cocyclic object
$(\widehat{Z}^{\,\ast}(\calD,w),\widehat{d}^{\,\ast},
\widehat{s}^{\,\ast},\widehat{w}^{\,\ast})$.
\end{proof}

\subsection{Examples of admissible septuples and transposition
  maps. Applications.} \label{sec:ex_adm_sep}
In this section we shall apply Theorem \ref{thm:coc_distr_law}
to several examples of admissible septuples. In this way we shall
show that the most known (co)cyclic objects in the literature can
be obtained as direct applications of the result obtained in Section
\ref{sec:main}.
\medskip

A functorial construction of a para-cyclic object in a category of
endofunctors, of a somewhat similar flavour to that in Theorem
\ref{thm:coc_distr_law}, was proposed in \cite{Sko:CycCom}. The following
example is its dual version.

\begin{example}\label{ex:dual_Skoda}
Let $\bs{T}=(\bs{T,m,u})$ be a monad on a category $\mathcal{M}$ and
  $\bs{t:TT\to
  TT}$ be a distributive law. Assume that $\bs{t}$ satisfies the Yang-Baxter
  relation
\begin{equation}\label{eq:YB}
\bs {tT}\circ \bs {Tt} \circ \bs {tT} = \bs {Tt} \circ \bs {tT} \circ \bs {Tt}
\end{equation}
of natural transformations $\bs{TTT \to TTT}$, and $\bs{m}\circ \bs {t}\circ
\bs{t}=\bs{m}$.
As a consequence of (\ref{eq:YB}), also $\bs{T}^0=(\bs{T,m}\circ\bs{t},\bs{u})$
is a monad, and $\bs{t}$ can be regarded as a distributive law $\bs{T^0 T\to
  TT^0}$. Furthermore, the datum $\calD:=(\mathcal{M},\mathcal{M}, \bs{T,T^0},
\mathrm{Id}_{\mathcal{M}}, \bs{t},  \mathrm{Id}_{\bs{T}})$ (where
$\mathrm{Id}_{\mathcal{M}}$ denotes the identity functor ${\mathcal{M}}\to
{\mathcal{M}}$ and $\mathrm{Id}_{\bs{T}}$ is the identity natural
transformation $\bs{T}\to \bs{T}$) is an admissible septuple. For any object
$X$ in $\mathcal{M}$, the identity morphism $\mathrm{Id}_{\bs{T}X}$ is a
transposition morphism. The corresponding para-cocyclic morphism is $\bs{t}_n$
in (\ref{eq:t_n}).
\end{example}

The simplest example of an admissible septuple can be obtained by
starting with a morphism $\varphi :R\rightarrow T$ of
$\mathbb{K}$-algebras. As in
Example \ref{ex:k_alg}, we define two monads on the category $\mathcal{M}
:=R$-$\mathrm{Mod}$-$R\ $by $\boldsymbol{T_l}:=T\otimes _{R}(-)$ and $
\boldsymbol{T_r}:=(-)\otimes _{R}T.$ The category $\calC$ is, by
definition, the category $\text{Mod-}\mathbb{K}$ of $\mathbb{K}$-modules. The
functor $\bs{\Pi}$ is constructed below.

\begin{definition}\label{def:Pi}
On the objects $X\in \rmr$, the functor $\bs{\Pi}:\rmr\to
\text{Mod-}\mathbb{K}$ is defined as a coequalizer
\begin{equation*}
\xymatrix{
X \otimes_{\mathbb K} R\ar@<2pt>[rr]^{x\ot r \mapsto x\cdot r}
\ar@<-2pt>[rr]_{x\ot r \mapsto r\cdot x}&&X
\ar[rr]^{\bs{p}X}&&\bs{\Pi}X
}.
\end{equation*}
For a morphism $f:X\rightarrow Y$ of $R$-bimodules, $\bs{\Pi} f$ is the
unique $\mathbb{K}$-linear map
such that $pY\circ f=\bs{\Pi} f\circ pX.$ Hence $p$ can be interpreted as a
natural epimorphism from the forgetful functor $U:$ $R$-$\mathrm{Mod}$-$R\to
\text{Mod-}\mathbb{K}$ to $\bs{\Pi}$.
\end{definition}

\begin{remark}\label{re:cyclic.tensor.product}
An $R$-bimodule $X$ can be considered as a left or right module for the
enveloping algebra $R^e:=R\ot_{\mathbb K} R^{op}$ of $R$. In terms of the
functor
$\bs{\Pi}$ in Definition \ref{def:Pi}, the {\em cyclic tensor product}
$X\wotR Y$ of two $R$-bimodules $X$ and $Y$ is defined by
$X\wotR Y:=\bs{\Pi}(X\ot_R Y)\cong X\ot_{R^e}Y$.
With this interpretation in mind, the $\mathbb{K}$-module $\bs{\Pi} X\cong
R\ot_{R^e} X$
can be seen as the cyclic tensor product of $R$ and $X$. For $R$-bimodules
$X_1,\dots, X_n$, the $n$-fold cyclic module tensor product is defined as
\[
X_1\wotR\dots\wotR X_n :=\bs{\Pi}(X_1\ot_R \dots \ot_R X_n) =
(X_1\ot_R \dots \ot_R X_k) \wotR (X_{k+1}\ot_R \dots \ot_R X_n),
\]
for $k\in\{1,\dots, n-1\}$.
It is generated by the cyclic tensor monomials $x_1 \wotR \dots \wotR x_n$.
It is well known that the symmetry $c_{X,Y}:X\otimes_{\mathbb K} Y\to
Y\ot_{\mathbb K} X$ induces
a natural isomorphism $\bs{i}_{X,Y}:X\wotR Y\cong Y\wotR X$.
In particular, there is a natural isomorphism
$\bs{i}_{X_1,\dots,X_n}:X_1\wotR\dots\wotR X_n\to X_2\wotR
\dots\wotR X_n\wotR X_1$, that maps a generator
$x_1\wotR\dots\wotR x_n$ to $x_2\wotR \dots\wotR x_n\wotR x_1$.
\end{remark}

As we have noticed in Example \ref{ex:DL_ass.constraint}, the
associativity constraint of the monoidal category $\rmr$ defines a
distributive law $\bs{t}:\bs{T_rT_l}\to\bs{T_lT_r}$. Thus, in this
particular case that we are investigating, $\bs{t}X$ is the
canonical isomorphism $(T\otimes_R X)\otimes_R T\cong
T\otimes_R(X\otimes_R T)$, for any $X$ in $\rmr$. Let us define
$\bs{i}X:\bs{\Pi}\bs{T_l}X\to\bs{\Pi}\bs{T_r}X$ by
$\bs{i}X:=\bs{i}_{T,X}$, as in Remark
\ref{re:cyclic.tensor.product}. In terms of these natural
transformations we can give one of the main examples of admissible
septuples.

\begin{proposition}\label{ex:bimod_Pi} \label{ex:bimod_i}\label{ex:alg_sep}
Let $\varphi:R\to T$ be a morphism of $\mathbb{K}$-algebras. The following
data:
\begin{itemize}
 \item the categories $\calM:=\rmr$ and $\calC:=\rm{Mod}\text{-}\mathbb{K}$,
 \item the monads $\bs{T_l}:=T\ot_R(-)$ and $\bs{T_r}:=(-)\ot_R T$,
 \item the functor $\bs{\Pi}:\calM\to\calC$, $\bs{\Pi} X:=R\ot_{R^e}X$,
 \item the natural transformation $\bs{tX}:(T\otimes_R X)\otimes_R T\to
   T\otimes_R(X\otimes_R T)$ defined by the canonical isomorphism;
 \item the natural transformation $\bs{i}X:T\wotR \, X\to X\wotR \, T$,
 $t\wotR\, x\mapsto x\wotR\, t$,
 \end{itemize}
define an admissible septuple $\calD_T$.
\end{proposition}

\begin{proof} Let $X$ be an $R$-bimodule. By definition,
 \[\bs{\Pi}\bs{T_l}^nX= R\ot_{R^e} (T^{\ot_R\, n} \ot_R X)
\cong T^{\wotR n}\wotR X\qquad
 \text{and}\qquad\bs{\Pi}\bs{T_r}^nX\cong  X\wotR T^{\wotR n}.\]
Via these identifications, $\bs{iT_l}X=\bs{i}_{T,T,X}$ and
$\bs{iT_r}X=\bs{i}_{T,X,T}$. So the conditions (\ref{eq:i}) take
the form
\[
\bs{i}_{T,X} \circ (\varphi \wotR X)= X\wotR \varphi
\qquad\textrm{and}\qquad \bs{i}_{T,X}  \circ (m_T \wotR X)=(X\wotR
m_T)\circ \bs{i}_{T\ot_R T,X}, \] identities which are obvious.
\end{proof}

Let $\calD_T$ be the  admissible septuple associated to an algebra
morphism $\varphi:R\to T$. A morphism of $R$-bimodules
$w:X\otimes_{ R} { T}\to { T}\otimes_{ R} X$ is a transposition map
in $\calW_{\calD_T}$ if, and only if, it satisfies the conditions
\begin{equation}\label{eq:bim_tr_map}
w\circ (X\otimes_{ R} {\varphi})={\varphi}\otimes_{ R} X\qquad
\textrm{and}\qquad w\circ (X\otimes_{ R} m_{ T})= ( m_{ T}
\otimes_{ R} X)\circ ({ T} \otimes_{ R} w)\circ (w\otimes_{ R} {
T}),
\end{equation}
where $m_{ T}:{ T} \otimes_{ R} { T}\to { T}$ denotes the
multiplication map $t\otimes_{ R} t'\mapsto tt'$. Note in passing
the similarity of conditions (\ref{eq:bim_tr_map}) to some of
those defining an entwining structure over ${ R}$.
({For the definition of entwining
structures see \cite[Definition 2.1]{BMa}, and for a reformulation over an
arbitrary base algebra $R$ see \cite[Section 2.3]{BB}.})
By Proposition \ref{prop:w_char}, there is a bijective correspondence between
transposition maps $w:X\otimes_{ R} { T}\to { T}\otimes_{ R} X$ on one hand,
and right $T$-actions on $T\ot_R X$, which are left $T$-module maps with
respect to the left $T$-action $t'\cdot (t\ot_R x)=t't\ot_R x$,
on the other hand.

\begin{theorem}
\label{thm:para_cyc}
Let $\calD_T$ be the  admissible septuple associated to an algebra
morphism $\varphi:R\to T$. Let $w:X\otimes_{ R} { T}\to {
T}\otimes_{ R} X$  be a transposition map in $\calW_{\calD_T}$, that
is, a morphism of $R$-bimodules satisfying (\ref{eq:bim_tr_map}).
Then there is a cocylic quotient $\widehat{Z}^\ast(\calD_T,w)$ of $T^{\wotR
  {\ast+1}}\wotR X$ such
that its cocyclic structure is induced by the para-cocyclic
morphisms ${w}_{n}:T^{\wotR {n+1}}\wotR X\rightarrow T^{\wotR
{n+1}}\wotR X$,
\begin{equation}\label{ec:t c} {w}_{n}:=\left({T}^{\wotR{n}}\wotR
    w\right)\circ \bs{i}_{T,\dots,T,X},
\end{equation}
where $\bs{i}_{T,\dots,T,X}$ is the $\mathbb{K}$-linear map defined in
Remark \ref{re:cyclic.tensor.product}.
\end{theorem}

\begin{proof}
Apply Theorem \ref{thm:coc_distr_law} for $\calD=\calD_T$. It yields
a para-cocyclic object ${Z}^\ast(\calD_T,w)$ whose para-cocyclic
operator is given in formula (\ref{eq:para_coc_op_distr_law}).
For $\calD_T$, the map $\bs{t}_n$ is the identity morphism of
$T^{\wotR {n+1}}\wotR X$, cf. (\ref{eq:t_n}). Hence $w_n$
satisfies (\ref{ec:t c}). We conclude the proof by applying Corollary
\ref{cor:cocyc_quotient}.
\end{proof}

\begin{corollary}\label{ex:no_coef}
Let $\mathcal{S}_T$ be the admissible septuple associated to a
$\mathbb{K}$-algebra homomorphism ${\varphi}: R\to T$ as in
Proposition \ref{ex:bimod_Pi}. Then the canonical isomorphism
$w_T:R\ot_R T\to T\ot_R R$ is a transposition map in
$\calW_{\calD_T}$. The corresponding cocyclic cosimplex $Z^\ast
(\mathcal{S}_T,w_T)$ has in degree $n$ the $\mathbb{K}$-module
$Z^n (\mathcal{S}_T,w_T)={{
  T}^{\widehat{\otimes}_{ R}n+1}}$. The coface and codegeneracy maps are
\begin{align*}
d_k(t_0 \widehat{\otimes}_{ R}\, t_1 \widehat{\otimes}_{ R} \dots
\widehat{\otimes}_{ R}\, t_{n-1})&= t_0 \widehat{\otimes}_{ R}\, t_1
\widehat{\otimes}_{ R} \dots \widehat{\otimes}_{ R}\, t_{k-1}
\widehat{\otimes}_{ R}\, 1_{ T} \widehat{\otimes}_{ R}\, t_{k}
\widehat{\otimes}_{ R} \dots
\widehat{\otimes}_{ R}\, t_{n-1}\\
s_k(t_0 \widehat{\otimes}_{ R}\, t_1 \widehat{\otimes}_{ R} \dots
\widehat{\otimes}_{ R}\, t_{n+1}) &= t_0 \widehat{\otimes}_{ R}\,
t_1 \widehat{\otimes}_{ R} \dots \widehat{\otimes}_{ R}\, t_{k-1}
\widehat{\otimes}_{ R}\, t_{k} t_{k+1} \widehat{\otimes}_{ R}\,
t_{k+2} \widehat{\otimes}_{ R} \dots \widehat{\otimes}_{ R}\,
t_{n+1},
\end{align*}
where {$k\in\{0,\ldots,n\}$}. The cocyclic operator is given
by
\[
w_n(t_0 \widehat{\otimes}_{ R} \, t_1 \widehat{\otimes}_{ R} \dots
\widehat{\otimes}_{ R}\, t_n)=t_1 \widehat{\otimes}_{ R}\, t_2
\widehat{\otimes}_{ R} \dots \widehat{\otimes}_{ R}\, t_n
\widehat{\otimes}_{ R}\, t_0.
\]
\end{corollary}

\begin{remark}  As $u_T:k\to T$, the unit of a $\mathbb{K}$-algebra
${ T}$, is an algebra map, we can apply Proposition
\ref{ex:bimod_Pi} to get an admissible septuple $\calD_{u_T}$.
Corresponding transposition maps were also considered by Kaygun in
\cite{Kay:UniHCyc} to construct cocyclic $\mathbb{K}$-modules.
{His approach should be considered, however, dual to our one (see
related remarks in the Introduction).} 
\end{remark}

It follows by an observation in \cite[page 11]{Bur:D-alg} that for the
admissible septuple $\calD_T$, associated to an algebra morphism $\varphi:R\to
T$ in Proposition \ref{ex:alg_sep}, any $R$-$T$ bimodule $Y$ admits a
transposition morphism  $w_Y:Y\ot_R T \to T\ot_R Y$, $y\ot_R t \mapsto 1_T
\ot_R y\cdot t$. In particular, for any $R$-bimodule  $X$, the pair $(X\ot_R
T, (u_T \ot_R X \ot_R T)\circ( X\ot_R m_T))$ is an object in
$\calW_{\calD_T}$. Corresponding para-cocyclic objects are given in the
following

\begin{example}
Let $\mathcal{S}_T$ be the admissible septuple associated to a
$\mathbb{K}$-algebra homomorphism ${\varphi}: R\to T$ as in
Proposition \ref{ex:bimod_Pi}. For any $R$-$T$ bimodule $Y$, there
is a para-cocyclic cosimplex $Z^\ast (\mathcal{S}_T,w_Y)$, given
in degree $n$ by the $\mathbb{K}$-module $Z^n
(\mathcal{S}_T,w_Y)={
{T}^{\widehat{\otimes}_{R}n+1}\widehat{\otimes}_{R} Y}$. The
coface and codegeneracy maps are
\begin{align*}
d_k(t_0 \widehat{\otimes}_{ R}\, t_1 \widehat{\otimes}_{ R} \dots
\widehat{\otimes}_{ R}\, t_{n-1}\widehat{\otimes}_{ R} \,y)&= t_0
\widehat{\otimes}_{ R}\, t_1 \widehat{\otimes}_{ R} \dots \widehat{\otimes}_{
  R}\, t_{k-1} \widehat{\otimes}_{ R}\, 1_{ T} \widehat{\otimes}_{ R}\, t_{k}
\widehat{\otimes}_{ R} \dots \widehat{\otimes}_{ R}\,
t_{n-1}\widehat{\otimes}_{ R} \,y\\
s_k(t_0 \widehat{\otimes}_{ R}\, t_1 \widehat{\otimes}_{ R} \dots
\widehat{\otimes}_{ R}\, t_{n+1}\widehat{\otimes}_{ R}\,y)&= t_0
\widehat{\otimes}_{ R}\, t_1 \widehat{\otimes}_{ R} \dots \widehat{\otimes}_{
  R}\, t_{k-1} \widehat{\otimes}_{ R}\, t_{k} t_{k+1} \widehat{\otimes}_{ R}\,
t_{k+2} \widehat{\otimes}_{ R} \dots \widehat{\otimes}_{ R}\,
t_{n+1}\widehat{\otimes}_{ R} \,y,
\end{align*}
where {$k\in\{0,\ldots,n\}$}. The para-cocyclic operator is
given by
\[
w_n(t_0 \widehat{\otimes}_{ R} \, t_1 \widehat{\otimes}_{ R} \dots
\widehat{\otimes}_{ R}\, t_n\widehat{\otimes}_{ R} \,y)=t_1 \widehat{\otimes}_{
  R} \dots \widehat{\otimes}_{ R}\, t_n\widehat{\otimes}_{ R} 1_T
\widehat{\otimes}_{ R} \, y \cdot t_0.
\]
Note that $w_n$ is degenerate in the sense that the cocyclic quotient of
$Z^\ast (\mathcal{S}_T,w_Y)$ (cf. Corollary \ref{cor:cocyc_quotient}) is given
by $\widehat{Z}^n (\mathcal{S}_T,w_Y)=R\,\widehat{\ot}_R\, Y$, in every degree
$n$.
\end{example}

Next we are going to associate an admissible septuple to every ribbon
algebra. Recall that a \emph{ribbon algebra} is an algebra
$(T,m_T,u_T)$ in a braided monoidal category $\calM$ together with
an automorphism $\sigma:T\to T$ in $\calM$ such that
\begin{equation}\label{eq:ribbon}
\sigma \circ u_T=u_T\qquad\text{and}\qquad \sigma \circ
m_T=m_T\circ (\sigma\otimes\sigma)\circ \bs{c}_{T,T}^2.
\end{equation}
Ribbon algebras appeared in \cite{AM}, where they are used to
define cyclic homology of quasialgebras (non-associative algebras
that are obtained by a cochain twist). We shall show that the
ribbon automorphism $\sigma$ can be used to define a certain
admissible septuple. For, we start with an algebra $(T,m_T,u_T)$
and an automorphism $\sigma:T\to T$ in a braided monoidal
category. It is easy to see that $T$ is also an associative and
unital algebra with respect to  {$m'_T:=m_T\circ\bs{c}_{T,T}$}
and $u'_T:=u_T$. To make distinction between $T$ and the new
algebra, latter one will be  denoted by $T'$.

Consider the monads $\bs{T_l'}$ and $\bs{T_l}$ on $\calM$, defined
as in Example \ref{ex:k_alg}. In the following, $\bs{T'_l}$ will
play the role of $\bs{T_r}$ in the definition of an admissible
septuple. We have seen in Example \ref{ex:braided} that {
$\bs{c}_{T,T}:T\ot T\to T\ot T$} induces a distributive law $\bs{
t}:\bs{T_l'T_l}\to\bs{T_lT_l'}$
\[\bs{ t}X:=\bs{a}_{T,T,X}\circ({\bs{c}_{T,T}}\ot
 X)\circ\bs{a}_{T,T,X}^{-1},
 \]
where $X$ is an arbitrary object in $\calM$. Furthermore, we take
$\calC=\calM$ and $\bs{\Pi}=\mathrm{Id}_\calM$.
By definition, the natural transformation $\bs{i}:\bs{T_l}\to\bs{T_l'}$ is
\[\bs{ i}X:=({\sigma}\ot
X).
 \]
It is not difficult to prove that the relations in (\ref{eq:i})
hold if, and only if, the identities in (\ref{eq:ribbon}) are
satisfied. Thus, we have the following

\begin{proposition}\label{pr:ribbon}
The algebra $(T,m_T,u_T)$ in a {braided} monoidal category
$\calM$ is a ribbon algebra with ribbon automorphism $\sigma$ if
and only if
$\calD_{T,\sigma}:=(\calM,\calM,\bs{T_l},\bs{T_l'},\mathrm{Id}_\calM,\bs{t},
\bs{i})$, the septuple constructed above, is admissible.
\end{proposition}

Let $(T,m_T,u_T)$ be a ribbon algebra with ribbon automorphism
$\sigma$. In view of Proposition \ref{pr:ribbon},  we can speak about
transposition morphisms with respect to $\calD_{T,\sigma}$.
A morphism $w:T\ot X\to T\ot X$ in $\calM$ is a transposition map in
$\calW_{\calD_{T,\sigma}}$ if, and only if
\begin{align}
w&\circ (u_{T}\otimes X)=u_{T}\otimes X,\label{eq:transition1.ribbon}\\
 {w}&{\circ (m_{T}\otimes X)\circ (\bs{c}_{T,T}\otimes
X)=(m_{T}\otimes X)\circ (T\otimes w)\circ
 (\bs{c}_{T,T}\otimes X)\circ
(T\ot w)}\label{eq:transition2.ribbon}.
\end{align}
Note that, in the second equation, we omitted the associativity
constraints, to make the formula as short as possible. In fact, in
view of the Coherence Theorem, we can omit bracketing in any
equality involving morphisms in an arbitrary monoidal category
$\calM$. {(For MacLane's Coherence Theorem consult e.g. \cite[Theorem 1, p.
162]{ML} or \cite[pp. 420-421]{Mj}).}

Relations (\ref{eq:transition1.ribbon}) and
(\ref{eq:transition2.ribbon}) already appeared in the definition
of braided twistors, structures that are used to construct new
associative and unitary multiplications on ${ T}\otimes { T}.$ For
details the reader is referred to \cite{LPvO}.

As an application of Theorem \ref{thm:coc_distr_law} we get
Proposition \ref{prop:ribbon} below. Note that, for simplifying the formulae
of coface, codegeneracy and para-cocyclic morphisms, we omitted the
associativity and unit constraints.

\begin{proposition}\label{prop:ribbon}
Let $(T,m_T,u_T,\sigma)$ be a ribbon algebra in a braided monoidal
category $\calM$. For every object $X$ in $\calM$, the sequence
$Z^n(\calD_{T,\sigma},X):=T^{\ot n+1}\ot X$ defines a cosimplicial
object, with respect to the coface and codegeneracy morphisms
\[d_k:=T^{\ot k}\ot u_T\ot T^{\ot n-k}\ot X\qquad\text{and}\qquad
s_k:=T^{\ot k}\ot m_T\ot T^{\ot n-k}\ot X,
 \]
{where $k\in\{0,\ldots, n\}$}. Moreover, if $w:T\ot X\to T\ot
X$ is a transposition map in $\calW_{\calD_{T,\sigma}}$, then
$Z^\ast(\calD_{T,\sigma},X)$ admits a para-cocyclic structure with
respect to the operator
\begin{equation*}{
w_{n}:=(T^{\otimes n}\otimes w)\circ (T^{\otimes n-1}\otimes
\bs{c}_{T,T}\otimes X)\circ \cdots \circ (\bs{c}_{T,T}\otimes
T^{\otimes n-1}\otimes X)\circ (\sigma \otimes T^{\otimes
n}\otimes X).}
\end{equation*}
We shall denote this
para-cocyclic object by $Z^\ast(\calD_{T,\sigma},w)$.
\end{proposition}

\begin{remarks}  {(i) Every ribbon algebra $(T,\sigma)$ in a braided
    category
$(\calM,\ot,\bs{a,l,r,c,1})$ can be seen as an algebra with ribbon
element $\sigma^{-1}$ in the opposite braided category of $\calM$.
Recall that the opposite of the braided category $\calM$ is
$(\calM,\ot,\bs{a,l,r,\widetilde{c},1})$, where
\[\widetilde{c}_{X,Y}={c}_{Y,X}^{-1}.\]}

(ii) To every (para)-cocyclic object with invertible para-cocyclic
  morphism, there corresponds a (para)-cyclic
object, namely its cyclic dual. Roughly speaking, the cyclic dual
is obtained by interchanging the coface and codegeneracy morphisms
and inverting the para-cocyclic operator. The interested reader
can find the definition of the cyclic dual in \cite{KR}. The
cyclic dual of ${Z^\ast(\calD_{T,\sigma^{-1}},\mathrm{Id}_T)}$
is, modulo a sign in the formula of $w_n$, the cyclic object in
\cite[Theorem 4]{AM}. Note that, via the identification $T\ot{\bf
1}\cong T$, the identity morphism $\mathrm{Id}_T$ can be regarded
as a transposition map in {$\calW_{\calD_{T,\sigma^{-1}}}$}.
Thus, for an arbitrary $w$ in {
$\calW_{\calD_{T,\sigma^{-1}}}$}, the cyclic dual of {
$Z^\ast(\calD_{T,\sigma^{-1}},w)$} may be interpreted as a
generalization of cyclic homology introduced in \cite{AM}.
\end{remarks}

Other examples of para-cocyclic objects, obtained as applications
of Theorem \ref{thm:coc_distr_law}, will be discussed in Theorem
\ref{thm:mod_alg} and Theorem \ref{thm:comod_alg}.

\subsection{The dual construction}\label{sec:dual}
In this section we turn to studying the situation dual to that in
Section \ref{sec:main}, i.e. application of Theorem
\ref{thm:coc_distr_law} to the opposite categories
$\mathcal{C}^{op}$ and $\mathcal{M}^{op}$. By $\mathcal{M}^{op}$
we mean the category with the same classes of objects and
morphisms in $\mathcal{M}$, with composition opposite to that in
$\mathcal{M}$. Note that any diagram expressing an identity of
morphisms in $\mathcal{M}$, yields a diagram in
$\mathcal{M}^{op}$, by interchanging the top and the bottom. In
particular, a \emph{comonad} on $\mathcal{M}$ is a monad on
$\mathcal{M}^{op}$. That is, a triple
$(\O,\Delta_\O,\boldsymbol{\varepsilon}_\O)$, consisting of a
functor $\O:\mathcal{M}\to \mathcal{M}$ and natural
transformations $\Delta_\O:\O\to \O^2$ and
$\boldsymbol{\varepsilon}_\O: \O\to \boldsymbol{Id}_{\mathcal M}$.
Their compatibility axioms are obtained by reversing the arrows in
the first two diagrams in Figure \ref{fig:monad}. For two comonads
$(\O,\Delta_\O,\boldsymbol\varepsilon_\O)$ and
$(\I,\Delta_\I,\boldsymbol\varepsilon_\I)$ on a category
$\mathcal{M}$, a \emph{dual distributive law} is a distributive
law for the monads $\O$ and $\I$ on $\mathcal{M}^{op}$. That is, a
natural transformation ${\boldsymbol{t}}:\O\I\to \I\O$ such that
the relations encoded in the up-down mirror images of the diagrams
in Figure \ref{fig:DLaw}
hold.

{To dualize admissible septuples we need two comonads $\O$ and $\I$
  on a category $\mathcal{M}$, a dual distributive law
  ${\boldsymbol{t}}:\O\I\to \I\O$, a covariant
functor $\bs{\Pi}:\mathcal{M}\to \mathcal{C}$ and a natural
transformation $\boldsymbol{i}:\bs{\Pi}\, \boldsymbol{T_r}\to
\bs{\Pi} \boldsymbol{T_l}$ that satisfy the identities
\begin{equation}\label{eq:i_dual}
\bs{\Pi}\boldsymbol{\varepsilon_{T_l}}\circ \boldsymbol{i}=
\bs{\Pi}\boldsymbol{\varepsilon_{T_r}}
\qquad \textrm{and}\qquad \bs{\Pi}\boldsymbol{\Delta_{T_l}}\circ
\boldsymbol{i}  = \boldsymbol{i} \O{}\circ \bs{\Pi}{\boldsymbol{t}}
\circ\boldsymbol{i} \I{}\
 \circ \bs{\Pi} \boldsymbol{\Delta_{T_r}}.
\end{equation}
Such a {\em dual admissible septuple} $\calD^0=({\mathcal M},{\mathcal C},
\boldsymbol{T_l}, \boldsymbol{T_r},\bs{\Pi}, {\boldsymbol{t}},\boldsymbol{i})$
determines a simplicial object $Z_{\ast
}(\calD^0,X)$ in $\mathcal{C}$, which in degree $n$ is given by
$Z_{n}(\calD^0 ,X)=\bs{\Pi} \O^{n+1}X$. Its face maps $
{d}_{k}:\bs{\Pi}\O^{n+1}X\rightarrow \bs{\Pi} \O^{n}X$ and
degeneracy maps ${s}_{k}:\bs{\Pi} \O ^{n+1}X\rightarrow \bs{\Pi}
\O^{n+2}X$ are
\begin{equation*}
{d}_{k}:=\bs{\Pi} \O^{k}{\varepsilon_\O}\O^{n-k}X, \qquad\qquad
{s}_{k}:=\bs{\Pi}\O^{k}{\Delta_\O}\O^{n-k}X,
\end{equation*}
for any $k\in\{0,\ldots,n\}$. An arrow $w: \O X \to \I X$ in
$\mathcal{M}$ is a \emph{transposition morphism} with respect to
the dual admissible septuple $\calD^0$ if, and only if
\begin{equation}\label{eq:w_dual}
\boldsymbol{\varepsilon_{T_r}} X \circ
w=\boldsymbol{\varepsilon_{T_l}} X \qquad \textrm{and}\qquad
\boldsymbol{\Delta_{T_r}} X\circ w =\I{} w \circ  {\boldsymbol{t}}
X \circ \O{} w\circ\boldsymbol{\Delta_{T_l}} X.
\end{equation}
Morphisms between transpositions maps can be easily defined by
duality. The category of transposition maps with respect to
$\calS^0$ will be denoted by $\calW_{\calS^0}$.

Note that $\bs{t}$ can be lifted to a natural transformation
$\boldsymbol{t}_n:\bs{\Pi}\O^n\I\rightarrow \bs{\Pi}\I\O^n$,
\begin{equation}\label{eq:t_n0}
\boldsymbol{t}_n^0:=\bs{\Pi}{\boldsymbol{t}}\O^{n-1} \circ
\bs{\Pi} \O{\boldsymbol{t}}\O^{n-2}  \circ \dots \circ \bs{\Pi}
\O^{n-2}{\boldsymbol{t}}\O  \circ\bs{\Pi}
\O^{n-1}{\boldsymbol{t}}.
\end{equation}
Now we can state, for future references, the dual of Theorem
\ref{thm:coc_distr_law}.

\begin{theorem}\label{thm:c_distr_law}
Consider a dual admissible septuple $\calD^0$ as above and a
transposition morphism $w: \O X \to \I X$ in $\calW_{\calD^0}$.
The simplex $Z_\ast(\calD^0,X)$ is para-cyclic with para-cyclic
morphism
\begin{equation}\label{eq:para_c_op_distr_law}
{w}_n:= \boldsymbol{i}\O^n X\circ \boldsymbol{t}_n^0 X\circ \bs{\Pi}
\O^n w.
\end{equation}
We shall denote this cyclic object by $Z_\ast(\calD^0,w)$.
For a morphism $f:(X,w)\to (X',w')$ in $\calW_{\calD^0}$, the morphisms
$\bs{\Pi} \O^{n+1} f:Z_n (\calD^0,w) \to Z_n(\calD^0,w')$ determine a morphism
of para-cyclic objects.
\end{theorem}}

Dually to Example \ref{ex:dual_Skoda}, we have the following

\begin{example}
Let $\bs{T}=(\bs{T,\Delta,\varepsilon})$ be a comonad on a category
  $\mathcal{M}$ and $\bs{t:TT\to TT}$ a dual distributive law. Assume that
  $\bs{t}$ satisfies the Yang-Baxter relation  \ref{eq:YB} and $\bs {t}\circ
\bs{t}\circ \bs{\Delta}=\bs{\Delta}$. Then
  $\bs{T}^0=(\bs{T,t}\circ\bs{\Delta},\bs{\varepsilon})$ is a comonad, and
  $\bs{t}$ can be regarded as a distributive law $\bs{T^0 T\to
  TT^0}$. Furthermore, the datum $\calD:=(\mathcal{M},\mathcal{M}, \bs{T^0,T},
  \mathrm{Id}_{\mathcal{M}}, \bs{t},  \mathrm{Id}_{\bs{T}})$ is a dual
  admissible septuple. For any
  object $X$ in $\mathcal{M}$, the identity morphism $\mathrm{Id}_{\bs{T}X}$
  is a transposition morphism. The corresponding para-cyclic morphism is
  $\bs{t_n^0}$ in \eqref{eq:t_n0}.
Note that if in addition $\bs{t}$ is an invertible morphism in
$\mathcal{M}$ then its properties assumed above are equivalent to
the premises in \cite[Theorem 1]{Sko:CycCom}.
\end{example}

Let ${ R}$ be an algebra over a commutative ring $\mathbb{K}$. It was recalled
in Example \ref{ex:k_alg} that $\rmr$, the category of $R$-bimodules, is
monoidal with respect to the tensor product $\ot_R$ and unit
object $R$. By definition, an {\em $R$-coring}
$({C},{{\Delta}},{\epsilon})$ is a coalgebra in $(\rmr,\ot_R,R)$.

\begin{proposition}\label{ex:coring}
Let  $({C},{{\Delta}_C},{\epsilon_C})$ be an $R$-coring. The
following data:
\begin{itemize}
 \item the category $\calC:=\mathrm{Mod}$-$\mathbb{K}$ of $\mathbb{K}$-modules
   and the category $\calM:={ R}$-$\mathrm{Mod}$-${ R}$ of ${R}$-bimodules,
 \item the comonads $\bs{T_l}:={ C}\otimes_{ R} (-)$ and
   $\bs{T_r}:=(-)\otimes_{ R} {
C}$ on ${  R}$-$\mathrm{Mod}$-${ R}$,
 \item the functor $\bs{\Pi}:{R}$-$\mathrm{Mod}$-${
   R}\to\mathrm{Mod}$-$\mathbb{K}$, $M\mapsto R\ot_{R^e} M$ in Definition
   \ref{def:Pi},
 \item the trivial dual distributive
law ${\bs{t}X}:{C}\otimes_{ R}(X \otimes_{ R} { C})\to(
{C}\otimes_{  R} X) \otimes_{ R} { C}$,
 \item the natural morphism $\bs{i}X:X\widehat{\otimes}_{ R}
{C} \to { C} \widehat{\otimes}_{ R} X$, given by the flip map $x\wotR c
\mapsto c\wotR x$,
\end{itemize}
define a dual admissible septuple $\calD_C$.
\end{proposition}

\begin{proof}{
We have to check the identities in \eqref{eq:i_dual}. Recall that,
for any $R$-bimodule $X$, the cyclic tensor product
$\Pi(X):=R\ot_{R^e}X$ is isomorphic to the quotient of $X$ modulo
the $\mathbb{K}$-submodule $[X,R]$ generated by all commutators
$[x,r]$, where $x\in X$ and $r\in R$. Hence
$\widehat{rx}=\widehat{xr}$, where $\widehat{z}$ denotes the class
of $z$ in the quotient module, for any $z\in X$. To
prove the first relation in \eqref{eq:i_dual}, note that
\[\big(\Pi\bs{\varepsilon_{T_l}}\circ\bs{i}\big)(x\widehat{\ot}_R\,
c)=\widehat{\varepsilon_C(c)x}=
\widehat{x\varepsilon_C(c)}=\Pi\bs{\varepsilon_{T_r}}(x\widehat{\ot}_R\,
c).\] For the coproduct in the coring $C$ we use a Sweedler type notation,
namely we write $\Delta_C(c)=c_{(1)}\ot_R c_{(2)}$, with implicit summation
understood.  A straightforward computation yields
\[\big(\bs{\Pi}\boldsymbol{\Delta_{T_l}}\circ
\boldsymbol{i}\big)(x\widehat{\ot}_R\, c)  =c_{(1)}\widehat{\ot}_R
c_{(2)}\widehat{\ot}_R x = \big( \boldsymbol{i} \O{}\circ
\bs{\Pi}{\boldsymbol{t}} \circ\boldsymbol{i} \I{}\
 \circ \bs{\Pi} \boldsymbol{\Delta_{T_r}}\big)(x\widehat{\ot}_R\, c),\]
for any $x\in X$ and $c\in C$. Thus the second relation in
\eqref{eq:i_dual} is also proven.}
\end{proof}

Let $\calD_C$ be the  dual admissible septuple associated to an
$R$-coring $({C},{{\Delta}_C},{\epsilon_C})$. In this particular
case, a map of $R$-bimodules $w:C\otimes_{ R} X\to X\otimes_{ R}
C$ is a transposition map in $\calW_{\calD_C}$ if, and only if,
\begin{equation}\label{eq:transition.map.dual}
(X\otimes_{ R} {\varepsilon_C})\circ w ={\varepsilon_C}\otimes_{
R} X\qquad \textrm{and}\qquad (X\otimes_{ R} \Delta_{ C})\circ w=
(w\otimes_{ R} { C})\circ ({ C} \otimes_{ R} w)\circ ( \Delta_{ C}
\otimes_{ R} X).
\end{equation}

\begin{theorem}\label{thm:para_cyc_dual}
Let $\mathcal{S}_C$ be the dual admissible septuple associated to
an ${ R}$-coring, as in Proposition \ref{ex:coring}.  Let
$w:C\otimes_{ R} X\to X\otimes_{ R}C$ be a transposition map in
$\calW_{\calD_C}$, that is, a morphism of $R$-bimodules satisfying
(\ref{eq:transition.map.dual}). Then there is a cyclic {
subobject} $\widehat{Z}_\ast(\calD_C,w)$ of $C^{\wotR
{\ast+1}}\wotR X$ whose cyclic structure is {restriction of the}
para-cyclic morphism ${w}_{n}:C^{\wotR {n+1}}\wotR X\rightarrow
C^{\wotR {n+1}}\wotR X$,
\begin{equation*}
{w}_{n}:=\bs{i}_{C,\dots,C,X}^{-1}\circ \left({C}^{\wotR{n}}\wotR  w\right),
\end{equation*}
where $\bs{i}_{C,\dots,C,X}$ is the isomorphism constructed in Remark
\ref{re:cyclic.tensor.product}.
\end{theorem}

\begin{proof}
Proceed as in the proof of Theorem \ref{thm:para_cyc}.
\end{proof}

Dually to Corollary \ref{ex:no_coef}, the following holds.

\begin{corollary}\label{ex:no.coef.dual}
Let $\mathcal{S}_C$ be the dual admissible septuple associated to
an ${ R}$-coring, as in Proposition \ref{ex:coring}. Then the
canonical isomorphism $w_C:C\ot_R R\to R\ot_R C$ is a
transposition map in $\calW_{\calD_C}$. The corresponding cyclic
object $Z_\ast (\mathcal{S}_C,w_C)$ has in degree $n$ the
$\mathbb{K}$-module {$Z_n (\mathcal{S}_C,w_C)=
{C}^{\widehat{\otimes}_{ R}n+1}$. The face and degeneracy maps are
 \begin{align*} d_k(c_0 \widehat{\otimes}_{ R}\, c_1
\widehat{\otimes}_{ R} \,\dots \widehat{\otimes}_{ R}\, c_{n})&=
 \left\{
\begin{array}{ll}c_0\widehat{\otimes}_{ R}\, \dots
\widehat{\otimes}_{ R}\, \epsilon_C(c_{k})c_{k+1}
 \widehat{\otimes}_{ R}\, \dots
\widehat{\otimes}_{ R}\, c_{n}, & \qquad\qquad\;\text{for } 0\leq k< n,\\
c_{0} \widehat{\otimes}_{ R}\,c_1 \widehat{\otimes}_{ R}\dots
\widehat{\otimes}_{ R}\,c_{n-2}\widehat{\otimes}_{ R}\,
c_{n-1}\epsilon_C(c_{n}), & \qquad\qquad\;\text{for } k=n,
\end{array} \right.\\
 s_k(c_0 \widehat{\otimes}_{ R}\, c_1
\widehat{\otimes}_{ R}\, \dots \widehat{\otimes}_{ R} c_{n})& =
c_0 \widehat{\otimes}_{ R}\, \dots \widehat{\otimes}_{ R}\,
c_{k-1} \widehat{\otimes}_{ R}\, \Delta_C(c_k)\,
\widehat{\otimes}_{ R}\, c_{k+1} \widehat{\otimes}_{ R}\, \dots
\widehat{\otimes}_{ R}\, c_{n},\quad\text{for } 0\leq k\leq n.
\end{align*}}
The cyclic operator is defined by
\[
w_n(c_0\, \widehat{\otimes}_{ R}\, \dots \, \widehat{\otimes}_{ R}
\, c_n)=  c_n \, \widehat{\otimes}_{ R}\, c_0 \,
\widehat{\otimes}_{ R}\, \dots \, \widehat{\otimes}_{ R}\,
c_{n-1}.
\]
\end{corollary}
A symmetrical version of the construction in Corollary
\ref{ex:no.coef.dual} is described in \cite[Proposition
3.1]{Rang:CycCor}.

\begin{example}\label{ex:Sweedler}
Let ${\varphi}:{ S}\to { T}$ be a homomorphism of algebras over a
commutative ring $\mathbb{K}$. It determines the Sweedler's ${ T}$-coring
${
  T}\otimes_{ S}{ T}$, where on ${
  T}\otimes_{ S}{ T}$ we take the obvious ${ T}$-bimodule
  structure. The
coproduct $\Delta_{T\ot_S T}$ and the counit $\epsilon_{T\ot_S T}$
are respectively defined by
\[
\begin{array}{lll}
\Delta_{T\ot_S T}:&{ T}\otimes_{ S}{ T} \to ({ T}\otimes_{ S}{
T})\otimes_{
  T} ({ T}\otimes_{ S}{ T}), \quad
&\Delta_{T\ot_S T}(t\otimes_{ S} t')= t\otimes_{ S}1_{ T} \otimes_{
S} t'\ ,\\
\,\epsilon\,_{T\ot_S T}:& { T}\otimes_{ S}{ T} \to { T},\quad
&\,\epsilon\,_{T\ot_S T}\,(t\otimes_{ S} t')=tt',
\end{array}
\]
where in the definition of $\Delta_{T\ot_R T}$ we identified
$({ T}\otimes_{ S}{ T})\otimes_T({ T}\otimes_{ S}{ T})$  and ${
T}\otimes_{ S}{T}\otimes_{ S}{ T}$.
\end{example}

{
For $S$-bimodules $T$ and $X$, let $v:T\otimes _{S}X\rightarrow X\otimes
_{S}T$ be an $S$-bimodule map. For $t\in T$ and $x\in X$ we shall use the
notation  $v(t\otimes _{S}x)=x_{v}\otimes _{S}t_{v},$ where in the
right hand side implicit summation is understood.

\begin{corollary}
\label{ex:dual_no_coef} Let ${\varphi }:{S}\rightarrow {T}$ be a
homomorphism of algebras over a commutative ring $\mathbb{K}$, $X$ be an $S$
-bimodule and $v:T\otimes _{S}X\rightarrow X\otimes _{S}T$ be an $S$
-bimodule map satisfying
\begin{equation*}
v\circ ({\varphi }\otimes _{S}X)=X \otimes _{S}\varphi\qquad \text{and}
\qquad v\circ (m_{T}\otimes _{S}X)=(X\otimes _{S}m_{T})\circ (v\otimes _{S}{T
})\circ ({T}\otimes _{S}v).
\end{equation*}
There is a cyclic object $Z_{\ast }(T/S,v)$, with $Z_{n}(T/S,v)=T^{{\hat{
\otimes}}_{S}\,n+1}{\hat{\otimes}}_{S}X$ whose face and degeneracy maps are
\begin{align*}
d_{k}(t_{0}{\widehat{\otimes }_{S}}\,t_{1}\widehat{\otimes
}_{S}\cdots \widehat{\otimes }_{S}\,t_{n}\widehat{\otimes }_{S}x)&
=\left\{
\begin{array}{ll}
t_{0}{\widehat{\otimes }_{S}}\,t_{1}\widehat{\otimes }{}_{S}\cdots \widehat{
\otimes }_{S}\,t_{k}t_{k+1}\widehat{\otimes }_{S}\cdots \widehat{\otimes }
_{S}\,t_{n}\widehat{\otimes }_{S}x, & \text{for }0\leq k<n, \\
(t_{n})_{v}t_{0}\widehat{\otimes }_{S}\,t_{1}\widehat{\otimes
}_{S}\cdots
\widehat{\otimes }_{S}\,t_{n-1}\widehat{\otimes }_{S}(x)_{v}\,, & \text{for }
k=n,
\end{array}
\right.  \\
s_{k}(t_{0}{\widehat{\otimes }_{S}}\,t_{1}{\widehat{\otimes }_{S}}\,\cdots {
\widehat{\otimes }_{S}}t_{n}\widehat{\otimes }_{S}x)& =t_{0}{\widehat{
\otimes }_{S}}\,t_{1}{\widehat{\otimes }_{S}}\,\cdots \widehat{\otimes }
_{S}\,t_{k}\widehat{\otimes }_{S}1_{T}{\widehat{\otimes
}_{S}}\cdots \widehat{\otimes }_{S}\,t_{n}\widehat{\otimes
}_{S}x,\quad \text{for }0\leq k\leq n.
\end{align*}
The para-cyclic map is given by 
\begin{equation*}
v_{n}(t_{0}\,\widehat{\otimes }_{S}\,\dots \,\widehat{\otimes }_{S}\,t_{n}
\widehat{\otimes }_{S}x)=(t_{n})_{v}\,\widehat{\otimes }_{S}\,t_{0}\,
\widehat{\otimes }_{S}\,\dots \,\widehat{\otimes }_{S}\,t_{n-1}\widehat{
\otimes }_{S}(x)_{v}.
\end{equation*}
\end{corollary}

\begin{proof}
In terms of the map $v$, we can equip $X\otimes _{S}T$ with a
$T$-bimodule structure by
\begin{equation*}
t_{1}(x\otimes _{S}t)t_{2}=v(t_{1}\otimes _{S}x)tt_{2}.
\end{equation*}
Moreover,
\begin{equation*}
v\otimes _{S}T:T\otimes _{S}X\otimes _{S}T\cong (T\otimes
_{S}T)\otimes _{T}(X\otimes _{S}T)\rightarrow (X\otimes
_{S}T)\otimes _{T}(T\otimes _{S}T)\cong X\otimes _{S}T\otimes
_{S}T
\end{equation*}
is a transposition map for the Sweedler's $T$-coring $T\otimes
_{S}T$ and the $T$-bimodule $X\otimes _{S}T$, in the sense of
\eqref{eq:w_dual}. Consequently, we can apply Corollary
\ref{ex:no.coef.dual} to Sweedler's coring in Example
\ref{ex:Sweedler}.  One proves that the corresponding para-cyclic
object has $Z_{\ast }(T/S,v)$ as underlying simplicial structure and $
v_{\ast }$ as para-cyclic map.
\end{proof}

\begin{remarks}
(i) Let us take $X=S.$ The canonical isomorphism $v:T\otimes
_{S}S\rightarrow S\otimes _{S}T$ satisfies the hypothesis of
Corollary \ref{ex:dual_no_coef}. The corresponding para-cyclic (in fact cyclic)
object was used in \cite{JaSt:CycHom} to define relative cyclic
homology. Moreover, this cyclic object and the cocyclic object
$Z^\ast (\mathcal{S}_T,w_T)$ in Corollary \ref{ex:no_coef} are (cyclic)
dual to each other. In the particular case when $R=k$ and
$\varphi=u_T$, we rediscover the cyclic object introduced by A.
Connes in order to define the cyclic homology of an algebra, cf.
\cite[p. 330]{We}. Thus $Z^\ast (\mathcal{S}_T,w_T)$ is the cyclic
dual of Connes' cyclic object.

(ii) Note that the construction of $Z_{\ast }(T/S,v)$ can be
performed for any algebra $T$ in a symmetric monoidal category $\mathcal{M}$,
by replacing everywhere $\widehat{\otimes }_{S}$ with $\otimes ,$ the tensor
product in $\mathcal{M}.$ Therefore, para-cyclic objects in \cite{Kay:UniHCyc}
and \cite{HaKhRaSo2} are examples of this type.

(iii) If $w:X\otimes _{S}T\rightarrow T\otimes _{S}X$ is an invertible $S$
-bimodule map satisfying \eqref{eq:bim_tr_map}, then $v=w^{-1}$
satisfies the relations in Corollary \ref{ex:dual_no_coef}.
Conversely,  in the case when the morphism $v$ in the above
construction is invertible, then its inverse is a transposition
map in the sense of \eqref{eq:bim_tr_map}. As a
matter of fact, the corresponding para-cocyclic object in Theorem \ref
{thm:para_cyc} is cyclic dual of the para-cyclic object in Corollary \ref
{ex:dual_no_coef}. This suggests a categorical approach to cyclic
duality, details of which will be studied elsewhere.
\end{remarks}}

Other examples of para-cyclic objects, obtained as applications of
Theorem \ref{thm:c_distr_law}, will be discussed in Theorem
\ref{thm:comod_coring} and Theorem \ref{thm:mod_coring}.

\section{Cyclic (co)homology of bialgebroids \label{sec:bialgebroid}}
\label{sec:bgd}

In this section we apply the categorical framework, obtained in
Section \ref{sec:adm.sep}, to examples provided by (co)module algebras and
(co)module corings of bialgebroids, and analyze the structure of the resulting
para-(co)cyclic objects.

\subsection{(Co)module algebras of bialgebroids}
\label{sec:(co)mod.alg}

In this section we consider admissible septuples
$\mathcal{S}_T$, coming from a $\mathbb{K}$-algebra homomorphism $
{\varphi}: { R}\to {  T}$ as in Proposition \ref{ex:bimod_Pi}. As
we have seen in Proposition \ref{ex:alg_sep}, $\mathcal{S}_T$
determines a cosimplex
\begin{equation}\label{eq:cosimp}
Z^n(\mathcal{S}_T, X)= T^{\wotR n+1}\wotR X,
\end{equation}
for any ${  R}$-bimodule $X$. Coface and codegeneracy maps are given by
\begin{equation}\label{eq:bim_coface&codeg}
\widehat{d}_k=T^{\wotR k}\,\wotR \,{\varphi}\,\wotR \,T^{\wotR
n-k}\,\wotR\, X\quad \textrm{and}\quad \widehat{s}_k=T^{\wotR
k}\,\wotR \, { m_T}\,\wotR \,T^{\wotR n-k}\,\wotR\, X,
\end{equation}
where {$k=0,\ldots,n$ and} $m_T$ denotes the multiplication
map $T\otimes_R T\to T$. Furthermore, by Theorem
\ref{thm:para_cyc}, the cosimplex $Z^\ast(\mathcal{S}_T, X)$ is
para-cocyclic provided that there is a transposition map
$w:X\otimes_{ R} {  T}\to { T}\otimes_{ R} X$. Conditions
(\ref{eq:bim_tr_map}) characterizing a transposition map are
reminiscent of some of the axioms of an entwining structure (over
an algebra ${  R}$), {cf. \cite[Section 2.3]{BB}. Main examples of
entwining structures over non-commutative algebras arise from
Doi-Koppinen data of bialgebroids (in the sense of \cite{BrCaMi}). In a
similar manner,} the aim
of this section is to construct canonical transposition maps in
the case when ${  T}$ is a (co)module algebra of an ${
R}$-bialgebroid $\mathcal{B}$ and $X$ is a
$\mathcal{B}$-(co)module.
\medskip

Bialgebroids can be thought of as a generalization of
bialgebras to arbitrary, non-commutative base algebras. The first
form of the structure that is known today as a left bialgebroid
was introduced by Takeuchi in \cite{Tak:bgd} under the name
$\times_{  R}$-bialgebra. Another definition and the name
`bialgebroid' was proposed by Lu in \cite{Lu:hgd}. The two
definitions were proven to be equivalent in \cite{BrzMil:bgd}.
`Left' and `right' versions of bialgebroids were defined in
\cite{KadSzl:D2bgd}.

\begin{definition}\label{def:left_bgd}
Consider an algebra ${  R}$ over a commutative ring $\mathbb{K}$. A
\emph{left bialgebroid} $\mathcal{B}$ over ${  R}$ consists of the
data $({ B}, {\xi}, {\zeta}, {\Delta},  {\epsilon})$. Here ${  B}$
is a $\mathbb{K}$-algebra and $ {\xi}$ and $ {\zeta}$ are $\mathbb{K}$-algebra
homomorphisms ${  R}\to { B}$ and ${  R}^{op}\to { B}$,
respectively, such that their ranges are commuting subalgebras in
${  B}$. In terms of the maps $ {\xi}$ and $ {\zeta}$, ${  B}$ can
be equipped with an   ${ R}$-bimodule structure as
\[
r_1\cdot b\cdot r_2:=  {\xi}(r_1)  {\zeta}(r_2) b,\qquad
  \textrm{for } r_1,r_2\in {    R}\quad \textrm{and } b\in {  B}.
\]
By definition, the coproduct $ {\Delta}:{  B}\to {  B}\otimes_{  R} { B}$ and
the counit $ {\epsilon}:{  B}\to {  R}$ equip this bimodule with an ${
R}$-coring structure.
Between the algebra and coring structures of ${  B}$ the
following compatibility axioms are required. For the coproduct we
introduce the index notation $ {\Delta}(b)= b_{(1)}\otimes_{  R}
b_{(2)}$, where implicit summation is understood.

\begin{itemize}
\item[{(i)}] $b_{(1)} {\zeta}(r)\otimes_{  R} b_{(2)}
  =b_{(1)}\otimes_{  R}
  b_{(2)} {\xi}(r)$, for $r\in {  R}$ and $b\in {  B}$.
\item[{(ii)}] $ {\Delta}(1_{  B})=1_{  B}\otimes_{  R} 1_{  B}$
and
  $ {\Delta}(bb')= b_{(1)}b'_{(1)}\otimes_{  R} b_{(2)} b'_{(2)}$,
  for $b,b'\in {  B}$.
\item[{(iii)}] $ {\epsilon}(1_{  B}) = 1_{  R}$ and $
{\epsilon}(bb')=  {\epsilon}
  \big(b  {\xi}( {\epsilon}(b')) \big)$, for $b,b'\in
  {  B}$.
\end{itemize}
\end{definition}

Axiom (i) in Definition \ref{def:left_bgd} needs to be imposed in
order for the second condition in axiom (ii) to make sense. Axiom
(iii) implies that also $ {\epsilon}(bb')=
 {\epsilon}\big(b
 {\zeta}( {\epsilon}(b')) \big)$, for $b,b'\in
{  B}$. It follows by the ${  R}$-module map properties, unitality
and multiplicativity of the coproduct $ {\Delta}$ that
\begin{equation}\label{eq:left_R_mods}
 {\Delta}\big( {\xi}(r_1)
 {\zeta}(r_2) b  {\xi}(r_3)  {\zeta}
(r_4)\big)=  {\xi}(r_1) b_{(1)} {\xi} (r_3) \otimes_{ R}
 {\zeta}(r_2) b_{(2)} {\zeta}(r_4),
\end{equation}
for $r_1,r_2,r_3, r_4\in {  R}$ and $b\in \mathcal{B}$. Since the
coproduct is coassociative, the Sweedler-Heynemann index notation
can be used. That is, for the iterated power of the coproduct we
write $( {\Delta}\otimes_{  R} \mathcal{B} \otimes_{ R} \dots
\otimes_{  R} \mathcal{B})\circ \dots \circ ( {\Delta}\otimes_{ R}
\mathcal{B})\circ  {\Delta}(b)=b_{(1)}\otimes_{  R} \dots
\otimes_{  R} b_{(n-1)}\otimes_{  R} b_{(n)}$, for any positive
integer $n$ and $b\in \mathcal{B}$.

Note that the axioms in Definition \ref{def:left_bgd} are not
invariant under changing the multiplication in ${  B}$ to the
opposite multiplication. Definition \ref{def:left_bgd} has a
symmetrical counterpart, known as a \emph{right
  bialgebroid}. For the details we refer to \cite{KadSzl:D2bgd}.

Definition \ref{def:left_bgd} is motivated by the following result
of Schauenburg. Consider two algebras ${  R}$ and ${
  B}$ over a commutative ring $\mathbb{K}$ and two $\mathbb{K}$-algebra
homomorphisms $ {\xi}:{
  R}\to {  B}$ and $ {\zeta}:{  R}^{op}\to {  B}$, whose
ranges are commuting subalgebras of ${  B}$. Clearly, in this
setting any (left or right) ${
  B}$ module can be equipped with an ${  R}$ bimodule structure using the
maps $ {\xi}$ and $ {\zeta}$. {For example, for a left ${
  B}$-module $V$ with action $\triangleright:{  B}\otimes_{\mathbb K} V\to
V$, one can define an ${R}$-bimodule structure as}
\[
r_1\otimes_{\mathbb K}  v \otimes_{\mathbb K}  r_2\mapsto  {\xi}(r_1)
 {\zeta}(r_2)\triangleright v.
\]
With respect to the resulting ${  R}$-actions, ${  B}$-module maps
are ${  R}$-bimodule maps. That is, there exists a forgetful
functor from the category of (left or right) ${  B}$-modules to
the category of ${  R}$-bimodules.

\begin{theorem}\label{thm:mod_cat_mon}
\cite[Theorem 5.1]{Scha:bia_nc} Consider two algebras ${  R}$ and
${
  B}$ over a commutative ring $\mathbb{K}$ and two $\mathbb{K}$-algebra
homomorphisms
$ {\xi}:{  R}\to {  B}$ and $ {\zeta}:{ R}^{op}\to { B}$, whose
ranges are commuting subalgebras of ${ B}$. There exists a right
(resp. left) bialgebroid $({
  B}, {\xi}, {\zeta}, {\Delta},
 {\epsilon})$ if and only if the  forgetful functor from
the
  category of right (resp. left) ${  B}$-modules to the category of ${
  R}$-bimodules is strict monoidal. That is, ${  R}$ is a right (resp. left)
  ${  B}$-module and the ${  R}$-module tensor product of two right
  (resp. left) ${  B}$-modules is a right (resp. left) ${  B}$-module.
\end{theorem}
{
Similarly to the case of a bialgebra, in Theorem \ref{thm:mod_cat_mon}
for a left ${  R}$-bialgebroid
$({B}, {\xi}, {\zeta}, {\Delta}, {\epsilon})$ the following ${  B}$-actions
are used
on ${  R}$, and on the ${  R}$-module tensor product
of two left ${B}$-modules $V$ and $W$, respectively}.
\begin{equation}\label{eq:mod_mon_prod}
b\triangleright r:= {\epsilon}\big(b
 {\xi}(r)\big)\qquad \textrm{and}\qquad b\triangleright
(v\otimes_{  R} w) :=  b_{(1)} \triangleright v \otimes_{
  R}  b_{(2)}\triangleright w,
\end{equation}
for $r\in {  R}$, $v\otimes_{  R} w\in V\otimes_{  R} W$ and $b\in
{B}$.
{It was proven in \cite[Theorem 5.1]{Scha:bia_nc} that the diagonal action
  in the second equation in (\ref{eq:mod_mon_prod}) is meaningful by axiom (i)
  in Definition \ref{def:left_bgd}.}

In light of Theorem \ref{thm:mod_cat_mon}, one can speak about
right (resp. left) module algebras of a right (resp. left)
bialgebroid $\mathcal{B}$, i.e. about algebras in the monoidal
category of right (resp. left) $\mathcal{B}$-modules.
\begin{definition}
Consider an algebra ${  R}$ over a commutative ring $\mathbb{K}$ and a left
${  R}$-bialgebroid $\mathcal{B}$. A \emph{left
  $\mathcal{B}$-module algebra} is a $\mathbb{K}$-algebra and left
$\mathcal{B}$-module $T$, with $\mathcal{B}$-action
$\triangleright:\mathcal{B} \ot_{\mathbb K} T\to T$, such that the
multiplication in $T$ is $R$-balanced and
\begin{equation}\label{eq:mod_alg}
b \triangleright 1_{  T} =
 {\xi}\big( {\epsilon}(b)\big)\triangleright
1_{  T} \qquad  \textrm{and}\qquad b \triangleright (tt') =
(b_{(1)} \triangleright t)(b_{(2)} \triangleright t'),
\end{equation}
for $b\in {  B}$ and $t,t'\in {  T}$.
\end{definition}

For example, the constituent algebra in a left bialgebroid
$\mathcal{B}$ is
  itself a (so called \emph{left  regular}) left
  $\mathcal{B}$-module algebra via the action given by the
  product in $\mathcal{B}$.

Note that for a left module algebra ${  T}$ of a left ${
R}$-bialgebroid $\mathcal{B}=({ B}, {\xi}, {\zeta},
 {\Delta}, {\epsilon})$, there is a canonical
$\mathbb{K}$-algebra homomorphism ${  R}\to {  T}$, $r\mapsto
 {\xi}(r)\triangleright 1_{  T}$. Hence there is a
corresponding admissible septuple as in Proposition
\ref{ex:alg_sep}.

A left comodule of a left ${  R}$-bialgebroid
$\mathcal{B}=(B,\xi,\zeta,\Delta,\epsilon)$ is defined as a left comodule of
the underlying $R$-coring $(\mathcal{B},\Delta,\epsilon)$. That is, a left
$R$-module $X$, together with a left $R$-module map $X\to \mathcal{B}
\otimes_R X$, $x\mapsto x_{[-1]}\otimes_R x_{[0]}$ (where implicit
summation is understood), satisfying coassociativity and counitality
axioms. Note that a left $\mathcal{B}$-comodule $X$ can be equipped with
an $R$-bimodule structure by introducing a right $R$-action
\begin{equation}\label{eq:right_R_mod}
x\cdot r:= {\epsilon}\big( x_{[-1]}{\xi}(r) \big)\cdot x_{[0]},\qquad
\textrm{for } r\in {  R}\quad \textrm{and}\quad x\in X.
\end{equation}
With respect to the resulting bimodule structure, $\mathcal{B}$-comodule maps
are $R$-bimodule maps.
In particular, the left $\mathcal{B}$-coaction on $X$ is an $R$-bimodule map
in the sense that, for $r,r'\in R$ and $x\in X$,
\begin{equation}\label{eq:coac_mod}
(r\cdot x \cdot r')_{[-1]}\otimes_R (r\cdot x \cdot r')_{[0]} =
\xi(r) x_{[-1]} \xi(r')\otimes_R x_{[0]}.
\end{equation}
Furthermore, for any $x\in X$ and $r\in R$,
\begin{equation}\label{eq:coac_centr}
x_{[-1]} \otimes_R x_{[0]} \cdot r=x_{[-1]} \zeta(r)\otimes_R x_{[0]}.
\end{equation}
\begin{theorem}\label{thm:mod_alg}
Let ${  R}$ be an algebra over a commutative ring $\mathbb{K}$ and let
$\mathcal{B}$ be a left bialgebroid over ${  R}$. Consider a left
$\mathcal{B}$-module algebra ${  T}$ with $\mathcal{B}$-action
$\triangleright$ and a left $\mathcal{B}$-comodule $X$ with
coaction $x\mapsto x_{[-1]}\otimes_{  R} x_{[0]}$ (where implicit
summation is understood). Then a transposition map for the admissible septuple
$\mathcal{S}_T$, associated via Proposition \ref{ex:alg_sep} to the
$\mathbb{K}$-algebra map $R\to T$, $r\mapsto \xi(r) \triangleright 1_T$, is
given by
\begin{equation}\label{eq:mod_alg_w}
w:X\otimes_{  R} {  T} \to {  T}\otimes_{  R} X,\qquad x\otimes_{
R} t \mapsto x_{[-1]}\triangleright t \otimes_{  R} x_{[0]}.
\end{equation}
Hence the cosimplex (\ref{eq:cosimp}) admits a para-cocyclic
structure
\begin{equation}\label{eq:mod_alg_hw}
{w}_n(t_0 \,\widehat{\otimes}_{  R}\,\dots\, \widehat{\otimes}_{
  R}\, t_{n}\, \widehat{\otimes}_{  R}\, x)=
t_1\,\widehat{\otimes}_{  R}\,  \dots \, \widehat{\otimes}_{ R}\,
t_{n}\,
  \widehat{\otimes}_{  R}\, x_{[-1]}\triangleright t_0  \,
  \widehat{\otimes}_{  R}\, x_{[0]}.
\end{equation}
\end{theorem}
\begin{proof}
The map (\ref{eq:mod_alg_w}) is a well defined left $R$-module homomorphism by
(\ref{eq:coac_mod}). Its right $R$-module map property follows by
(\ref{eq:coac_centr}). Conditions (\ref{eq:bim_tr_map}) follow by definition
(\ref{eq:mod_alg}) of a module algebra {as it is shown below. Denote by
  $\varphi$
  the algebra homomorphism $R\to T$, $r\mapsto \xi(r)\triangleright
  1_T=\zeta(r)\triangleright 1_T$. Omitting canonical isomorphisms $R \ot_R X
  \cong X \cong X \ot_R R$, we have
\begin{eqnarray*}
&&\big(w\circ(X\ot_R \varphi)\big)(x)=
x_{[-1]}\triangleright 1_T \ot_R x_{[0]}=
\zeta\big(\epsilon(x_{[-1]})\big)\triangleright 1_T \ot_R x_{[0]}=
1_T \ot_R \epsilon(x_{[-1]}) \cdot x_{[0]}=\\
&&\qquad 1_T \ot_R x =
\big((\varphi\ot_R X)\big)(x),\\
&&\big((m_T \ot_R X)\circ (T \ot_R w) \circ (w\ot_R T)\big)(x\ot_R t
\ot_R t')=
(x_{[-1]}\triangleright t)(x_{[0][-1]}\triangleright t') \ot_R x_{[0][0]} =\\
&&\qquad (x_{[-1](1)}\triangleright t)(x_{[-1](2)}\triangleright t') \ot_R
x_{[0]} =
x_{[-1]}\triangleright tt'\ot_R x_{[0]} =
\big(w\circ(X\ot_R m_T)\big)(x\ot_R t \ot_R t').
\end{eqnarray*}
}
\end{proof}

Analogously to (\ref{eq:right_R_mod}), also a right comodule $V$ of a right
  $R$-bialgebroid $\mathcal{B}=(B,\xi,\zeta,\Delta,\epsilon)$
can be equipped with an $R$-bimodule structure by
  introducing a left $R$-action
\begin{equation}\label{eq:left_R_mod}
r\cdot v:=v^{[0]}\cdot
 {\epsilon}\big( {\xi}(r)v^{[1]}\big),\qquad
\textrm{for } r\in {  R}\quad \textrm{and}\quad v\in V,
\end{equation}
where $v\mapsto v^{[0]}\otimes_R v^{[1]}$ denotes the right
$\mathcal{B}$-coaction on $V$, with implicit summation understood.
Hence there exists a forgetful functor from the category of right
$\mathcal{B}$-comodules to the category of ${
  R}$-bimodules. With this observation in mind, the next theorem follows by
a symmetrical form of \cite[Proposition 5.6]{Scha:bia_nc}.

\begin{theorem}\label{thm:comod_cat_mon}
Consider an algebra ${  R}$ over a commutative ring $\mathbb{K}$ and a
right ${
  R}$-bialgebroid $\mathcal{B}$.
Then the forgetful functor from the category of right
$\mathcal{B}$-comodules
  to the category of ${  R}$-bimodules is strict monoidal. That is, ${
  R}$ is a right $\mathcal{B}$-comodule and the ${  R}$-module tensor
  product of two right $\mathcal{B}$-comodules is a right
  $\mathcal{B}$-comodule.
\end{theorem}

{Similarly to a bialgebra, in Theorem \ref{thm:comod_cat_mon} for a right ${
  R}$-bialgebroid $\mathcal{B}=({B}, {\xi}, {\zeta}, {\Delta}, {\epsilon})$
the following $\mathcal{B}$-coactions on ${  R}$, and on the
  ${  R}$-module tensor product $V\otimes_{  R} W$ of two right
  $\mathcal{B}$-comodules $V$ and $W$, are used.}
\[
{  R}\to {  R}\otimes_{  R} {\mathcal{B}} \cong \mathcal{B},\quad
r\mapsto  {\xi}(r) \quad \textrm{and}\quad V\otimes_{ R} W \to
V\otimes_{  R} W\otimes_{  R} {\mathcal{B}}, \quad v\otimes_{  R}
w\mapsto v^{[0]}\otimes_{  R} w^{[0]}\otimes_{ R} v^{[1]} w^{[1]}.
\]
The coaction on $V\otimes_R W$ is well defined by the right bialgebroid
versions of properties (\ref{eq:coac_mod}) and (\ref{eq:coac_centr}), i.e. the
identities
\begin{equation}\label{eq:right_coac_mod}
(r\cdot v\cdot r')^{[0]}\otimes_R (r\cdot v\cdot r')^{[1]} =
v^{[0]}\otimes_R \xi(r) v^{[1]} \xi(r')
\qquad \textrm{and}\qquad
r\cdot v^{[0]} \otimes_R v^{[1]}=v^{[0]}\otimes_R \zeta(r) v^{[1]},
\end{equation}
for $r,r'\in R$ and $v\in V$.
Symmetrically, the forgetful functors from the category of left
comodules of a right ${  R}$-bialgebroid, and from the categories
of right or left comodules of a left ${  R}$-bialgebroid, to the
category of $R$-bimodules are strict (anti-)monoidal.

In light of Theorem \ref{thm:comod_cat_mon}, one can speak about
right comodule algebras of a right bialgebroid $\mathcal{B}$, i.e.
about algebras in the monoidal category of right
$\mathcal{B}$-comodules.
\begin{definition}
Consider an algebra ${  R}$ over a commutative ring $\mathbb{K}$ and a
right ${  R}$-bialgebroid $\mathcal{B}$. A \emph{right
  $\mathcal{B}$-comodule algebra} is a $\mathbb{K}$-algebra and right
$\mathcal{B}$-comodule ${  T}$, with coaction $t\mapsto
t^{[0]}\otimes_{
  R} t^{[1]}$, such that the multiplication in ${  T}$ is
${  R}$-balanced and, for $t,t'\in {  T}$,
\begin{equation}\label{eq:com_alg}
1_{  T}^{[0]} \otimes_{  R} 1_{  T}^{[1]} = 1_{  T} \otimes_{ R}
1_{  T} \qquad \textrm{and}\qquad (tt')^{[0]} \otimes_{
R}(tt')^{[1]} = t^{[0]} t^{\prime [0]} \otimes_{
  R} t^{[1]} t^{\prime [1]}.
\end{equation}
\end{definition}

For example, the constituent algebra in a right bialgebroid
$\mathcal{B}$ is
  itself a (so called \emph{right  regular}) right
  $\mathcal{B}$-comodule algebra via the coaction given by the
  coproduct in $\mathcal{B}$.

Note that the second condition in (\ref{eq:com_alg}) is meaningful
since the multiplication in ${  T}$ is ${  R}$-balanced and the second
condition in (\ref{eq:right_coac_mod}) holds.
For a right comodule algebra ${  T}$ of a right ${ R}$-bialgebroid
$\mathcal{B}$, there is a canonical $\mathbb{K}$-algebra homomorphism ${
R}\to {
  T}$, $r\mapsto r\cdot 1_{  T}=1_{  T}\cdot r$, in terms of which $r\cdot
t=( r\cdot 1_{  T}) t$ and $t\cdot r = t( r\cdot 1_{  T})$. Hence
there is a corresponding admissible septuple as in Proposition
\ref{ex:alg_sep}.

\begin{theorem}\label{thm:comod_alg}
Let ${  R}$ be an algebra over a commutative ring $\mathbb{K}$ and let
$\mathcal{B}$ be a right bialgebroid over ${  R}$. Consider a
right $\mathcal{B}$-comodule algebra ${  T}$, with coaction
$t\mapsto t^{[0]}\otimes_{  R} t^{[1]}$ (where implicit summation
is understood) and a right $\mathcal{B}$-module $X$ with action
$\triangleleft$. Then a transposition map for the admissible septuple
$\mathcal{S}_T$, associated via Proposition \ref{ex:alg_sep} to the
$\mathbb{K}$-algebra map $R\to T$, $r\mapsto r\cdot 1_T=1_T\cdot r$, is given
by
\[
w:X\otimes_{  R} {  T} \to {  T}\otimes_{  R} X,\qquad x\otimes_{
R} t \mapsto t^{[0]} \otimes_{  R} x\triangleleft t^{[1]}.
\]
Hence the cosimplex (\ref{eq:cosimp}) admits a para-cocyclic
structure
\[
{w}_n( t_0 \,\widehat{\otimes}_{  R}\,\dots\, \widehat{\otimes}_{
  R}\, t_{n}\, \widehat{\otimes}_{  R}\, x)=
t_1\, \widehat{\otimes}_{  R}\,
  \dots \, \widehat{\otimes}_{  R} t_{n}\, \widehat{\otimes}_{  R}\,
t_0^{[0]}\, \widehat{\otimes}_{  R}\, x\triangleleft t_0^{[1]}.
\]
\end{theorem}
\begin{proof}
The map $w$ is a well defined $R$-bimodule homomorphism by
(\ref{eq:right_coac_mod}).
Conditions (\ref{eq:bim_tr_map}) follow by definition
(\ref{eq:com_alg}) of a comodule algebra.
\end{proof}

\subsection{(Co)module corings of bialgebroids}
\label{sec:(co)mod_coring}

In this section we consider comonads $\O:={  C}\otimes_{  R} (-)$ on the
category of ${  R}$-bimodules as in Proposition \ref{ex:coring}, determined by
a coring $({ C}, { {\Delta}}, {\epsilon})$ over an algebra ${ R}$. Let
$\bs{\Pi}$ be the functor in Definition \ref{def:Pi}. By Proposition
\ref{ex:coring}, for any ${R}$-bimodule $X$, there is an
associated simplicial module $Z_\ast(\mathcal{S}_C, X)$, which in
degree $n$ is given by
\begin{equation}\label{eq:coring_simp}
{Z_n(\mathcal{S}_C, X)=C^{\wotR n+1}\wotR X.}
\end{equation}
{Face and degeneracy  maps in degree $n$}  are
\begin{equation}\label{eq:bim_face&deg}
\widehat{d}_k=C^{\wotR k}\,\wotR \,{\epsilon_C}\,\wotR \,C^{\wotR
n-k}\,\wotR\, X\quad \textrm{and}\quad \widehat{s}_k=C^{\wotR
k}\,\wotR \, { \Delta_C}\,\wotR \,C^{\wotR n-k}\,\wotR\, X,
\end{equation}
{for $k\in \{0,\dots,n\}$.}
As we have seen in Proposition \ref{ex:coring}, we can choose a
second comonad $\I:=(-)\otimes_{  R} {  C}$ on the category of ${
R}$-bimodules, and a natural transformation $\boldsymbol{i}: (-)\,
\widehat{\otimes}_{  R}\, {
  C} \to {  C}\, \widehat{\otimes}_{  R} (-)$, given by the flip map.
Taking the trivial natural transformation ${\boldsymbol{t}}: {
  C}\otimes_{  R} (-)\otimes_{  R} {  C} \to {  C}\otimes_{  R}
(-)\otimes_{  R} {  C}$, the conditions in (\ref{eq:i_dual}) hold.
By Theorem \ref{thm:para_cyc_dual}, the simplex
$Z_\ast(\mathcal{S}_C, X)$ is para-cyclic provided that there
exists a morphism $w:{ C}\otimes_{  R} X \to X \otimes_{  R} {
C}$, satisfying (\ref{eq:transition.map.dual}). Note the
similarity of conditions (\ref{eq:transition.map.dual}) to some of
the axioms of an entwining structure over ${  R}$. Similarly to
the way Doi-Koppinen data determine entwining structures, in this
section we construct dual transposition maps
in the case when ${ C}$ is a
(co)module coring of a bialgebroid $\mathcal{B}$ and $X$ is a
$\mathcal{B}$-(co)module.
\medskip

In light of (a symmetrical version of) Theorem
\ref{thm:comod_cat_mon}, one can speak about left comodule corings
of a left bialgebroid $\mathcal{B}$, i.e. about coalgebras in the
monoidal category of left $\mathcal{B}$-comodules.
\begin{definition}
Consider an algebra ${  R}$ over a commutative ring $\mathbb{K}$ and a left
$R$-bialgebroid $\mathcal{B}=(B,\xi,\zeta,\Delta,\epsilon)$. A
\emph{left $\mathcal{B}$-comodule coring} is an $R$-coring and left
$\mathcal{B}$-comodule ${  C}$, with one and the same underlying $R$-bimodule
structure,
such that for $c\in { C}$
\begin{equation}\label{eq:comod_coring}
 {\xi}\big( {\epsilon}_{  C}(c)\big) =
 {\zeta}\big( {\epsilon}_{  C}(c_{[0]})\big)
c_{[-1]} \quad \textrm{and}\quad c_{[-1]}\stac {  R} c_{[0](1)}
\stac {  R} c_{[0](2)} = c_{(1)[-1]} c_{(2)[-1]}\stac {  R}
c_{(1)[0]}\stac {  R}c_{(2)[0]},
\end{equation}
where $\epsilon_C$ is the counit and
$\Delta_{  C}:c\mapsto c_{(1)}\otimes_{  R} c_{(2)}$ is the coproduct of $C$
and $c\mapsto c_{[-1]} \otimes_{  R} c_{[0]}$ is the $\mathcal{B}$-coaction on
$C$.
\end{definition}
The second condition in (\ref{eq:comod_coring}) is meaningful by
(\ref{eq:coac_mod}) and (\ref{eq:coac_centr}).
For example, the constituent ${  R}$-coring in a left ${
R}$-bialgebroid
  $\mathcal{B}$ is itself a (so called \emph{left  regular}) left
  $\mathcal{B}$-comodule coring via the coaction given by the coproduct.

\begin{theorem}\label{thm:comod_coring}
Let ${  R}$ be an algebra over a commutative ring $\mathbb{K}$ and let
$\mathcal{B}$ be a left bialgebroid over ${  R}$. Consider a left
$\mathcal{B}$-comodule coring ${  C}$ with coaction $c\mapsto
c_{[-1]}\otimes_{  R} c_{[0]}$, and a left $\mathcal{B}$-module
$X$ with action $\triangleright$. Then a transposition map $w:{ C}\otimes_{ R}
X \to X \otimes_{  R} {  C}$
for the dual admissible septuple $\mathcal{S}_C$, associated via Proposition
\ref{ex:coring} to the $R$-coring $C$, is given by  $w(c\otimes_{  R}
x):=c_{[-1]}\triangleright x \otimes_{  R} c_{[0]}$. Hence the
simplex (\ref{eq:coring_simp}) admits a para-cyclic structure
\[
{w}_n( c_0 \,\widehat{\otimes}_{  R}\,\dots\, \widehat{\otimes}_{
  R}\, c_{n}\, \widehat{\otimes}_{  R}\, x )=
c_{n [0]}\, \widehat{\otimes}_{  R}\, c_0\,
  \widehat{\otimes}_{  R}\, \dots \widehat{\otimes}_{  R}\, c_{n-1}\,
  \widehat{\otimes}_{  R}\,  c_{n [-1]}\triangleright x.
\]
\end{theorem}
\begin{proof}
The map $w$ is a well defined $R$-bimodule homomorphism by (\ref{eq:coac_mod})
and (\ref{eq:coac_centr}).
Conditions (\ref{eq:transition.map.dual}) follow by definition
(\ref{eq:comod_coring}) of a comodule coring.
\end{proof}

In light of Theorem \ref{thm:mod_cat_mon}, one can speak about
right (resp. left) module corings of a right (resp. left)
bialgebroid $\mathcal{B}$, i.e. about coalgebras in the monoidal
category of right (resp. left) $\mathcal{B}$-modules.
\begin{definition}
Consider an algebra $R$ over a commutative ring $\mathbb{K}$ and a right
 $R$-bialgebroid $\mathcal{B}=(B,\xi,\zeta,\Delta,\epsilon)$. A
\emph{right $\mathcal{B}$-module coring} is an $R$-coring and right
$\mathcal{B}$-module $C$, with one and the same underlying $R$-bimodule
structure,
such that for $c\in C$ and $b\in B$
\begin{equation}\label{eq:mod_coring}
\epsilon_{  C}(c\triangleleft b)= \epsilon\big(  {\xi}( \epsilon_{
C}(c) ) b \big) \qquad \textrm{and}\qquad (c\triangleleft
b)^{(1)}\otimes_{  R} (c\triangleleft b)^{(2)} =
c^{(1)}\triangleleft b^{(1)}\otimes_{  R} c^{(2)}\triangleleft
b^{(2)},
\end{equation}
where ${\epsilon}_{ C}$ is the counit and $\Delta_C:c\mapsto c^{(1)}\otimes_{
  R} c^{(2)}$ is the coproduct in $C$, the symbol $\triangleleft$ denotes the
  $\mathcal{B}$-action on $C$ and for the coproduct in $\mathcal{B}$
the index notation $\Delta:b\mapsto b^{(1)}\otimes_{  R} b^{(2)}$ is used,
  with implicit summation understood.
\end{definition}

For example, the constituent ${  R}$-coring in a right ${
R}$-bialgebroid
  $\mathcal{B}$ is itself a (so called \emph{right  regular}) right
  $\mathcal{B}$-module coring via the action given by the product in
  $\mathcal{B}$.

\begin{theorem}\label{thm:mod_coring}
Let ${  R}$ be an algebra over a commutative ring $\mathbb{K}$ and let
$\mathcal{B}$ be a right bialgebroid over ${  R}$. Consider a
right $\mathcal{B}$-module coring ${  C}$ with action
$\triangleleft$, and a right $\mathcal{B}$-comodule $X$ with
coaction $x\mapsto x^{[0]}\otimes_{  R} x^{[1]}$. Then a transposition map
$w:{C}\otimes_{  R} X \to X \otimes_{  R} {  C}$,
for the dual admissible septuple $\mathcal{S}_C$, associated via Proposition
\ref{ex:coring} to the $R$-coring $C$, is given by $w(c\otimes_{  R}
x):=x^{[0]} \otimes_{  R} c\triangleleft x^{[1]}$. Hence the
simplex (\ref{eq:coring_simp}) admits a para-cyclic structure
\begin{equation}\label{eq:mod_coring_hw}
{w}_n( c_0 \,\widehat{\otimes}_{  R}\,\dots\,
\widehat{\otimes}_{
  R}\, c_{n}\, \widehat{\otimes}_{  R}\, x )=
c_n\triangleleft x^{[1]}\, \widehat{\otimes}_{  R}\, c_0\,
  \widehat{\otimes}_{  R}\, \dots \, \widehat{\otimes}_{  R}\, c_{n-1}\,
  \widehat{\otimes}_{  R}\,   x^{[0]}.
\end{equation}
\end{theorem}
\begin{proof}
The map $w$ is a well defined $R$-bimodule homomorphism by
(\ref{eq:right_coac_mod}).
Conditions (\ref{eq:transition.map.dual}) follow by definition
(\ref{eq:mod_coring}) of a module coring.
\end{proof}

\subsection{Stable anti Yetter-Drinfel'd modules of $\times_R$-Hopf algebras}
\label{sec:st.a.Yet.Dri} For a right module coring $C$ and right
comodule $X$ of a right bialgebroid
$\mathcal{B}=(B,\xi,\zeta,\Delta,\epsilon)$ over a $\mathbb{K}$-algebra $R$,
there is a $\mathbb{K}$-module isomorphism
 \[C\,^{\widehat{\otimes}_Rn+1}\widehat{\otimes}_R X
\cong C\,^{{\otimes}_{ R}n+1} \otimes_{R^e} X.\] Assume that $X$
has an additional left $\mathcal{B}$-module structure. In this
case, corresponding to the
    $\mathbb{K}$-algebra homomorphism ${  R}^e \to {  B}$,
    $r\otimes_{\mathbb K} r'\mapsto  {\xi}(r)  {\zeta}(r') =
     {\zeta}(r')  {\xi}(r)$, there is a canonical
    epimorphism $C^{\widehat{\otimes}_{  R}\, n+1} \widehat{\otimes}_{  R}
    \,X\to {  C}^{\otimes_{  R}\, n+1}\otimes_{\mathcal{B}} X$, where
    $C^{\otimes_{  R} \, n+1}$ is understood to be a right
    $\mathcal{B}$-module via the diagonal action $(c_1\ot_R\dots \ot_R
 c_{n+1})\triangleleft b := c_1 \triangleleft b^{(1)} \ot_R\dots \ot_R c_{n+1}
 \triangleleft b^{(n+1)}$,
 given by the iterated coproduct in $\mathcal{B}$.
{It is a well defined action by (a symmetrical version of) \cite[Theorem
    5.1]{Scha:bia_nc}.}

\begin{lemma}\label{lem:dual_B_mod_prod}
Let ${  R}$ be an algebra over a commutative ring $\mathbb{K}$ and let
$\mathcal{B}$ be a right bialgebroid over ${  R}$. Consider a
right $\mathcal{B}$-module coring ${  C}$ and a left
$\mathcal{B}$-module right $\mathcal{B}$-comodule $X$.
Then the simplex in Theorem \ref{thm:mod_coring} projects to a
simplex ${  C}^{\otimes_{  R}\, n+1}\otimes_{\mathcal{B}} X$.
\end{lemma}
\begin{proof}
Since the coproduct $\Delta_C$ and the counit $\epsilon_C$ of $C$ are
  right $\mathcal{B}$-module maps by definition, face and degeneracy maps of
  the simplex in Example \ref{ex:no.coef.dual} are right $\mathcal{B}$-module
  maps with respect to the diagonal action. Hence we can take their tensor
  product with the identity map on $X$ over the algebra $\mathcal{B}$,
  yielding a simplex as stated.
\end{proof}

The task of this section is to find criteria for the cyclicity of the quotient
simplex in Lemma \ref{lem:dual_B_mod_prod}. In order to do so, some
restriction on the involved bialgebroid is needed.

\begin{definition}\label{def:x_R_Hopf}
\cite[Theorem and Definition 3.5]{Scha:dual.double.bgd} Let ${ R}$
be an algebra over a commutative ring $\mathbb{K}$. A right ${
  R}$-bialgebroid $\mathcal{B}=({
  B}, {\xi}, {\zeta}, {\Delta},
 {\epsilon})$ is said to be a right \emph{$\times_{
R}$-Hopf algebra} provided that the map
\begin{equation}\label{eq:x_R_Hopf_a}
\vartheta: \mathcal{B}\otimes_{{  R}^{op}} \mathcal{B}\to
\mathcal{B}\otimes_{  R} \mathcal{B}, \qquad b \otimes_{{ R}^{op}}
b'\mapsto bb^{\prime (1)} \otimes_{  R} b^{\prime(2)}
\end{equation}
is bijective. In the domain of the map in (\ref{eq:x_R_Hopf_a}),
${R}^{op}$-module structures are given by right and left multiplication by
  $ {\zeta}(r)$, for $r\in {  R}$. In the codomain of the map in
  (\ref{eq:x_R_Hopf_a}), ${  R}$-module structures are given by right
  multiplication by $ {\xi}(r)$ and $ {\zeta}(r)$, for
  $r\in {  R}$.
\end{definition}

The notion of a $\times_{  R}$-Hopf algebra extends that of a Hopf
algebra. Indeed, if $\mathcal{B}$ is a bialgebra over a {\em
commutative} ring ${  R}$, with coproduct $b\mapsto
b^{(1)}\otimes_{  R} b^{(2)}$, then the map (\ref{eq:x_R_Hopf_a})
is bijective if and only if $\mathcal{B}$ is a Hopf algebra. In
this case the inverse is given in terms of the antipode $S$ as
$\vartheta^{-1}(b\otimes_{  R} b'):= bS(b^{\prime (1)}) \otimes_{
R} b^{\prime(2)}$.

For an algebra ${  R}$, consider a right $\times_{  R}$-Hopf
algebra $\mathcal{B}$. Since the map $\vartheta$ in
(\ref{eq:x_R_Hopf_a}) is a left $\mathcal{B}$-module map, its
inverse is determined by the restriction
$\vartheta^{-1}(1_{\mathcal B} \otimes_{  R} b)=: b_{-}\otimes_{{
R}^{op}} b_{+}$, where implicit summation is understood. Lemma
\ref{lem:x_R_Hopf_a}, which is essentially a symmetrical version
of \cite[Proposition 3.7]{Scha:dual.double.bgd}, collects
properties of the map $b\mapsto b_{-}\otimes_{{  R}^{op}} b_{+}$.

\begin{lemma}\label{lem:x_R_Hopf_a}
For an algebra ${  R}$, consider a right $\times_{  R}$-Hopf
algebra $\mathcal{B}=({ B}, {\xi}, {\zeta}, {\Delta},
 {\epsilon})$. Write
$ {\Delta}(b)=:b^{(1)}\otimes_{  R} b^{(2)}$ for the coproduct and
in terms of the map (\ref{eq:x_R_Hopf_a}) put $b_{-}\otimes_{{
R}^{op}} b_{+} := \vartheta^{-1}(1_{\mathcal B} \otimes_{  R} b)$.
The following identities hold, for $b,b'\in \mathcal{B}$ and $r\in
{  R}$.
\begin{itemize}
\item[{(i)}] $b_{-} {b_{+}}^{(1)}\otimes_{  R} {b_{+}}^{(2)} =
1_{\mathcal
  B} \otimes_{  R} b$
\item[{(ii)}] $b^{(1)} {b^{(2)}}_{-}\otimes_{{
R}^{op}}{b^{(2)}}_{+} =
  1_{\mathcal B} \otimes_{{  R}^{op}} b$
\item[{(iii)}] $(bb')_{-}\otimes_{{  R}^{op}}(bb')_{+} =b'_{-}
b_{-}
  \otimes_{{  R}^{op}} b_{+} b'_{+}$
\item[{(iv)}] $1_{\mathcal{B} \, -} \otimes_{{  R}^{op}}
1_{\mathcal{B} \, +} =1_{\mathcal{B}} \otimes_{{  R}^{op}}
1_{\mathcal{B}}$ \item[{(v)}] $b_{-} \otimes_{{  R}^{op}}
{b_{+}}^{(1)}\otimes_{  R}
  {b_{+}}^{(2)} = {b^{(1)}}_{-} \otimes_{{  R}^{op}} {b^{(1)}}_{+}
  \otimes_{  R} b^{(2)}$
\item[{(vi)}] ${b_{-}}^{(1)}\otimes_{  R} {b_{-}}^{(2)} \otimes_{{
    R}^{op}} b_{+} = b_{+ -} \otimes_{  R} b_{-}\otimes_{{  R}^{op}}
    b_{++}$
\item[{(vii)}] $b=
 {\zeta}\big( {\epsilon}(b_{-})\big)
  b_{+}$
\item[{(viii)}] $b_{-} b_{+} =
   {\xi}\big( {\epsilon}(b)\big)$
\item[{(ix)}] $ {\zeta}(r) b_{-} \otimes_{{  R}^{op}} b_{+} =
b_{-} \otimes_{{  R}^{op}} b_{+} {\zeta}(r)$.
\end{itemize}
\end{lemma}

Next Definition \ref{def:Y-D_mod} extends \cite[Definition
4.1]{JaSt:CycHom} or \cite[Definition 2.1]{HaKhRaSo1}.

\begin{definition}\label{def:Y-D_mod}
For an algebra ${  R}$, consider a right $\times_{  R}$-Hopf
algebra $\mathcal{B}=({ B}, {\xi}, {\zeta}, {\Delta},
 {\epsilon})$. Let $X$ be a right $\mathcal{B}$-comodule
and left $\mathcal{B}$-module. Denote the right
$\mathcal{B}$-coaction on $X$ by $x\mapsto x^{[0]}\otimes_{  R}
x^{[1]}$, for $x\in X$ (where implicit summation is understood)
and denote the left $\mathcal{B}$-action by $b\triangleright x$,
for $b\in \mathcal{B}$ and $x\in X$. We say that $X$ is an
\emph{anti Yetter-Drinfel'd module} provided that the following
conditions hold.

(i) The ${  R}$-bimodule structures of $X$, underlying its module
  and comodule structures, coincide. That is, for $x\in X$ and $r\in {  R}$,
\begin{equation}\label{eq:Y-D_R_mod}
x\cdot r =  {\zeta}(r)\triangleright x\qquad \textrm{and}\qquad
r\cdot x =  {\xi}(r)\triangleright x,
\end{equation}
where $x\cdot r$ denotes the right ${  R}$-action on the right
$\mathcal{B}$-comodule $X$ and $r\cdot x$
is the canonical left ${   R}$-action (\ref{eq:left_R_mod}) coming
from the right $\mathcal{B}$-coaction.

(ii) For $b\in \mathcal{B}$ and $x\in X$,
\begin{equation}\label{eq:Y-D_comp}
(b\triangleright x)^{[0]}\otimes_{  R}(b\triangleright x)^{[1]} =
{b^{(1)}}_{+} \triangleright x^{[0]}\otimes_{  R} b^{(2)} x^{[1]}
{b^{(1)}}_{-},
\end{equation}
where for the coproduct $ {\Delta}$ and the inverse of the map
(\ref{eq:x_R_Hopf_a}) the respective index notations (with
implicit summation), $ {\Delta}(b)=b^{(1)}\otimes_{ R} b^{(2)}$,
and $\vartheta^{-1}(1_{\mathcal B}\otimes_{  R} b)=
b_{-}\otimes_{{  R}^{op}} b_{+}$ are used, for $b\in \mathcal{B}$.

The anti Yetter-Drinfel'd module $X$ is said to be \emph{stable}
if in addition, for any $x\in X$,
\begin{equation}\label{eq:Y-D_stable}
x^{[1]}\triangleright x^{[0]} =x.
\end{equation}
\end{definition}

We need to show that condition (ii) in Definition
\ref{def:Y-D_mod} is meaningful, i.e. the expression
on the right hand side of (\ref{eq:Y-D_comp}) is well defined.
This follows by the following

\begin{lemma}\label{lem:Y-D_well_def}
For an algebra ${  R}$, consider a right $\times_{  R}$-Hopf
algebra $\mathcal{B}=({ B}, {\xi}, {\zeta}, {\Delta},
 {\epsilon})$. Let $X$ be a right $\mathcal{B}$-comodule
and left $\mathcal{B}$-module. Keeping the notations in Definition
\ref{def:Y-D_mod}, assume that axiom (i) in Definition
\ref{def:Y-D_mod} holds. Then the following hold.
\begin{enumerate}
\item Considering $\mathcal{B}$ as a left ${  R}$-module via
  $ {\zeta}$, the ${  R}$-module tensor product $X\otimes_{  R}
  \mathcal{B}$ is a left $\mathcal{B}$-module via the action
\[
b'\blacktriangleright (x\otimes_{  R} b):= b'_{+}\triangleright
x\otimes_{  R} bb'_{-}.
\]
\item $X\otimes_{  R}  \mathcal{B}$ is a
$\mathcal{B}$-$\mathcal{B}^{op}$
  bimodule, via the left $\mathcal{B}$-action in part (i) and the right
  $\mathcal{B}^{op}$-action
\[
(x\otimes_{  R} b)\blacktriangleleft b':= x\otimes_{  R} b'b.
\]
\item For any elements $x\in X$ and $r\in {  R}$,
\[
 {\xi}(r) \blacktriangleright \big(x^{[0]} \otimes_{ R}
x^{[1]}\big) =  \big(x^{[0]} \otimes_{  R}
x^{[1]}\big)\blacktriangleleft  {\zeta}(r).
\]
\end{enumerate}
\end{lemma}

\begin{proof}
By the first condition in (\ref{eq:Y-D_R_mod}) and Lemma
\ref{lem:x_R_Hopf_a} (ix), there is a well defined map
\[
\mathcal{B}\otimes_{\mathbb K} (X\otimes_{  R} \mathcal{B}) \to X\otimes_{
R} \mathcal{B}, \qquad b'\otimes_{\mathbb K} (x\otimes_{  R} b)\mapsto
b'_{+}\triangleright x \otimes_{  R} b\, b'_{-}.
\]
It is an associative left $\mathcal{B}$-action by part (iii) of
Lemma \ref{lem:x_R_Hopf_a} and it is unital by part (iv). This
proves claim (1). Claim (2) is obvious. Claim (3) follows by the
following computation, for $r\in {  R}$ and $x\in X$.
\begin{eqnarray*}
 {\xi}(r) \blacktriangleright \big(x^{[0]} \otimes_{ R}
x^{[1]}\big) &=&  {\xi}(r) \triangleright x^{[0]} \otimes_{  R}
x^{[1]} = x^{[0][0]}\cdot
 {\epsilon}( {\zeta}(r)
x^{[0][1]})\otimes_{  R} x^{[1]}  \\
&=&x^{[0]} \otimes_{  R}x^{[1](2)}
 {\zeta}\big( {\epsilon}( {\zeta}(r)
x^{[1](1)})\big) = x^{[0]} \otimes_{  R} {\zeta}(r) x^{[1]}\\
&=&\big(x^{[0]} \otimes_{  R} x^{[1]}\big)\blacktriangleleft
 {\zeta}(r).
\end{eqnarray*}
In the first equality we used that, by unitality and right ${
R}$-module map property of the coproduct, the map
(\ref{eq:x_R_Hopf_a}) satisfies $\vartheta(1_{\mathcal{B}}
\otimes_{  R^{op}}  {\xi}(r)) = 1_{\mathcal{B}} \otimes_{  R}
 {\xi}(r)$. Hence $ {\xi}(r)_{-}\otimes _{{
R}^{op}}  {\xi}(r)_{+} = 1_{\mathcal{B}} \otimes_{ R^{op}}
 {\xi}(r)$. The second equality follows by the second
condition in (\ref{eq:Y-D_R_mod}). The third equality follows by
coassociativity of the $\mathcal{B}$-coaction on $X$.
The fourth equality follows by counitality of $ {\Delta}$ and the
right bialgebroid version of (\ref{eq:left_R_mods}), i.e. the
identity
\begin{equation}\label{eq:R_mods}
 {\Delta}\big( {\xi}(r_1)
 {\zeta}(r_2) b  {\xi}(r_3)  {\zeta}
(r_4)\big)=  {\zeta}(r_2) b^{(1)} {\zeta} (r_4) \otimes_{  R}
 {\xi}(r_1) b^{(2)} {\xi}(r_3),
\end{equation}
for $r_1,r_2,r_3, r_4\in {  R}$ and $b\in \mathcal{B}$.
\end{proof}

Using the notations in Lemma \ref{lem:Y-D_well_def}, the right
hand side of (\ref{eq:Y-D_comp}) is equal to the well defined
expression $b^{(1)} \blacktriangleright \big(x^{[0]} \otimes_{ R}
x^{[1]}\big) \blacktriangleleft b^{(2)}$.

\begin{remark}
Note that a Hopf algebra $H$ over a commutative ring ${\mathbb K}$ is a
$\times_{\mathbb K}$-Hopf algebra. The bialgebroid structure is given by the
equal source and target maps ${\mathbb K}\to H$, $\kappa \mapsto \kappa 1_H$,
and the coproduct and counit in $H$. The canonical map \eqref{eq:x_R_Hopf_a}
has an inverse $\vartheta^{-1}(h' \otimes_{\mathbb K} h)= h' S(h^{(1)})
\otimes_{\mathbb K} h^{(2)}$, where $h \mapsto h^{(1)} \otimes_{\mathbb K}
h^{(2)}$ is the usual Sweedler index notation for the coproduct, with implicit
summation understood. 
{That is, $h_- \otimes_{\mathbb K} h_+= S(h^{(1)}) \otimes_{\mathbb K}
h^{(2)}$. Clearly, in this case  
\eqref{eq:Y-D_R_mod} 
becomes the trivial condition 
$$
x \kappa = (\kappa 1_H)\triangleright x = \kappa x, 
$$
for $\kappa \in {\mathbb K}$ and any element $x$ of a left $H$-module right
$H$-comodule $X$. This condition simply expresses the requirement that the
left and right ${\mathbb K}$-actions on $X$ are equal, and are induced by the 
$H$-module structure.
The second condition \eqref{eq:Y-D_comp} in Definition \ref{def:Y-D_mod}
reduces to  
$$
(b \triangleright x)^{[0]} \otimes (b\triangleright x)^{[1]}=
b^{(2)}\triangleright x^{[0]} \otimes b^{(3)} x^{[1]} S(b^{(1)}),
$$
for $b\in H$ and $x\in X$, which is the defining property of a (left-right)
anti Yetter-Drinfel'd module $X$ of a Hopf algebra $H$ in \cite[Definition
4.1]{JaSt:CycHom} or \cite[Definition 2.1]{HaKhRaSo1}. So we conclude that
Definition \ref{def:Y-D_mod} generalizes these definitions.}
\end{remark}

\begin{example}
For an algebra ${  R}$, consider a right $\times_{  R}$-Hopf algebra
 $\mathcal{B}=({ B}, {\xi}, {\zeta}, {\Delta},  {\epsilon})$. Note that left
 $B$-actions $\triangleright$ on $R$, satisfying $\zeta(r')\triangleright r
 =rr'$, are in bijective correspondence with maps $\chi:B\to R$, obeying the
 following properties, for $r\in R$, $b,b'\in B$.
\[
\chi(\zeta(r)b)=\chi(b)r\, \qquad
\chi(bb')=\chi\big(b\zeta(\chi(b'))\big), \qquad
\chi(1_B)=1_R.
\]
Indeed, in terms of such a map $\chi$, one can put $b\triangleright
r:=\chi(b\zeta(r))$. Furthermore, right $B$-coactions on $R$, with underlying
right regular $R$-module structure, are in bijective correspondence with
grouplike elements in $B$, i.e. $g\in B$ such that
$\Delta(g)=g\otimes_R g$ and
$\epsilon(g)=1_R$.
Indeed, in terms of a grouplike element $g$, a right $B$-coaction on $R$ is
given by
$r \mapsto 1_R \ot_R g\xi(r)$. One checks that the left $B$-module determined
by $\chi$ and the right $B$-comodule determined by $g$ combine into an anti
Yetter-Drinfel'd module on $R$ if and only if, for $r\in R$ and $b\in B$,
\[
\epsilon\big(\xi(r) g\big)=\chi\big(\xi(r)\big),\qquad \textrm{and}\qquad
g \xi\big(\chi(b)\big)  = b^{(2)} g {b^{(1)}}_-
\zeta\big(\chi({b^{(1)}}_+)\big).
\]
The anti Yetter-Drinfel'd module $R$ is stable if in addition $\chi \big(g
\xi(r)\big)= r$, for all $r\in R$. The pair $(\chi,g)$ with these properties
generalizes the notion of a {\em modular pair in involution} for a Hopf
algebra in \cite{ConMos:DiffCyc} or a weak Hopf algebra in \cite{Vecsi}.

\end{example}

\begin{proposition}\label{prop:dual_YD_proj}
Let $\mathcal{B}$ be a right $\times_{  R}$-Hopf algebra over an
algebra ${  R}$. Consider a right $\mathcal{B}$-module coring ${
C}$ with $\mathcal{B}$-action $\triangleleft$.
\begin{enumerate}
\item For an anti Yetter-Drinfel'd module $X$ of $\mathcal{B}$,
the
  para-cyclic object in Theorem \ref{thm:mod_coring}
  projects to a para-cyclic structure on ${
  C}^{\otimes_{  R}\, n+1}\otimes_{\mathcal{B}} X$.
\item For a stable anti Yetter-Drinfel'd module $X$ of $\mathcal{B}$, the
  para-cyclic object ${   C}^{\otimes_{  R}\, n+1}\otimes_{\mathcal{B}} X$ in
  part (1) is cyclic, in which case it will be denoted by $Z_\ast(C,M)$.
\end{enumerate}
\end{proposition}

\begin{proof}
We need to show that the composite map
\[
\xymatrix{
C^{\widehat{\otimes}_R\, n+1}\, \widehat{\otimes}_{  R} \,X
\ar[r]^{w_n}&
C^{\widehat{\otimes}_R\, n+1}\, \widehat{\otimes}_{  R} \,X
\ar@{>>}[r]& {  C}^{\otimes_{  R}\, n+1}{\otimes}_{\mathcal{B}}
\,X }
\]
is $\mathcal{B}$-balanced, i.e. that
\begin{equation*}
\big(c_n \triangleleft b^{(n+1)} x^{[1]} \stac{  R}
c_0\triangleleft b^{(1)} \stac {  R}  \dots \stac {  R}
c_{n-1}\triangleleft b^{(n)}\big) \stac {\mathcal{B}} x^{[0]} =
\big(c_n \triangleleft (b\triangleright x)^{[1]} \stac {  R} c_0
\stac {  R} \dots \stac {   R} c_{n-1} \big) \stac {\mathcal{B}}
(b\triangleright x)^{[0]},
\end{equation*}
for $b\in \mathcal{B}$, $x\in X$ and $c_0\,{\otimes}_{
R}\,\dots\,{\otimes}_{   R}\, c_{n}\in {  C}^{\otimes_{  R}\,
n+1}$. By counitality of the coproduct in $\mathcal{B}$, the left
hand side is equal to
\begin{align*}
c_n \triangleleft b^{(n+2)}& x^{[1]} \otimes_ {  R}
c_0\triangleleft
  b^{(2)} {\zeta}\big( {\epsilon}(b^{(1)})\big)
  \otimes_ {  R} c_1 \triangleleft b^{(3)}
  \otimes_ {  R} \dots \otimes_ {  R}c_{n-1}\triangleleft b^{(n+1)}\,
  \otimes_ {\mathcal{B}} x^{[0]} \\
&= c_n \triangleleft b^{(n+2)}
  x^{[1]} {\xi}\big( {\epsilon}(b^{(1)})\big)
  \otimes_ {  R} c_0\triangleleft b^{(2)}
  \otimes_ {  R} c_1 \triangleleft b^{(3)}
  \otimes_ {  R} \dots\otimes_ {  R} c_{n-1}\triangleleft b^{(n+1)}
  \otimes_ {\mathcal{B}} x^{[0]}\\
&= c_n \triangleleft b^{(n+2)} x^{[1]}{b^{(1)}}_{-}{b^{(1)}}_{+}
  \otimes_ {  R} c_0\triangleleft b^{(2)}
  \otimes_ {  R} c_1 \triangleleft b^{(3)}
  \otimes_ {  R} \dots \otimes_ {  R}c_{n-1}\triangleleft b^{(n+1)}
  \otimes_ {\mathcal{B}} x^{[0]}\\
&= c_n \triangleleft {b_{+}}^{(n+2)} x^{[1]}{b}_{-}{b_{+}}^{(1)}
  \otimes_ {  R} c_0\triangleleft  {b_{+}}^{(2)}
  \otimes_ {  R} \dots \otimes_ {  R}c_{n-1}\triangleleft {b_{+}}^{(n+1)}
  \otimes_ {\mathcal{B}} x^{[0]}\\
&= \big(c_n \triangleleft {b_{+}}^{(2)} x^{[1]}{b}_{-}
  \otimes_ {  R} c_0 \otimes_ {  R} \dots \otimes_ {  R}c_{n-1}
  \big)\triangleleft
  {b_{+}}^{(1)} \otimes_ {\mathcal{B}} x^{[0]}\\
&= \big(c_n \triangleleft {b_{+}}^{(2)} x^{[1]}{b}_{-}
  \otimes_ {  R} c_0 \otimes_ {  R} \dots \otimes_ {  R}c_{n-1} \big)
\otimes_ {\mathcal{B}}
{b_{+}}^{(1)} \triangleright x^{[0]}\\
&= \big(c_n \triangleleft (b\triangleright x)^{[1]} \otimes_ { R}
c_0 \otimes_ {  R} \dots \otimes_ {   R} c_{n-1} \big) \otimes_
{\mathcal{B}} (b\triangleright x)^{[0]}.
\end{align*}
The second equality follows by part (viii) in Lemma
\ref{lem:x_R_Hopf_a} and the third one follows by part (v). The
last equality is a consequence of Lemma \ref{lem:x_R_Hopf_a} (v)
and (\ref{eq:Y-D_comp}). This proves that the para-cyclic map
${w}_n$ in Theorem \ref{thm:mod_coring} factors to a map
$\widehat{w}_n:{
  C}^{\otimes_{  R}\,
  n+1}\otimes_\mathcal{B} X\to {  C}^{\otimes_{
    R}\,n+1}\otimes_\mathcal{B} X$, hence it defines a para-cyclic
structure on the simplex in Lemma \ref{lem:dual_B_mod_prod}. This
completes the proof of part (1). Furthermore,
$(\widehat{w}_n)^{n+1}$
takes an element $\big(c_0\otimes_{  R} \dots \otimes_{
  R} c_n\big)\otimes_{\mathcal{B}} x\in C^{\otimes_{  R}\,
  n+1}\otimes_{\mathcal{B}} X$ to
\[
\big(c_0\triangleleft x^{[1](1)}\stac {  R} \dots \stac {  R}
c_n\triangleleft x^{[1](n+1)} \big) \stac {\mathcal{B}} x^{[0]} =
\big(c_0 \stac {  R} \dots \stac{  R} c_n\big)\triangleleft
x^{[1]}  \stac {\mathcal{B}}  x^{[0]} = \big(c_0 \stac {  R} \dots
\stac {  R} c_n\big)\stac {\mathcal{B}} x^{[1]} \triangleright
x^{[0]}.
\]
Hence if $X$ is a stable anti Yetter-Drinfel'd module, i.e.
condition (\ref{eq:Y-D_stable}) holds, then
$(\widehat{w}_n)^{n+1}$
is the identity map. Thus we have claim (2) proven.
\end{proof}

\begin{remark}
Let $\mathcal{B}=({  B}, {\xi}, {\zeta},  {\Delta},
 {\epsilon})$ be a right bialgebroid over an algebra ${
R}$. Consider $\mathcal{B}$ as an ${  R}^{op}$-bimodule via right
multiplications by $ {\xi}$ and $ {\zeta}$. This bimodule has an
${  R}^{op}$-coring structure with counit $ {\epsilon}$ and
coproduct $ {\Delta}_{cop}:b \mapsto b^{(2)}\otimes_{{ R}^{op}}
b^{(1)}$, co-opposite to $ {\Delta}$. Together with the opposite
algebra ${ B}^{op}$, they constitute a left ${
R}^{op}$-bialgebroid $\mathcal{B}^{op}_{cop} =({ B}^{op}, {\xi},
{\zeta},  {\Delta}_{cop} ,  {\epsilon})$. A right
$\mathcal{B}$-module algebra ${ T}$ determines a left
$\mathcal{B}^{op}_{cop}$-module algebra ${ T}^{op}$, canonically.
Furthermore, a right $\mathcal{B}$-comodule $X$ can be looked at
as a left $\mathcal{B}^{op}_{cop}$-comodule. Application of
Theorem \ref{thm:mod_alg} to the left ${ R}^{op}$-bialgebroid
$\mathcal{B}^{op}_{cop}$, the left $\mathcal{B}^{op}_{cop}$-module
algebra ${  T}^{op}$ and the left
$\mathcal{B}^{op}_{cop}$-comodule $X$, yields a para-cocyclic
cosimplex {that has in degree $n$
\[
\big({T}^{op}\big)\, ^{\widehat{\otimes}_{{  R}^{op}}n+1}\,
\widehat{\otimes}_{{  R}^{op}}\, X \cong {  T}\,
^{\widehat{\otimes}_{  R}n+1} \, \widehat{\otimes}_{  R} \, X,
\]}
where the isomorphism is given by reversing the order of the
factors, i.e.
\[
t_0 \,\widehat{\otimes}_{{  R}^{op}} \, t_1 \widehat{\otimes}_{{
R}^{op}} \, \dots \, \widehat{\otimes}_{{  R}^{op}} \, t_n \,
\widehat{\otimes}_{{  R}^{op}} \, x \quad \mapsto  \quad t_n \,
\widehat{\otimes}_{{  R}} \, t_{n-1} \, \widehat{\otimes}_{{ R}}
\, \dots \, \widehat{\otimes}_{{  R}} \, t_0 \,
\widehat{\otimes}_{{
    R}} \, x.
\]
Resulting coface and codegeneracy maps on ${ T}\widehat{\otimes}_{
R} \dots \widehat{\otimes}_{  R} {  T}\widehat{\otimes}_{  R} X$
are the maps in (\ref{eq:bim_coface&codeg}) and the para-cocyclic
map comes out as
\begin{equation}\label{eq:op_mod_alg_hw}
t_0 \,\widehat{\otimes}_{  R}\,\dots\, \widehat{\otimes}_{
  R}\, t_{n}\, \widehat{\otimes}_{  R}\, x \quad \mapsto \quad
t_n \triangleleft x^{[1]}\, \widehat{\otimes}_{  R}\, t_0\,
  \widehat{\otimes}_{
  R}\,  \dots \, \widehat{\otimes}_{  R}\, t_{n-1}\, \widehat{\otimes}_{
  R}\, x^{[0]}.
\end{equation}
By the right $\mathcal{B}$-module map property of the maps $
{\varphi}:{  R}\to {  T}$, $r\mapsto 1_{ T}\triangleleft
 {\xi}(r)$, and the multiplication map $m_T:T\otimes_R T\to T$, the coface and
codegeneracy maps (\ref{eq:bim_coface&codeg}) project to ${
  T}^{\otimes_{  R}\, n+1}\otimes_{\mathcal{B}} X$. Note moreover that the
para-cocyclic map (\ref{eq:op_mod_alg_hw}) and the para-cyclic map
(\ref{eq:mod_coring_hw}) are of the same form. Hence it follows by
the computation in the proof of Proposition
\ref{prop:dual_YD_proj} that also the para-cocyclic map
(\ref{eq:op_mod_alg_hw}) projects to ${  T}^{\otimes_{
    R}\, n+1}\otimes_{\mathcal{B}} X$, whenever $X$ is an anti
Yetter-Drinfel'd module of $\mathcal{B}$. That is, in this case ${
  T}^{\otimes_{  R}\, \ast+1}\otimes_{\mathcal{B}} X$ is a para-cocyclic
object, which is cocyclic if the anti Yetter-Drinfel'd module $X$ is
stable.
\end{remark}

\subsection{Galois extensions of $\times_{  R}$-Hopf algebras}
\label{sec:Galois}

For a right $\times_R$-Hopf algebra $\mathcal{B}$ over an
algebra $R$, consider a right comodule algebra $T$ with coaction $t\mapsto
t^{[0]}\otimes_{  R} t^{[1]}$ (where implicit summation is
understood). The subalgebra $S$ of {\em coinvariants} in $T$
consists of those elements $s\in T$ for which
$s^{[0]}\otimes_{  R} s^{[1]}=s\otimes_ R 1_{\mathcal{B}}$. To
the inclusion map ${  S}\hookrightarrow {  T}$ one can associate a
cyclic simplex
\begin{equation}\label{eq:no_coef}
Z_\ast(T/S)=
T^{\widehat{\otimes}_S \,  \ast +1}
\end{equation}
as in Corollary \ref{ex:dual_no_coef}. On the other hand, it follows
by Proposition \ref{prop:dual_YD_proj} that, regarding
$\mathcal{B}$ as a right $\mathcal{B}$-module coring, for any
stable anti Yetter-Drinfel'd module $X$ of $\mathcal{B}$ there is
a cyclic simplex
\begin{equation}\label{eq:Gal_cosim}
Z_\ast(B,X)={\mathcal{B}}^{\otimes_{  R}\, \ast +1}\otimes_{\mathcal{B}}X,
\end{equation}
where ${\mathcal{B}}^{\otimes_{  R}\, n+1}$ is understood to be a
right $\mathcal{B}$-module via the diagonal action. In this
section, under the additional assumption that ${  T}$ is a
$\mathcal{B}$-\emph{Galois extension} of ${  S}$, we construct a
stable anti Yetter-Drinfel'd module $X:={  T}/\{\ st-ts\ |\ s\in {
S},\ t\in {  T}\ \}\cong S\ot_{S^e} T$ of $\mathcal{B}$, such that the cyclic
simplices (\ref{eq:no_coef}) and (\ref{eq:Gal_cosim}) are
isomorphic. This extends \cite[Theorem
  3.7]{JaSt:CycHom}.
\medskip

In a right comodule algebra ${  T}$ of a right bialgebroid
$\mathcal{B}$ over an algebra ${  R}$, we denote the coaction by
$\varrho:t\mapsto t^{[0]}\otimes_{  R} t^{[1]}$, where implicit
summation is understood. For the iterated power of the coaction we
write $(\varrho\otimes_{  R} {\mathcal{B}^{\otimes_{  R}n-1}})
\circ \dots \circ (\varrho\otimes_{ R} \mathcal{B})\circ
\varrho(t)=: t^{[0]}\otimes_{  R} \dots \otimes_{ R}
t^{[n-1]}\otimes_{  R} t^{[n]}$.

\begin{definition}\label{def:B_Gal}
Let $\mathcal{B}$ be a right bialgebroid over an algebra $R$.
A right $\mathcal{B}$-comodule algebra ${T}$ is said to be a
$\mathcal{B}$-Galois extension of its coinvariant subalgebra $S$ if
the canonical map
\begin{equation}\label{eq:B_can}
\mathrm{can}: {  T}\otimes_{  S} {  T} \to {  T}\otimes_{  R}
\mathcal{B},\qquad t'\otimes_{  S} t\mapsto t't^{[0]}\otimes_{ R}
t^{[1]}
\end{equation}
is bijective.
\end{definition}
For example, if $\mathcal{B}=({  B},  {\xi},  {\zeta},
 {\Delta},  {\epsilon})$ is a right $\times_{
R}$-Hopf algebra, then the right regular $\mathcal{B}$-comodule
algebra is a $\mathcal{B}$-Galois extension of the coinvariant
subalgebra $ {\zeta}({  R})\cong { R}^{op}$.

Let $\mathcal{B}$ be a right ${  R}$-bialgebroid and ${  T}$ a
right $\mathcal{B}$-comodule algebra. It follows by the right ${
R}$-module map property of a right $\mathcal{B}$-coaction and
(\ref{eq:com_alg}) that, for a coinvariant $s \in {  S}$ and $r\in
{  R}$,
\[
((1_{T}\cdot r)s)^{[0]}\otimes_{R} ((1_{T}\cdot r)s)^{[1]} = s \otimes_{R}
 {\xi}(r) = (s(1_{T}\cdot r))^{[0]}\otimes_{R} (s( 1_{T}\cdot r))^{[1]}.
\]
Hence, applying the counit of $\mathcal{B}$ to the second factor
on both sides, we conclude that the elements $s\in {  S}$ commute
with $1_{T}\cdot r$, for all $r\in {  R}$. Hence ${T}\otimes_{R} \mathcal{B}$
is a right ${  S}$-module, with action
\[
(t\otimes_{  R} b)\cdot s := ts\otimes_{  R} b.
\]

Consider a right $\times_{  R}$-Hopf algebra $\mathcal{B}$ over an
algebra
  ${  R}$, and a $\mathcal{B}$-Galois extension ${  S}\subseteq {  T}$.
As in Lemma \ref{lem:x_R_Hopf_a}, in terms of the maps
(\ref{eq:x_R_Hopf_a}) and (\ref{eq:B_can}), introduce the
  index notations
\begin{equation}\label{eq:index_notation}
\mathrm{can}^{-1}(1_{  T}\otimes_{  R} b)=: b^{\{-\}}\otimes_{ S}
b^{\{+ \}}  \qquad \textrm{and}\qquad
\vartheta^{-1}(1_{\mathcal{B}}\otimes_{  R} b)=: b_{-}\otimes_{{
R}^{op}} b_{+},
\end{equation}
for $b\in \mathcal{B}$, where in both cases implicit summation is
understood. The following lemma is a right bialgebroid version of
\cite[Lemma 4.1.21]{Hobst:PhD}. It extends Lemma
\ref{lem:x_R_Hopf_a} (vi).
\begin{lemma}\label{lem:pentagon}
Consider a right $\times_{  R}$-Hopf algebra $\mathcal{B}=({
  B}, {\xi}, {\zeta}, {\Delta},
   {\epsilon})$ over an algebra
  ${  R}$, and a $\mathcal{B}$-Galois extension ${  S}\subseteq {  T}$.
Using the notations in (\ref{eq:index_notation}), for any $b\in
\mathcal{B}$ the following pentagonal equation holds in
  $({  T}\otimes_{  R} \mathcal{B})\otimes_{  S} {  T}$.
\[
b^{\{-\}[0]} \otimes_{  R} b^{\{-\}[1]} \otimes_{  S} b^{\{+\}} =
{b_{+}}^{\{-\}}\otimes_{  R} b_{-} \otimes_{  S} {b_{+}}^{\{+\}}.
\]
\end{lemma}
\begin{proof}
The second condition in (\ref{eq:right_coac_mod}) implies that there is a well
defined bijection
\[
\mathrm{can}_{13}: ({  T}\otimes_{  R} \mathcal{B})\otimes_{  S} {
  T}\to {  T}\otimes_{  R} (\mathcal{B}\otimes_{{  R}^{op}}
  \mathcal{B}),\qquad
(t'\otimes_{  R} b)\otimes_{  S} t \mapsto t' t^{[0]}\otimes_{
R}(b
  \otimes_{{  R}^{op}} t^{[1]}),
\]
where in the ${  R}^{op}$-module tensor product
$\mathcal{B}\otimes_{{  R}^{op}} \mathcal{B}$ the ${
R}^{op}$-actions ${  R}^{op}\otimes_{\mathbb K} \mathcal{B}\otimes_{\mathbb K}
{R}^{op}\to \mathcal{B}$, $r_1 \otimes_{\mathbb K} b \otimes_{\mathbb K}
r_2\mapsto {\zeta}(r_2)b {\zeta}(r_1)$ are used, and
$\mathcal{B}\otimes_{{  R}^{op}} \mathcal{B}$ is meant to be a
left ${
  R}$-module via the action $r\cdot(b\otimes_{{  R}^{op}} b')=b\otimes_{{
    R}^{op}} b'  {\zeta}(r)$. A straightforward computation using
(\ref{eq:com_alg}) shows that the ${  S}$-bimodule map and right
${
  R}$-module map (\ref{eq:B_can}), and the left ${  R}$-module map
(\ref{eq:x_R_Hopf_a}) satisfy
\[
(\mathrm{can}\otimes_{  R} \mathcal{B})\circ ({  T}\otimes_{  S}
\mathrm{can})= ({  T}\otimes_{  R} \vartheta)\circ
\mathrm{can}_{13}\circ (\mathrm{can}\otimes_{  S} {  T}).
\]
Hence, by bijectivity of all involved maps,
\[
(\mathrm{can}\otimes_{  S} {  T})\circ ({  T}\otimes_{  S}
\mathrm{can}^{-1}) = \mathrm{can}_{13}^{-1} \circ ({ T}\otimes_{
R} \vartheta^{-1}) \circ (\mathrm{can}\otimes_{  R} \mathcal{B}).
\]
Application of this identity to $1_{  T}\otimes_{  S} 1_{  T}
\otimes_{  R} b$ yields the claim.
\end{proof}

For an algebra extension $S\subseteq T$, $T$ has a canonical $S$-bimodule
structure. Hence application of the functor
$\bs{\Pi}:S\text{-}\mathrm{Mod}\text{-}S \to \mathrm{Mod}\text{-}\mathbb{K}$ in
Definition \ref{def:Pi} to $T$ yields a $\mathbb{K}$-module $\bs{\Pi} T \cong
S\widehat{\otimes}_S\, T$.

\begin{proposition}\label{prop:T_S}
Consider a right $\times_{  R}$-Hopf algebra $\mathcal{B}$ over an
algebra
  ${  R}$, and a $\mathcal{B}$-Galois extension ${  S}\subseteq {  T}$.
Then the quotient
\begin{equation}\label{eq:T_S}
T_S:= S\widehat{\otimes}_S\, T
\end{equation}
is a stable anti Yetter-Drinfel'd module.
\end{proposition}

\begin{proof}
Since the ${  S}$-, and ${  R}$-actions on ${  T}$ commute (cf.
second paragraph following Definition \ref{def:B_Gal}), there is a
unique ${
  R}$-bimodule structure on ${  T}_{  S}$ such that the epimorphism
$p_{  T}:{  T}\to {  T}_{  S}$ is an ${  R}$-bimodule map.
Furthermore, by the ${  S}$-bimodule map property of the coaction
$\varrho:t\mapsto t^{[0]}\otimes_{  R} t^{[1]}$ in ${  T}$, the
map $(p_{  T}\otimes_{  R} \mathcal{B})\circ \varrho:{  T}\to {
T}_{
  S}\otimes_{  R} \mathcal{B} $ coequalizes the left and right ${
  S}$-actions on ${  T}$. Hence there exists a unique right
$\mathcal{B}$-comodule structure on ${  T}_{  S}$ such that $p_{
  T}:{  T}\to
{  T}_{  S}$ is a right $\mathcal{B}$-comodule map.

The algebra map ${  S}\hookrightarrow {  T}$ equips ${  T}$ with
an $S$-bimodule structure. The center $(T\otimes_S T)^S$ of the $S$-bimodule
$T\otimes_S T$ is
  an algebra, with multiplication $(\sum_i u_i\otimes_{  S} u'_i)(\sum_j
  v_j\otimes_{  S} v'_j)=\sum_{i,j} v_j u_i \otimes_{  S} u'_i v'_j$.
Recall from \cite[Section 2.2]{JaSt:CycHom} that for any any ${
  T}$-bimodule $M$, the quotient $S\widehat{\otimes}_S\, M\cong
  M/\{\ s\cdot m-m\cdot s\ |\  s\in {  S},\ m\in M\ \}$ is a left
  $({  T}\otimes_{  S} {  T})^{  S}$-module via the action
\[
({  T}\otimes_{  S} {  T})^{  S} \otimes_{\mathbb K} M_{  S}\to M_{  S},
\qquad (\sum_i u_i\otimes_{  S} u'_i) \otimes_{\mathbb K} p_M (m) \mapsto
\sum_i p_M (u'_i m u_i),
\]
where $p_M:M \to S\widehat{\otimes}_S\, M$ denotes the canonical
epimorphism. In particular, ${  T}_{  S}$ is a left $({  T}\otimes_{  S}
  {  T})^{  S}$-module.

On the other hand, for a Galois extension ${  S}\subseteq {  T}$
by a right ${  R}$-bialgebroid $\mathcal{B}$, using the notation in
  (\ref{eq:index_notation}), the map
\begin{equation}\label{eq:transl_map}
\mathcal{B} \to ({   T}\otimes_{  S} {  T})^{  S},\qquad b \mapsto
b^{\{-\}}\otimes_{  S} b^{\{+\}}
\end{equation}
is an algebra homomorphism. Indeed, by the ${  S}$-bimodule map
  property of the coaction on ${  T}$ it follows that
\[
\mathrm{can}(sb^{\{-\}}\otimes_{  S} b^{\{+\}})=s \otimes_{  R} b
= \mathrm{can}(b^{\{-\}}\otimes_{  S} b^{\{+\}}s).
\]
Hence, by bijectivity of $\mathrm{can}$, $b^{\{-\}}\otimes_{  S}
  b^{\{+\}}\in ({   T}\otimes_{  S} {  T})^{  S}$, for all $b \in
  \mathcal{B}$. The map (\ref{eq:transl_map}) is unital and multiplicative
  since by (\ref{eq:com_alg}) for all $b,b'\in \mathcal{B}$,
\[
\mathrm{can}(1_{  T}\otimes_{  S}1_{  T})=1_{  T}\otimes_{  S}
  1_{\mathcal{B}}\qquad \textrm{and}\qquad
\mathrm{can}(b^{\prime \{-\}}b^{\{-\}}\otimes_{  S} b^{\{+\}}
b^{\prime
  \{+\}}) = 1_{  T}\otimes_{  R} bb'
\]
and $\mathrm{can}$ is bijective. This
  proves that ${  T}_{  S}$ is a left $\mathcal{B}$-module, with a so called
  Miyashita-Ulbrich type action
\begin{equation}\label{eq:T_S_mod}
b\triangleright p_{  T}(t)=p_{  T}(b^{\{+\}} t b^{\{-\}}).
\end{equation}
It remains to check the compatibility conditions in Definition
\ref{def:Y-D_mod} between the $\mathcal{B}$-module and
$\mathcal{B}$-comodule structures on ${  T}_{  S}$. It follows by
the ${  R}$-bimodule map property of (\ref{eq:B_can}) that
\begin{equation}\label{eq:xi-zeta}
 {\zeta}(r)^{\{-\}}\otimes_{
  S}  {\zeta}(r)^{\{+\}} = r\cdot 1_{  T}\otimes_{  S} 1_{
  T}\qquad \textrm{and}\qquad {\xi}(r)^{\{-\}}\otimes_{
  S}  {\xi}(r)^{\{+\}} = 1_{  T}\otimes_{  S} 1_{
  T}\cdot r.
\end{equation}
Hence
\[
 {\zeta}(r) \triangleright p_{  T}(t)=p_{
T}\big(t(r\cdot 1_{  T})\big) = p_{  T}(t\cdot r)=p_{ T}(t)\cdot r
\ \  \textrm{and}\ \  {\xi}(r) \triangleright p_{ T}(t)=p_{
T}\big((1_{  T}\cdot r)t\big) = p_{  T}(r\cdot t)=r\cdot p_{
T}(t).
\]
Furthermore, for $b\in \mathcal{B}$ and $t\in {  T}$,
\begin{align*}
\left(b\triangleright p_{  T}(t)\right)^{[0]}\otimes_{  R}
\left(b\triangleright p_{  T}(t)\right)^{[1]} &= p_{
T}\left(b^{\{+\}[0]} t^{[0]} b^{\{-\}[0]} \right) \otimes_{  R}
b^{\{+\}[1]} t^{[1]} b^{\{-\}[1]} \\
&=p_{  T}\left(b^{(1)\{+\}} t^{[0]} b^{(1)\{-\}[0]} \right)
\otimes_{  R}
b^{(2)} t^{[1]} b^{(1)\{-\}[1]} \\
&= p_{  T}\left({{b^{(1)}}_{+}}^{\{+\}} t^{[0]}
{{b^{(1)}}_{+}}^{\{-\}}\right) \otimes_{  R} b^{(2)} t^{[1]} {b^{(1)}}_{-}\\
&={{b^{(1)}}_{+}} \triangleright p_{  T}(t)^{[0]}\otimes_{  R}
b^{(2)} p_{  T}(t)^{[1]} {b^{(1)}}_{-}\ .
\end{align*}
The first equality follows by (\ref{eq:T_S_mod}), the second
condition in (\ref{eq:com_alg}) and the comodule map property of
$p_{  T}$. The second equality is a consequence of the right
$\mathcal{B}$-comodule map property of
(\ref{eq:B_can}) hence of the map (\ref{eq:transl_map}). The third
equality is resulted by the application of Lemma
\ref{lem:pentagon}. The last equality follows by
(\ref{eq:T_S_mod}) and the comodule map property of $p_{  T}$
again. Thus we proved that ${  T}_{  S}$ is an anti
Yetter-Drinfel'd module. Finally, by the comodule map property of
$p_{  T}$ and the identity
\begin{equation}\label{eq:can_ii}
t^{[0]}t^{[1]\{-\}}\otimes_{  S} t^{[1]\{+\}} =
\mathrm{can}^{-1}\big( \mathrm{can}(1_{  T}\otimes_{  S}
t)\big)=1_{
    T}\otimes_{  S} t, \qquad \textrm {for } t\in {  T},
\end{equation}
it follows that
\[
p_{  T}(t)^{[1]}\triangleright p_{  T}(t)^{[0]}=
t^{[1]}\triangleright p_{  T}(t)^{[0]}= p_{
T}\big(t^{[1]\{+\}}t^{[0]}t^{[1]\{-\}}\big)=p_{  T}(t).
\]
That is, the anti Yetter-Drinfel'd module ${  T}_{  S}$ is stable.
\end{proof}

\begin{lemma}\label{lem:alpha}
For an algebra ${  R}$ over a commutative ring $\mathbb{K}$, consider a
right ${  R}$-bialgebroid $\mathcal{B}$ and a $\mathcal{B}$-Galois
extension ${  S}\subseteq {  T}$. Using the notation in
(\ref{eq:index_notation}), for any non-negative integer
  $n$ there exist ${  S}$-bimodule isomorphisms ${  T}^{\otimes_{  S}\,
  n+1}\cong {  T}\otimes_{  R}  \mathcal{B}^{\otimes_{  R}\, n}$,
\begin{align*}
&\alpha_n(
t_0\stac{  S} \dots \stac{  S}  t_n )=t_0 t_1^{[0]} t_2^{[0]}\dots
t_n^{[0]}\stac{  R} t_1^{[1]} t_2^{[1]}\dots t_n^{[1]}\stac{  R}
t_2^{[2]}t_3^{[2]}\dots t_n^{[2]}\stac{  R} \dots \stac{  R}
t_{n-1}^{[n-1]} t_{n}^{[n-1]} \stac{  R}  t_n^{[n]},\\
&\alpha_n^{-1}(
t\,\stac{  R}\,  b_1\,\stac{  R}\, \dots \,\stac{  R}\,  b_n)=
t\, {b_1}^{\{-\}} \,\stac{  S} {b_1}^{\{+\}} {b_2}^{\{-\}}
\,\,\stac{  S} {b_2}^{\{+\}} {b_3}^{\{-\}} \, \stac{
  S}  \dots \,\stac{  S} {b_{n-1}}^{\{+\}} {b_n}^{\{-\}} \,\stac{  S}
        {b_n}^{\{+\}}.
\end{align*}
Projections of the above isomorphisms yield $\mathbb{K}$-module
isomorphisms
$\widehat{\alpha}_n:{  T}^{\widehat{\otimes}_{  S}\, n+1}\to {
  T}_{  S}  \otimes_{  R}  \mathcal{B}^{\otimes_{  R}\, n}$,
where ${  T}_{  S}$ is the ${  R}$-bimodule (\ref{eq:T_S}).
\end{lemma}

\begin{proof}
It follows by the $S$-bimodule map property of the right
$\mathcal{B}$-coaction on $T$ that $\alpha_n$ is a well  defined $S$-bimodule
map. The to-be-inverse $\alpha_n^{-1}$ is well defined by
(\ref{eq:xi-zeta}) and Lemma \ref{lem:x_R_Hopf_a} (iii). We prove
by induction that the maps $\alpha_n$ and $\alpha_n^{-1}$ are
mutual inverses. For $n=0$, both $\alpha_0$ and $\alpha_0^{-1}$
are equal to the identity map on $T$, hence they are mutual
inverses. It follows by the second condition in (\ref{eq:com_alg})
that, for all values of $n$,
\[
\alpha_{n+1}=(\alpha_n\otimes_R \mathcal{B})\circ (
T^{\otimes_S\, n}\otimes_S \mathrm{can})\qquad
\textrm{and}\qquad \alpha_{n+1}^{-1}= (T^{\otimes_S\,
 n}\otimes_S \mathrm{can}^{-1}) \circ (\alpha_n^{-1}\otimes_R \mathcal{B}).
\]
Hence if $\alpha_n^{-1}$ is the inverse of $\alpha_n$ then
$\alpha_{n+1}$ is also an ${  S}$-bimodule isomorphism with
inverse $\alpha_{n+1}^{-1}$. Applying the functor in Definition
\ref{def:Pi}, from the category of ${  S}$-bimodules to the
category of $\mathbb{K}$-modules, it takes $\alpha_n$ to the required
$\mathbb{K}$-module isomorphism $\widehat{\alpha}_n:{
  T}^{\widehat{\otimes}_{  S}\, n+1}\to {  T}_{  S}
  \otimes_{  R}  \mathcal{B}^{\otimes_{  R}\, n}$.
\end{proof}

\begin{lemma}\label{lem:beta}
Let $\mathcal{B}$ be a right $\times_{  R}$-Hopf algebra over an
algebra
  ${  R}$.
For the inverse of the canonical map (\ref{eq:x_R_Hopf_a}) use the
index notation in (\ref{eq:index_notation}). Then, for any
non-negative integer $n$, there exist right $\mathcal{B}$-module
isomorphisms $\mathcal{B}^{\otimes_{
    R}\ n}\otimes_{{  R}^{op}} \mathcal{B}\cong \mathcal{B}^{\otimes_{
    R}\ n+1}$,
\begin{align*}
\beta_n \big( b_1\otimes_{  R}\dots \otimes_{  R} b_n\otimes_{{
R}^{op}} b'\big)&=  b_1 b^{\prime (1)}\otimes_{  R}\dots \otimes_{
R} b_n
b^{\prime    (n)}\otimes_{  R}  b^{\prime (n+1)}\\
\beta_n^{-1} \big(  b_1\otimes_{  R}\dots \otimes_{  R}
b_{n}\otimes_{ R} b'\big)&= b_1 {b'_{-}}^{(1)} \otimes_{ R}\dots
\otimes_{  R} b_n
  {b'_{-}}^{(n)} \otimes_{{  R}^{op}} b'_{+}.
\end{align*}
\end{lemma}

\begin{proof}
For any right $\mathcal{B}$-module $N$, the $\mathcal{B}$-bimodule
isomorphism (\ref{eq:x_R_Hopf_a}) induces a right
$\mathcal{B}$-module isomorphism $N\otimes_{\mathcal{B}}
\vartheta:N\otimes_{\mathcal{B}} \mathcal{B}\otimes_{{  R}^{op}}
\mathcal{B} \cong N \otimes_{{  R}^{op}} \mathcal{B} \to
N\otimes_{\mathcal{B}} \mathcal{B}\otimes_{{  R}} \mathcal{B}
\cong N \otimes_{{  R}} \mathcal{B}$. Consider
$\mathcal{B}^{\otimes_{  R}\, n}$ as a right $\mathcal{B}$-module
via the diagonal action. Then $\beta_n = \mathcal{B}^{\otimes_{
R}\, n} \otimes_{\mathcal{B}} \vartheta$ is a right
$\mathcal{B}$-module isomorphism as stated.
\end{proof}
\begin{theorem}\label{thm:main_iso}
Let $\mathcal{B}$ be a right $\times_{  R}$-Hopf algebra over an
  algebra ${  R}$ and let ${  S}\subseteq {  T}$ be a
  $\mathcal{B}$-Galois extension. Consider the right regular
  $\mathcal{B}$-module coring and the stable anti Yetter-Drinfel'd
  module ${  T}_{  S}$ in Proposition \ref{prop:T_S}. Then the associated
  cyclic simplex $Z_\ast(B,T_S)$ in Proposition \ref{prop:dual_YD_proj}
is isomorphic to the cyclic simplex $Z_\ast(T/S)$ in Corollary
\ref{ex:dual_no_coef}.
\end{theorem}
\begin{proof}
In terms of the maps in Lemmata \ref{lem:alpha} and
\ref{lem:beta}, for any non-negative integer $n$, one constructs a
$\mathbb{K}$-module isomorphism $\omega_n$ as the composition of the
following morphisms
$$
\xymatrix{
T^{\widehat{\otimes}_{  S}\, n+1} \ar[r]^-{\widehat{\alpha}_n} &
T_S\stac R \mathcal{B}^{\otimes_{  R}\, n} \ar[r]&
\mathcal{B}^{\otimes_R\, n}\stac {R^{op}} T_S\ar[r]^-{\cong} &
\mathcal{B}^{\otimes_R\, n}\stac {R^{op}} \mathcal{B}
  \stac {\mathcal{B}} T_S\ar[r]^-{\beta_n \otimes_{\mathcal{B}} T_S} &
\mathcal{B}^{\otimes_{  R}\, n+1}\stac {\mathcal{B}} {  T}_{  S}
}.
$$
Explicitly,
\begin{equation}\label{eq:main_iso}
\omega_n\big( t_0\,\cten{  S}\, t_1 \,\cten{  S}\,\dots \,\cten{
S}\, t_n\big)= \big(t_1^{[1]}t_2^{[1]}\dots t_n^{[1]}\stac {
R}t_2^{[2]}\dots t_n^{[2]}\stac {  R}\dots \stac {  R}
t_n^{[n]}\stac {  R} 1_{\mathcal{B}} \big)\stac {\mathcal{B}} p_{
T}(t_0 t_1^{[0]}t_2^{[0]} \dots t_n^{[0]}),
\end{equation}
where $p_{  T}:{  T}\to {  T}_{  S}\cong S \widehat{\otimes}_S\, T$ denotes
the canonical epimorphism. We show that (\ref{eq:main_iso}) is a homomorphism
of cyclic simplices.

Denote the counit in $\mathcal{B}$ by $ {\epsilon}$. By
multiplicativity of the coproduct in $\mathcal{B}$, for any
integer $0\leq k<n$,
\begin{align*}
\big((\mathcal{B}&^{\otimes_{  R}\, k}
  \stac {  R}  {\epsilon} \stac {  R}\mathcal{B}^{\otimes_{
  R}\, n-k})
\stac {\mathcal{B}} {  T}_{  S}\big)\circ \omega_n (t_o\, \cten {
  S}\, \dots \, \cten {  S}\, t_n)=\\
&= (t_1^{[1]}\dots t_n^{[1]}
\stac {  R} \dots \stac {  R} t_k^{[k]}\dots t_n^{[k]}
 \stac {  R} t_{k+2}^{[k+1]}\dots t_n^{[k+1]}\stac {  R} \dots \stac {
  R} t_n^{[n-1]}\stac {  R} 1_{\mathcal{B}} ) \stac {\mathcal{B}} p_{
  T}(t_0 t_1^{[0]}\dots t_n^{[0]})\\
&= \omega_{n-1} (t_o\, \cten {  S}\, \dots \, \cten {  S}\,
t_{k-1} \cten
  {  S} t_k t_{k+1} \cten {  S} t_{k+2}\cten {  S}\dots \cten{  S}
  t_n).
\end{align*}
Furthermore, using the form of the (diagonal) right
$\mathcal{B}$-action on $\mathcal{B}^{\otimes_{  R}\, n}$ (in the
second equality) and the form of the left $\mathcal{B}$-action
(\ref{eq:T_S_mod}) on ${  T}_{  S}$ (in the third equality), one
computes
\begin{align}
\big((\mathcal{B}^{\otimes_{  R}\, n}&\stac {  R}
 {\epsilon})
  \stac {\mathcal{B}} {  T}_{  S}\big)\circ \omega_n (t_o\, \cten {
  S}\, \dots \, \cten {  S}\, t_n)=\nonumber\\
&=\big(t_1^{[1]}\dots t_n^{[1]}\stac {  R}t_2^{[2]}\dots t_n^{[2]}
\stac {  R} \dots \stac {  R} t_n^{[n]}\big)\stac {\mathcal{B}}
p_{
  T}(t_0 t_1^{[0]}\dots t_n^{[0]})\nonumber\\
&= \big(t_1^{[1]}\dots t_{n-1}^{[1]}\stac {  R}t_2^{[2]}\dots
t_{n-1}^{[2]} \stac {  R} \dots \stac {  R} t_{n-1}^{[n-1]}\stac {
R}
  1_{\mathcal{B}}\big) \triangleleft t_n^{[1]} \stac {\mathcal{B}} p_{
  T}(t_0 t_1^{[0]}\dots t_n^{[0]})\nonumber\\
&= \big(t_1^{[1]}\dots t_{n-1}^{[1]}\stac {  R}t_2^{[2]}\dots
t_{n-1}^{[2]} \stac {  R} \dots \stac {  R} t_{n-1}^{[n-1]}\stac {
R}
  1_{\mathcal{B}}\big) \stac {\mathcal{B}} p_{
  T}(t_n^{[1]\{+\}} t_0 t_1^{[0]}\dots t_n^{[0]}t_n^{[1]\{-\}})\nonumber\\
&= \big(t_1^{[1]}\dots t_{n-1}^{[1]}\stac {  R}t_2^{[2]}\dots
t_{n-1}^{[2]} \stac {  R} \dots \stac {  R} t_{n-1}^{[n-1]}\stac {
R}
  1_{\mathcal{B}}\big) \stac {\mathcal{B}} p_{
  T}(t_n t_0 t_1^{[0]}\dots t_{n-1}^{[0]})\nonumber\\
&= \omega_{n-1} (t_n t_o\, \cten{  S} t_1 \, \cten {  S}\, \dots
\, \cten
  {  S}\,  t_{n-1}),
\label{eq:k=n}
\end{align}
where the penultimate equality follows by (\ref{eq:can_ii}).
This proves that the map $\omega_n$ is compatible with the face maps.

Denote the coproduct in $\mathcal{B}$ by $ {\Delta}$. It follows
by its multiplicativity that, for any integer $0\leq k<n$,
\begin{align*}
\big((\mathcal{B}&^{\otimes_{  R}\, k}
  \stac {  R}  {\Delta} \stac {  R}\mathcal{B}^{\otimes_{
  R}\, n-k})
\stac {\mathcal{B}} {  T}_{  S}\big)\circ \omega_n (t_o\, \cten {
  S}\, \dots \, \cten {  S}\, t_n)=\\
&= (t_1^{[1]}\dots t_n^{[1]}
\stac {  R} \dots \stac {  R} t_{k+1}^{[k+1]}\dots t_n^{[k+1]}
 \stac {  R} t_{k+1}^{[k+2]}\dots t_n^{[k+2]}\stac {  R} \dots \stac {
  R} t_n^{[n+1]}\stac {  R} 1_{\mathcal{B}} ) \stac {\mathcal{B}} p_{
  T}(t_0 t_1^{[0]}\dots t_n^{[0]})\\
&= \omega_{n+1} (t_o\, \cten {  S}\, \dots \, \cten {  S}\, t_{k}
\cten
  {  S} 1_{  T} \cten{  S} t_{k+1} \cten {  S}\dots \cten{  S}
  t_n).
\end{align*}
Furthermore, by unitality of $ {\Delta}$ and of the
$\mathcal{B}$-coaction on ${  T}$,
\begin{align*}
\big((\mathcal{B}&^{\otimes_{  R}\, n}\stac {  R}
 {\Delta})
  \stac {\mathcal{B}} {  T}_{  S}\big)\circ \omega_n (t_o\, \cten {
  S}\, \dots \, \cten {  S}\, t_n)=\\
&=\big(t_1^{[1]}t_2^{[1]}\dots t_n^{[1]}\stac { R}t_2^{[2]}\dots
t_n^{[2]}\stac {  R}\dots \stac {  R} t_n^{[n]}\stac {  R}
1_{\mathcal{B}} \stac {  R} 1_{\mathcal{B}} \big)\stac
{\mathcal{B}} p_{
  T}(t_0 t_1^{[0]}t_2^{[0]} \dots t_n^{[0]}) \\
&= \omega_{n+1} (t_o\, \cten {  S}\, \dots \, \cten {  S}\, t_n
\cten
  {  S} 1_{  T}).
\end{align*}
This proves that the map $\omega_n$ is compatible with the
degeneracy maps.

Finally, by similar steps as in (\ref{eq:k=n}), the cyclic map
$\widehat{w}_n$ in Proposition \ref{prop:dual_YD_proj} is checked
to satisfy
\begin{align*}
\widehat{w}_n\circ &\omega_n (t_o\, \cten {
  S}\, \dots \, \cten {  S}\, t_n)=\\
&=\big(t_0^{[1]}t_1^{[1]}t_2^{[1]}\dots t_n^{[1]}\stac {
R}t_1^{[2]}\dots t_n^{[2]}\stac {  R}\dots \stac {  R}
t_{n-1}^{[n]}t_n^{[n]}\stac {  R}
 t_n^{[n+1]} \big)\stac {\mathcal{B}} p_{
  T}(t_0^{[0]} t_1^{[0]}t_2^{[0]} \dots t_n^{[0]}) \\
&=\big(t_0^{[1]}t_1^{[1]}t_2^{[1]}\dots t_{n-1}^{[1]}\stac {
  R}t_1^{[2]}\dots t_{n-1}^{[2]}\stac {  R}\dots \stac {  R} t_{n-1}^{[n]}
\stac {  R} 1_{\mathcal{B}} \big) \stac {\mathcal{B}} p_{ T}(t_n
t_0^{[0]} t_1^{[0]}t_2^{[0]} \dots
  t_{n-1}^{[0]}) \\
&=\omega_n (t_n \, \cten {  S}\, t_0\, \cten {  S}\, \dots \,
\cten {
  S}\, t_{n-1}).
\end{align*}
Hence $\omega_n$ is compatible with the cyclic maps as well, what
proves the claim.
\end{proof}

\section{Cyclic homology of groupoids}

Consider a groupoid $\mathcal G$ (i.e. a small category in which all morphisms
are invertible) with a {\em finite} set ${\mathcal G}^0$ of objects and an
arbitrary set ${\mathcal G}^1$ of morphisms. Via the map associating to
$x\in {\mathcal G}^0$ the identity morphism $x\to x$, we consider ${\mathcal
G}^0$ as a subset of ${\mathcal G}^1$. Composition in ${\mathcal G}$ is
denoted by $\circ$ while the source and target maps ${\mathcal G}^1 \to
{\mathcal G}^0$ are denoted by $s$ and $t$, respectively.
For any field $\mathbb K$, the $\mathbb K$-vector space $B:= {\mathbb K}
{\mathcal G}^1$, spanned by the elements of ${\mathcal G}^1$, has a right
$\times_R$-Hopf algebra structure over the commutative base algebra
$R:={\mathbb K} {\mathcal G}^0$. Structure maps are the following.
Multiplication in the ${\mathbb K}$-algebra $B$ is given on basis elements
$g,g'\in {\mathcal G}^1$ by $g\circ g'$, if $g$ and $g'$ are composable,
i.e. $s(g)=t(g')$, and zero otherwise. We denote by juxtaposition this
product, linearly extended to all elements of $B$. Unit element is
$1_B=\sum_{x\in
{\mathcal G}^0} x$. Similarly, $R$ is a commutative ${\mathbb K}$-algebra with
minimal orthogonal idempotents $\{\ x\in {\mathcal G}^0\ \}$. Both
algebra maps $\xi$ and $\zeta:R \to B$ are induced by the inclusion map
${\mathcal G}^0 \hookrightarrow {\mathcal G}^1$. That is, $B$ is an $R$-module
(or $R$-bimodule, with coinciding left and right actions) via multiplication
on the right. Coproduct is diagonal on the basis elements $g\in {\mathcal
G}^1$, i.e. $\Delta(g):=g\otimes_R g$. Counit $\epsilon$ maps $g\in {\mathcal
G}^1$ to $s(g)\in {\mathcal G}^0$. The canonical map \eqref{eq:x_R_Hopf_a} has
the explicit form on the generating set $\{\ g\ot_R g'\  | \ g,g'\in {\mathcal
  G}^1\ \}$,
\[
\vartheta: B \ot_R B \to B \ot_R B,\qquad  g\ot_R g'\mapsto g g'\ot_R g' .
\]
(Note that $R$-module structures in the domain and codomain are different,
cf. Definition \ref{def:x_R_Hopf}.)
It obviously has an inverse $\vartheta^{-1}(g\ot_R g')=gg^{\prime -1}\ot_R
g'$.

In this final section we apply the theory developed in the earlier
sections to the groupoid bialgebroid ${\mathcal B}$ and its stable anti
Yetter-Drinfel'd modules. {In this way, we obtain expressions for
  Hochschild and cyclic homologies of a groupoid with finitely many
  objects. Describing then any groupoid as a direct limit of groupoids with
  finitely many objects, we extend the computation of cyclic homology to
  arbitrary groupoids. Similar arguments don't seem to apply in case of
  Hochschild homology.} 

\subsection{Anti Yetter-Drinfel'd modules for groupoids}\label{sec:gr_SAYD}
The subject of this section is a complete characterization of (stable) anti
Yetter-Drinfel'd modules of a groupoid bialgebroid. As a first step, we study
comodules.

\begin{proposition}\label{prop:comod_decomp}
Let ${\mathcal G}$ be a small groupoid with finite set of objects. Let
$\mathcal B$ be the groupoid $R$-bialgebroid associated to ${\mathcal G}$.
Then any right ${\mathcal B}$-comodule $M$ has a direct sum decomposition
$M\cong \oplus_{g\in {\mathcal G}^1} M_g$ (as an $R$-module) such that the
$R$-action $\cdot$ and the $\mathcal B$-coaction $\varrho$ satisfy the
conditions
\begin{itemize}
\item[{(i)}] $m\cdot x = \delta_{x,s(g)} m\quad $
and
\item[{(ii)}] $\varrho(m) = m\ot_R g$,
\end{itemize}
for $x\in {\mathcal G}^0$, $g\in {\mathcal G}^1$ and $m\in M_g$.
Conversely, on an $R$-module $M\cong \oplus_{g\in {\mathcal G}^1} M_g$,
which is subject to condition (i), there is a unique right $\mathcal
B$-coaction satisfying (ii).
\end{proposition}

\begin{proof}
Recall that by definition $R$ acts on $B$ by right multiplication.
For any $g\in {\mathcal G}^1$ there is an $R$-module map $\chi_g:B\to R$,
$h\mapsto \delta_{g,h} s(g)$. Introduce the map
\[
\pi_g:=(M\ot_R \chi_g)\circ \varrho: M\to M.
\]
We claim that $M$ is isomorphic to a direct sum of the $R$-modules $M_g:=
\mathrm{Im} (\pi_g)$. Since $B$ is a free ${\mathbb K}$-module, there exist
(non-unique) elements $\{\ m_g\ |\ g\in {\mathcal G}^1\ \}$ in $M$, in terms of
which $\varrho(m)=\sum_{g\in {\mathcal G}^1} m_g \ot_R g$, hence $\pi_g(m)= m_g
\cdot s(g)$.
By construction, for a given element $m\in M$ there are only finitely many
  elements $g\in {\mathcal G}^1$ such that $m_g\neq 0$, hence $\pi_g(m)\neq
  0$.
By coassociativity of $\varrho$,
\[
\sum_{g,h\in {\mathcal G}^1} (m_h)_g \ot_R g\ot_R h =
\sum_{g  \in {\mathcal G}^1}  m_g \ot_R g\ot_R g.
\]
Hence applying $M\ot_R \chi_{g'}\ot_R \chi_{h'}$, for some $g',h'\in {\mathcal
  G}^1$, we conclude that $(m_{h'})_{g'}\cdot s(g')s(h') = \delta_{g',h'}
  m_{g'}\cdot s(g')$, i.e. $\pi_{g'} \circ \pi_{h'} = \delta_{g',h'}\pi_{g'}$.
By counitality of $\varrho$, for all $m\in M$,
\[
\sum_{g  \in {\mathcal G}^1} \pi_g(m) =\sum_{g  \in {\mathcal G}^1} m_g \cdot
s(g) =m.
\]
This proves the direct sum decomposition of $M$ as an $R$-module. Condition
(i) follows by the computation, for $m\in M$, $g\in {\mathcal G}^1$ and $x\in
{\mathcal G}^0$,
\[
\pi_g(m)\cdot x = m_g \cdot s(g) x = \delta_{s(g),x} m_g \cdot s(g) =
\delta_{s(g),x} \pi_g(m).
\]
Moreover, for $m\in M$,
\[
\varrho(m) =\sum_{g  \in {\mathcal G}^1} m_g \ot_R g =
\sum_{g  \in {\mathcal G}^1} m_g \ot_R g s(g)=
\sum_{g  \in {\mathcal G}^1} m_g\cdot s(g) \ot_R g =
\sum_{g  \in {\mathcal G}^1} \pi_g(m) \ot_R g.
\]
Hence, by orthogonality of the projections $\pi_g$,
\[
\varrho(\pi_h(m))=\sum_{g  \in {\mathcal G}^1} \pi_g(\pi_h(m)) \ot_R g=
\pi_h(m) \ot_R h,
\]
which proves condition (ii).
Conversely, assume that for an $R$-module $M\cong \oplus_{g\in {\mathcal G}^1}
M_g$ condition (i) holds. Put $\varrho_g:M_g \to M_g \ot_R B$, $m\mapsto
m\ot_R g$. One can check that it makes $M_g$ to a right ${\mathcal
  B}$-comodule. By universality of a direct sum, this defines a unique right
${\mathcal B}$-coaction $\varrho$ on $M$, such that condition (ii) holds.
\end{proof}

Note that property (ii) in Proposition \ref{prop:comod_decomp} characterizes
uniquely the elements $m$ of a component $M_g$. Indeed, if for $m\in M$,
$\varrho(m)=m\ot_R g$, then by counitality of the coaction $\varrho$ we obtain
$\pi_g(m)=m\cdot s(g)=m$.

We are ready to characterize (stable)  anti Yetter-Drinfel'd modules of
groupoids.

\begin{theorem}\label{thm:gr_SAYD}
Let ${\mathcal G}$ be a small groupoid with finite set of objects. Let
$\mathcal B$ be the groupoid bialgebroid associated to ${\mathcal G}$.
A left $\mathcal B$-module $M\cong \oplus_{g\in {\mathcal G}^1} M_g$, with
action $\triangleright$, is an anti Yetter-Drinfel'd $\mathcal B$-module if and
only if the following conditions hold.
\begin{itemize}
\item[{(i)}] For $x\in {\mathcal G}^0$, $g\in {\mathcal G}^1$ and $m\in M_g$,
$\delta_{x,s(g)} m =x\triangleright m = \delta_{x,t(g)} m$.
\item[{(ii)}] For $g,h\in {\mathcal G}^1$ and $m\in M_g$, the element
  $h\triangleright m$ is zero if $hgh^{-1}=0$ in $B$ and it
belongs to $M_{hgh^{-1}}$ if $hgh^{-1}\neq 0$ in $B$.
\end{itemize}
The  anti Yetter-Drinfel'd $\mathcal B$-module $M$ is stable if and
only if in addition $g\triangleright m =m$, for all $g\in {\mathcal G}^1$ and
$m\in M_g$.
\end{theorem}

\begin{proof}
By Proposition \ref{prop:comod_decomp} there is a unique ${\mathcal
  B}$-coaction on $M$ corresponding to the given direct sum decomposition. For
  $g\in {\mathcal G}^1$, it takes $m\in M_g$ to $m\ot_R g$. The
$R$-bimodule structure corresponding to this coaction comes out, for $x,y\in
{\mathcal G}^0$, $g\in {\mathcal G}^1$ and $m\in M_g$, as
\[
x\cdot m \cdot y = m \cdot\epsilon(xgy) =  \delta_{s(g),y}\, \delta_{t(g),x} \
m \cdot s(g)=\delta_{s(g),y} \, \delta_{t(g),x} \, m.
\]
Hence axiom \eqref{eq:Y-D_R_mod} of an anti Yetter-Drinfel'd module translates
to condition (i) in the theorem. A straightforward computation shows that for
a groupoid bialgebroid, axiom \eqref{eq:Y-D_comp} of an anti Yetter-Drinfel'd
module takes the form
\begin{equation}\label{eq:gr_YD}
\varrho(h\triangleright m) = h\triangleright m\ot_R hgh^{-1},
\end{equation}
for $h,g\in {\mathcal G}^1$ and $m\in M_g$. If $hgh^{-1}=0$ in $B$, then the
right hand side of \eqref{eq:gr_YD} vanishes. Since $\varrho$ is a
monomorphism of $R$-modules (split by $M\ot_R \epsilon$), this is equivalent to
$h\triangleright m=0$. If $hgh^{-1}\neq 0$ in $B$, then
\eqref{eq:gr_YD} is equivalent to $h\triangleright m\in M_{hgh^{-1}}$
(cf. Proposition \ref{prop:comod_decomp}, and discussions following it). This
completes the proof.
\end{proof}

As a consequence of Theorem \ref{thm:gr_SAYD} (condition (i)),
in the direct sum decomposition of an anti Yetter-Drinfel'd
module $M$ of a groupoid bialgebroid,
$M_g$ is non-zero for only those elements
$g\in {\mathcal G}^1$ for which $s(g)=t(g)$, i.e. which are {\em
  loops}. That is, introducing the notation ${\mathcal L}({\mathcal G}):=\{\
g\in {\mathcal G}^1 \ |\ s(g)=t(g)\ \}$, one can write
$M\cong \oplus_{l\in {\mathcal L}({\mathcal G})} M_l$.

Our next aim is to decompose an anti Yetter-Drinfel'd module $M$ of  a
groupoid bialgebroid ${\mathcal B}$ as a direct sum {\em of anti
    Yetter-Drinfel'd modules}.
For a loop $l\in {\mathcal L}({\mathcal G})$, denote by $[l]$ the orbit of $l$
in ${\mathcal L}({\mathcal G})$ for the adjoint action, that is, the set of
different non-zero elements of the form $glg^{-1}$, as $g$ runs through
${\mathcal G}^1$. This gives a (${\mathcal G}$-invariant) partition ${\mathcal
T}({\mathcal G})$ of ${\mathcal L}({\mathcal G})$.
Using Theorem \ref{thm:gr_SAYD} one concludes that $M_{[l]}:=\oplus_{l'\in
  [l]} M_{l'}$ is an anti Yetter Drinfel'd ${\mathcal B}$-module, and
\begin{equation}\label{eq:gr_SAYD_decomp}
M \cong \bigoplus_{[l]\in
  {\mathcal T}({\mathcal G})} M_{[l]},
\end{equation}
as anti Yetter-Drinfel'd modules.
Let us give an alternative description of the anti Yetter-Drinfel'd module
$M_{[l]}$. Introduce the following subalgebras of the groupoid algebra $B$.
For $l\in {\mathcal L}({\mathcal G})$,
let $B_l$ be the group algebra of the centralizer ${\mathcal G}^1_l$ of $l$ in
the group
$\{\ l'\in {\mathcal L}({\mathcal G})\ |\ s(l')=s(l)\ \}$. That is,
\begin{equation}\label{eq:B_l}
B_l:= {\mathbb K}{\mathcal G}^1_l\equiv {\mathbb K} \{\ l'\in {\mathcal
  L}({\mathcal G})\ |\ l' l l^{\prime -1} =l\ \}.
\end{equation}
The algebra $B_l$ possesses a unit (the element $s(l)=t(l)$), different from
the unit element of $B$.
For $x\in {\mathcal G}^0$, let $B(x)$ denote the vector space spanned by the
elements $g\in {\mathcal G}^1$, such that $s(g)=x$.
Consider the group bialgebra structure of $B_l$.
Clearly, for any $l\in {\mathcal L}({\mathcal G})$, the component $M_l$ in the
direct sum decomposition of an anti Yetter Drinfel'd ${\mathcal B}$-module $M$
is
an anti Yetter-Drinfel'd module of the group bialgebra $B_l$ (in the sense of
\cite{JaSt:CycHom}, cf. Theorem \ref{thm:gr_SAYD}). Moreover, $B(s(l))$ is a
$B$-$B_l$ bimodule and a right $B_l$-module coalgebra (with coproduct induced
by the map $g\mapsto g\ot_{\mathbb K} g$, for $g\in {\mathcal G}^1$ such that
$s(g)=s(l)$.)

\begin{lemma} \label{lem:B(x)_free}
Let ${\mathcal G}$ be a small groupoid with a finite set of objects. Let
$\mathcal B$ be the groupoid bialgebroid over a field ${\mathbb K}$,
associated to ${\mathcal G}$.
Keeping the notation introduced above, $B(s(l))$ is free as a right
$B_l$-module, for any $l\in {\mathcal L}({\mathcal G})$.
\end{lemma}

\begin{proof}
Since any morphism in a groupoid is invertible, the right action
of the group ${\mathcal G}^1_l$ on the set $\{\ g\in {\mathcal
G}^1\ | \ s(g)=s(l)\ \}$ is faithful in the sense that $gl'=g$
implies $l'=s(l)$. Hence the claim follows by the fact that {
${\mathbb K}X$} is a free module for a group algebra ${\mathbb K}
G$ whenever $G$ acts faithfully on the set $X$. Indeed, fix a
${\mathbb K}$-basis $\{\ e_x\ |\ x\in X\ \}$ of {${\mathbb
K}X$}, such that $e_x \cdot g = e_{x\cdot g}$. Fix a section $f$
of the canonical epimorphism from $X$ to the set of $G$-orbits
$X/G$. By construction, $\{\ e_{f({\mathcal O})}\ |\ {\mathcal O}
\in X/G\ \}$ is a generating set of the ${\mathbb K} G$-module
{${\mathbb K}X$}. It is also linearly independent over
${\mathbb K} G$, by the following reasoning. Assume that, for some
coefficients $a_{\mathcal O}=\sum_{g\in G} \alpha_{{\mathcal
O},g}\,  g\in {\mathbb K} G$,
\[
0=\sum_{{\mathcal O} \in X/G} e_{f({\mathcal O})} \cdot a_{\mathcal O}=
\sum_{{\mathcal O} \in X/G}\sum_{g\in G}\alpha_{{\mathcal O},g} \,
e_{f({\mathcal O}) \cdot g}.
\]
Since $G$ acts on $X$ faithfully, in the above sum each element
$e_{f({\mathcal O})}$ appears exactly once. Hence $\alpha_{{\mathcal O},g}=0$,
for all ${\mathcal O} \in X/G$ and $g\in G$. Thus we have the claim proven.
\end{proof}

\begin{proposition}\label{prop:M_prod}
Let ${\mathcal G}$ be a small groupoid with finite set of objects. Let
$\mathcal B$ be the groupoid bialgebroid associated to ${\mathcal G}$.
Let $M$ be an anti Yetter-Drinfel'd ${\mathcal B}$-module. Keeping the
notation introduced above, there is an isomorphism of anti Yetter-Drinfel'd
${\mathcal B}$-modules
\[
M_{[l]} \cong B(s(l)) \ot_{B_l} M_l,\qquad
\textrm{for all }l\in {\mathcal L}({\mathcal G}).
\]
\end{proposition}
\begin{proof}
First we construct a left $B$-module isomorphism
\begin{equation}\label{eq:phi}
\varphi_l:B(s(l)) \ot_{B_l} M_l \to M_{[l]}, \qquad g\ot_{B_l} m \mapsto
g\triangleright m.
\end{equation}
For any $g\in {\mathcal G}^1$ such that $s(g)=s(l)$, consider the map
\[
\psi_g:M_{glg^{-1}}\to B(s(l)) \ot_{B_l} M_l,\qquad
g\mapsto g \ot_{B_l} g^{-1}\triangleright m.
\]
In order to see that the map $\psi_g$ does not depend on $g$, only on $g l
g^{-1}$,
choose another element $h\in {\mathcal G}^1$, such that $g l g^{-1}= h l
h^{-1}$. Note that this implies in particular $t(g)=t(h)$. Then, for $m\in
M_{glg^{-1}}$,
\begin{eqnarray*}
\psi_h(m)&\!\!\!\!=&\!\!\!\!
h \ot_{B_l} h^{-1}\triangleright m =
h \ot_{B_l} h^{-1} t(h)\triangleright m =
h \ot_{B_l} h^{-1} t(g)\triangleright m =
h \ot_{B_l} h^{-1}gg^{-1}\triangleright m \\ &\!\!\!\!=&\!\!\!\!
h h^{-1}g \ot_{B_l} g^{-1}\triangleright m =
t(h) g \ot_{B_l} g^{-1}\triangleright m =
t(g) g \ot_{B_l} g^{-1}\triangleright m =
g \ot_{B_l} g^{-1}\triangleright m =
\psi_g(m),
\end{eqnarray*}
where in the fifth equality we used that $h^{-1}g$ is an element of $B_l$.
Thus we  conclude by universality of a direct sum on the existence of a
map $\psi:M_{[l]} \cong \oplus_{l'\in [l]} M_{l'} \to B(s(l)) \ot_{B_l} M_l$,
mapping $m\in M_{glg^{-1}}$ to $\psi_g(m)=g \ot_{B_l} g^{-1}\triangleright
m$. A straightforward computation shows that $\psi$ is the inverse of
$\varphi_l$ in \eqref{eq:phi}.

Next we show that $B(s(l)) \ot_{B_l} M_l$ is an anti Yetter-Drinfel'd
module with respect to the direct sum decomposition $B(s(l)) \ot_{B_l} M_l
\cong \oplus_{l'\in [l]} \psi(M_{l'})$, hence
\eqref{eq:phi} is an isomorphism of anti Yetter-Drinfel'd modules, as stated.
For $m\in M_l$, $l'\in {\mathcal L}({\mathcal G})$, $h,h'\in {\mathcal G}^1$
and $y\in {\mathcal G}^0$, such that $hlh^{-1}=l'$ and $s(h')=t(h)$,
\begin{eqnarray*}
y\triangleright (h\ot_{B_l} m)&=&yh \ot_{B_l} m =\delta_{t(h),y}\, h \ot_{B_l}
m = \delta_{s(l'),y}\, h \ot_{B_l} m, \\
h'\triangleright (h\ot_{B_l} m)&=& h' h \ot_{B_l}  m \in \psi(M_{h'
  l' h^{\prime -1}}).
\end{eqnarray*}
In view of Theorem \ref{thm:gr_SAYD} this implies that $B(s(l)) \ot_{B_l} M_l$
is an anti Yetter-Drinfel'd module with respect to the given
decomposition, hence it completes the proof.
\end{proof}

\subsection{Hochschild and cyclic homology with coefficients}
\label{sec:gr_hom_coef}
In view of Proposition \ref{prop:dual_YD_proj}, there is a cyclic simplex
associated to a groupoid bialgebroid ${\mathcal B}$, the right regular
${\mathcal B}$-module coring and any stable anti Yetter-Drinfel'd ${\mathcal
B}$-module $M$. At degree $n$, it is given by $Z_n(B,M)=B^{\ot_R\, n+1}\ot_B
M$ (where $B$ acts on $B^{\ot_R\, n+1}$ via the diagonal right action).
In this section we compute its Hochschild and cyclic homologies.

With an eye on the decomposition of $M$ in Section \ref{sec:gr_SAYD},
computations start with following

\begin{lemma}\label{lem:B_prod}
Let ${\mathcal G}$ be a small groupoid with finite set of objects. Let
$\mathcal B$ be the groupoid bialgebroid over a field $\mathbb K$, associated
to ${\mathcal G}$.
Let $M$ be an anti Yetter-Drinfel'd ${\mathcal B}$-module. Using notations
introduced in Section \ref{sec:gr_SAYD}, there is an isomorphism of right
$B_l$-modules
\[
B^{\ot_R\, n+1}\ot_B B(x) \cong B(x)^{\ot_{\mathbb K}\, n+1},
\]
for all $l\in {\mathcal L}({\mathcal G})$ and $x:=s(l)$. Here $B^{\ot_R\,
  n+1}\ot_B B(x)$ is understood to be a right $B_l$-module via the last factor
  and the group algebra $B_l$ acts on $B(x)^{\ot_{\mathbb K}\, n+1}$ via the
  diagonal action.
\end{lemma}

\begin{proof}
Since $R={\mathbb K}{\mathcal G}^0$ is a separable
${\mathbb K}$-algebra, $B^{\ot_R\, n+1}$ is
isomorphic to the subspace $B^{\times\, n+1}$ of $B^{\ot_{\mathbb K}\, n+1}$,
spanned by those elements $g_0\ot_{\mathbb K}\dots \ot_{\mathbb K} g_n$ for
which $s(g_i)=s(g_{i+1})$, for all $i=0\dots n$. Thus it suffices to
prove $B^{\times\, n+1}\ot_B B(x) \cong B(x)^{\ot_{\mathbb K}\, n+1}$.
We claim that the right $B_l$-module map
\begin{equation}\label{eq:iso}
B(x)^{\ot_{\mathbb K}\, n+1} \to B^{\times\, n+1}\ot_B B(x), \qquad
g_0\ot_{\mathbb K} \dots \ot_{\mathbb K} g_n\mapsto
(g_0\ot_{\mathbb K} \dots \ot_{\mathbb K} g_n)\ot_B x
\end{equation}
is an isomorphism. Since the map
\[
B^{\times\, n+1}\ot_{\mathbb K} B(x) \to B(x)^{\ot_{\mathbb K}\, n+1},\qquad
(g_0\ot_{\mathbb K}\dots \ot_{\mathbb K} g_n) \ot_{\mathbb K} h \mapsto
g_0 h \ot_{\mathbb K}\dots \ot_{\mathbb K} g_n h\qquad
\]
factorizes through $B^{\times\, n+1}\ot_B B(x)$, it defines a unique map
\begin{equation}\label{eq:inv}
B^{\times\, n+1}\ot_{B} B(x) \to B(x)^{\ot_{\mathbb K}\, n+1},\qquad
(g_0\ot_{\mathbb K}\dots \ot_{\mathbb K} g_n) \ot_B h \mapsto
g_0 h \ot_{\mathbb K}\dots \ot_{\mathbb K} g_n h.
\end{equation}
Obviously, \eqref{eq:iso} and  \eqref{eq:inv} are mutual inverses.
\end{proof}

\begin{proposition}\label{prop:red_to_grp}
Let ${\mathcal G}$ be a small groupoid with finite set of objects. Let
$\mathcal B$ be the groupoid bialgebroid over a field $\mathbb K$, associated
to ${\mathcal G}$.
Let $M$ be a stable anti Yetter-Drinfel'd ${\mathcal B}$-module, with
decomposition
\eqref{eq:gr_SAYD_decomp}.
For all $l\in {\mathcal L}({\mathcal G})$, the cyclic simplex
$Z_\ast(B,M_{[l]})=B^{\ot_R\, \ast+1}\ot_B M_{[l]}$ is isomorphic to
$Z_\ast(B(s(l)),M_l)=B(s(l))^{\ot_{\mathbb K}\, \ast+1}\ot_{B_l} M_l $,
corresponding (as
in \cite[Remark 4.16]{JaSt:CycHom}) to the module coalgebra $B(s(l))$, and
stable anti Yetter-Drinfel'd module $M_l$, of the group bialgebra $B_l$.
\end{proposition}

\begin{proof}
Combining the isomorphisms in Proposition \ref{prop:M_prod} and Lemma
\ref{lem:B_prod}, we obtain an isomorphism of vector spaces
\[
\begin{array}{rcl}
B(s(l))^{\ot_{\mathbb K}\, n+1}\ot_{B_l} M_l \cong
& B^{\ot_R\, n+1}\ot_B B(s(l)) \ot_{B_l} M_l \cong
& B^{\ot_R\, n+1}\ot_B M_{[l]} \\
(g_0\ot_{\mathbb K} \dots \ot_{\mathbb K} g_n)\ot_{B_l} m
& \longmapsto
& (g_0\ot_R \dots \ot_R g_n)\ot_{B} m,
\end{array}
\]
with inverse $B^{\ot_R\, n+1}\ot_B M_{hlh^{-1}} \ni (g_0\ot_R \dots \ot_R
  g_n)\ot_{B} m \mapsto (g_0 h\ot_{\mathbb K} \dots \ot_{\mathbb K} g_n
  h)\ot_{B_l} h^{-1}\triangleright m$.
It is left to the reader to check that it is an isomorphism of cyclic objects.
\end{proof}

Hochschild and cyclic homologies of a group, with coefficients in a stable anti
Yetter-Drinfel'd module, were computed in \cite[Corollary 5.13]{JaSt:CycHom}.
Hence in what follows we relate Hochschild and cyclic homologies of the cyclic
object $Z_\ast (B(s(l)),M_l)\cong Z_\ast(B,M_{[l]})$ in Proposition
\ref{prop:red_to_grp} to the known respective homology of $Z_\ast(B_l,M_l)$.

\begin{lemma}\label{lem:free_res}
Let $(C,\Delta,\epsilon)$ be a coalgebra over a field ${\mathbb K}$.
Consider the corresponding simplex $C^{\ot_{\mathbb
    K}\, \ast +1}$ (with face maps $\partial_i$ induced by $\epsilon$ and
degeneracy maps $\mu_i$ induced by $\Delta$). Then the associated complex
${\widetilde C}_\ast(C)$ is acyclic, i.e. $H_n({\widetilde C}_\ast(C))=
\delta_{n,0} {\mathbb K}$, for any non-negative integer $n$. Moreover,
if $(C,\Delta,\epsilon)$ is a right module coalgebra of a ${\mathbb K}$-Hopf
algebra $H$ and free as a right $H$-module then ${\widetilde C}_\ast (C)$ is a
free resolution of ${\mathbb K}$ in the category of right $H$-modules.
\end{lemma}

\begin{proof}
We need to show that the chain complex
\[
\xymatrix{
\dots \ar[r]^-{\delta_{n+1}}& C^{\ot_{\mathbb K}\, n +1}\ar[r]^{\delta_{n}}&
C^{\ot_{\mathbb K}\, n}\ar[r]&\dots \ar[r]^{\delta_{1}}& C
\ar[r]^{\delta_{0}=\epsilon}& {\mathbb K} \ar[r]&0}
\]
is exact, where $\delta_n =\sum_{i=0}^n (-1)^i \partial_i$. Indeed, the map
$C\ot_{\mathbb K} \epsilon$ is surjective, having a section $\Delta$. Since
$C$ is a faithfully flat module of the field ${\mathbb K}$, this implies
  surjectivity of $\epsilon$. Take an element $g\in C$ such that
  $\epsilon(g)=1_{\mathbb K}$. One easily checks
that, for $z\in C^{\ot_{\mathbb K}\, n +1}$ such that $\delta_n(z)=0$,
$\delta_{n+1}(g\ot_{\mathbb K} z)=z-(g\ot_{\mathbb K}
\delta_n(z))=z$.

Assume now that $(C,\Delta,\epsilon)$ is a right module
coalgebra of a ${\mathbb K}$-Hopf algebra $H$. Then ${\widetilde C}_\ast (C)$
is an acyclic complex in the category of right H-modules. It remains to show
that $C^{\otimes_{\mathbb K}\, n+1}$ is a free $H$-module whenever $C$ is
so. Indeed,
in this case $C^{\otimes_{\mathbb K}\, n+1}$ is free as a right module of
$H^{\otimes_{\mathbb K}\, n+1}$ via factorwise action.
By \cite[Lemma 2.10]{JaSt:CycHom} (cf. Lemma \ref{lem:beta}),
$H^{\otimes_{\mathbb K}\, n+1}$ is free as a right $H$-module via the
diagonal action.
Hence $C^{\otimes_{\mathbb K}\, n+1}\cong C^{\otimes_{\mathbb
    K}\, n+1} \ot_{H^{\otimes_{\mathbb K}\, n+1}}H^{\otimes_{\mathbb K}\,
  n+1}$ is a free right $H$-module.
\end{proof}

\begin{lemma} \label{lem:hom_iso}
Let $H$ be a Hopf algebra over a field ${\mathbb K}$ and let $i:C\to C'$ be a
morphism of right $H$-module coalgebras. Let $M$ be a stable anti
Yetter-Drinfel'd $H$-module. If both $C$ and $C'$
are free as right $H$-modules, then the induced morphism
$i_\ast:Z_\ast(C,M)\to Z_\ast (C',M)$ of cyclic objects gives rise to
isomorphisms both of Hochschild and cyclic homologies. That is,
\[
\mathrm{HH}_\ast(C,M)\cong \mathrm{HH}_\ast (C',M)\qquad
\textrm{and}\qquad \mathrm{HC}_\ast(C,M)\cong \mathrm{HC}_\ast
(C',M).
\]
\end{lemma}

\begin{proof}
By Lemma \ref{lem:free_res}, both ${\widetilde C}_\ast (C)$ and ${\widetilde
  C}_\ast (C')$ are free resolutions of ${\mathbb K}$ over $H$. Thus
  $i^{\otimes_{\mathbb K}\, \ast+1}:{\widetilde C}_\ast (C) \to {\widetilde
  C}_\ast (C')$ is a quasi-isomorphism. Then so is $i^{\otimes_{\mathbb K}\,
  \ast+1}\ot_H M:{\widetilde C}_\ast (C)\ot_H M\to {\widetilde C}_\ast
  (C')\ot_H M$, yielding an isomorphism
\[
\mathrm{HH}_n(C,M)\equiv \mathrm{H}_n({\widetilde C}_\ast (C)\ot_H
M)\cong \mathrm{H}_n({\widetilde
  C}_\ast (C')\ot_H M) \equiv \mathrm{HH}_n(C',M),
\]
for all non-negative integer $n$. The isomorphism of Hochschild homologies
implies the isomorphism of cyclic homologies, see e.g. (a dual form of)
\cite[2.2.3]{Loday}.
\end{proof}

\begin{theorem}\label{thm:gr_main}
Let ${\mathcal G}$ be a small groupoid with finite set of objects. Let
$\mathcal B$ be the groupoid bialgebroid over a field $\mathbb K$, associated
to ${\mathcal G}$.
Let $M$ be a stable anti Yetter-Drinfel'd ${\mathcal B}$-module, with
decomposition \eqref{eq:gr_SAYD_decomp}.
For all $l\in {\mathcal L}({\mathcal G})$,
\[
\mathrm{HH}_\ast(B,M_{[l]})\cong \mathrm{HH}_\ast (B_l,M_l)\qquad
\textrm{and}\qquad \mathrm{HC}_\ast(B,M_{[l]})\cong
\mathrm{HC}_\ast (B_l,M_l).
\]
\end{theorem}

\begin{proof}
By Proposition \ref{prop:red_to_grp}, the cyclic objects
$Z_\ast(B,M_{[l]})$ and $Z_\ast(B(s(l)),M_l)$ are isomorphic.
Hence $\mathrm{HH}_\ast(B,M_{[l]})\cong
\mathrm{HH}_\ast(B(s(l)),M_l)$. Furthermore, the right
$B_l$-module coalgebra $B(s(l))$ is a free right $B_l$-module by
Lemma \ref{lem:B(x)_free}. Hence we conclude by Lemma
\ref{lem:hom_iso} that the inclusion map $B_l\hookrightarrow
B(s(l))$ induces an isomorphism
$\mathrm{HH}_\ast(B(s(l)),M_l)\cong \mathrm{HH}_\ast(B_l,M_l)$.
Combination of these isomorphisms proves the theorem.
\end{proof}

{Let us assume now that $\mathbb{K}$ is a field of
characteristic zero.} In view of Theorem \ref{thm:gr_main}, one
can compute also
\[
\mathrm{HH}_\ast(B,M)\cong \bigoplus_{[l]\in {\mathcal
T}({\mathcal G})}\mathrm{HH}_\ast \left(B_l,M_l\right)\qquad
\textrm{and}\qquad \mathrm{HC}_\ast(B,M)\cong \bigoplus_{[l]\in
{\mathcal T}({\mathcal G})} \mathrm{HC}_\ast \left(B_l,M_l\right).
\]
Since $B_l$ is a group algebra of ${\mathcal G}^1_l$ (cf.
\eqref{eq:B_l}), $\mathrm{HH}_\ast(B_l,M_l)=H_\ast({\mathcal
G}^1_l,M_l)$ is the group homology of ${\mathcal G}^1_l$ with
coefficients in $M_l$. Applying \cite[Corollary
5.13]{JaSt:CycHom}, for the cyclic homology we get
\[
\mathrm{HC}_\ast (B_l,M_l) = \left\{
\begin{array}{ll}
\bigoplus_{i\geq 0} \mathrm{H}_{\ast -2i}\big({{\mathcal
G}^1_l}/{\langle l\rangle} ,M_l\big),
        \quad &\textrm{if } l \textrm{ has finite order}\\
\mathrm{H}_{\ast}\big({{\mathcal G}^1_l}/{\langle l\rangle}
,M_l\big), \quad
        &\textrm{if } l \textrm{ has infinite order}.
\end{array}
\right.
\]

\subsection{Cyclic homology of groupoids}
The results in Section \ref{sec:gr_hom_coef} can be specialized further to
stable anti Yetter-Drinfel'd modules provided by groupoid Galois extensions,
cf. Proposition \ref{prop:T_S}. This enables us, in particular, to compute
ordinary (i.e. non-relative) Hochschild and cyclic homologies of a groupoid.

Let ${\mathcal G}$ be a small groupoid with finite set of objects. Let
$\mathcal B$ be the groupoid bialgebroid over a field $\mathbb K$,
associated to ${\mathcal G}$.
Recall from \cite[Section 3]{CaDeGr:WHA_Gal} that a Galois extension
$S\subseteq T$ by $\mathcal B$ has an equivalent description as follows.
$T$ is a strongly ${\mathcal G}$-graded ${\mathbb K}$-algebra, that is,
$T\cong \oplus_{g\in {\mathcal G}^1} T_g$ (as a vector space),
\[
T_gT_{g'}=\left\{
\begin{array}{ll}
T_{g\circ g'}\quad &\textrm{if } s(g)=t(g')\\
0\quad &\textrm{if } s(g)\neq t(g'),
\end{array}
\right.
\]
$1_T=\sum_{x\in {\mathcal G}^0} 1_{T_x}$, and $S$ is equal to the
subalgebra $\oplus_{x\in {\mathcal G}^0} T_x$.
Note that each direct summand $T_g$ is an $R={\mathbb K}{\mathcal G}^0$-module
via $v\cdot x:= v 1_{T_x}=\delta_{s(g),x}\, v$, for $v\in T_g$ and $x\in
{\mathcal G}^0$.

By Proposition \ref{prop:T_S}, $T_S:=S \widehat{\ot}_S T$ is a stable anti
Yetter-Drinfel'd module. In the following lemma its structure is investigated.

\begin{lemma}\label{lem:gr_Gal_SAYD}
Let ${\mathcal G}$ be a small groupoid with finite set of objects. Let
$\mathcal B$ be the groupoid bialgebroid over a field $\mathbb K$, associated
to ${\mathcal G}$. Let $S\subseteq T$ be a Galois extension by $\mathcal
B$. Using the notations introduced in Section \ref{sec:gr_SAYD}, there is an
isomorphism of anti Yetter-Drinfel'd modules
\[
T_S \cong \bigoplus_{l\in \mathcal{L}({\mathcal G})} S \widehat{\ot}_S T_l.
\]
The (Miyashita-Ulbrich) action is given, for $g\in {\mathcal G}^1$,
$l\in\mathcal{L}({\mathcal G})$ and $v\in T_l$, by
\[
g\triangleright p_{T_l}(v)= \sum_{i=1}^n p_{T_l}(b_i v a_i),
\]
where $p_{T_l}: T_l \to S \widehat{\ot}_S T_l$ denotes the canonical
epimorphism and the elements $\{a_1,\dots ,a_n\} \subseteq T_{g^{-1}}$ and
$\{b_1,\dots , b_n\}\subseteq T_g$ satisfy $\sum_{i=1}^n a_i b_i
=1_{T_{s(g)}}$.
\end{lemma}
\begin{proof}
Since $T\cong\oplus_{g\in {\mathcal G}^1} T_g$ is ${\mathcal G}$-graded, each
direct summand $T_g$ is an $S$-bimodule. Hence
\[
S {\widehat{\otimes}}_S T \cong \bigoplus_{g\in {\mathcal G}^1} S
{\widehat{\otimes}}_S T_g.
\]
We claim that only those elements $g\in {\mathcal G}^1$ give non-zero
contribution to this direct sum, for which $s(g)= t(g)$.
Recall that $S \widehat{\ot}_S T_g$ is isomorphic to
the quotient of $T_g$ with respect to the commutator subspace $[S,T_g]=\{
qv-vq\ |\ q\in S,\ v\in T_g\ \}$.
Take an element $g\in {\mathcal G}^1$ such that $s(g)\neq t(g)$ and an element
$v\in T_g$. By strong grading of $T$,
there exist elements $\{a_1,\dots ,a_n\} \subseteq T_{g^{-1}}$ and
$\{b_1,\dots , b_n\}\subseteq T_g$ such that $\sum_{i=1}^n a_i b_i
=1_{T_{s(g)}}$. Then
\[
\sum_{i=1}^n [v a_i, b_i]=
\sum_{i=1}^n v a_i b_i - b_i v a_i =v 1_{T_{s(g)}} =v,
\]
where the penultimate equality follows by the fact that, for all values of
$i$, $ b_i v a_i \in T_g T_g T_{g^{-1}}$ is zero by the assumption that
$s(g)\neq t(g)$. Since for all values of $i$, $v a_i\in T_{g^{-1}} T_g =
T_{s(g)} \subseteq S$, we conclude that $T_g\subseteq [S,T_g]$. Since the
converse inclusion is obvious, we have $S \widehat{\ot}_S T_g\cong
T_g/[S,T_g]=0$ proven.

In order to write down the Miyashita-Ulbrich action, consider again $g\in
{\mathcal G}^1$ and elements
$\{a_1,\dots ,a_n\} \subseteq T_{g^{-1}}$ and
$\{b_1,\dots , b_n\}\subseteq T_g$ such that $\sum_{i=1}^n a_i b_i
=1_{T_{s(g)}}$. The canonical map $\mathrm{can}:T\ot_S T \to T\ot_R B$
satisfies
\[
\mathrm{can}(\sum_{i=1}^n a_i \ot_S b_i)=\sum_{i=1}^n a_i b_i \ot_R g =
1_{T_{s(g)}}\ot_R g =1_T\cdot s(g) \ot_R g= 1_T \ot_R g.
\]
Hence the stated form of the Miyashita-Ulbrich action follows by
\eqref{eq:T_S_mod}.
\end{proof}

\begin{corollary}\label{cor:rel_gal}
Let ${\mathcal G}$ be a small groupoid with finite set of objects. Let
$\mathcal B$ be the groupoid bialgebroid over a field $\mathbb K$, associated
to ${\mathcal G}$. Let $S\subseteq T$ be a Galois extension by $\mathcal
B$. Using the notations introduced in Section \ref{sec:gr_SAYD}, the
$S$-relative Hochschild and cyclic homologies of $T$ are given, respectively,
by
\[
\mathrm{HH}_\ast(T/S)= \bigoplus_{[l]\in {\mathcal T}({\mathcal
G})} \mathrm{HH}_\ast (B_l, S \widehat{\ot}_S T_l) \qquad
\textrm{and}\qquad \mathrm{HC}_\ast(T/S)= \bigoplus_{[l]\in
{\mathcal T}({\mathcal G})} \mathrm{HC}_\ast (B_l, S
\widehat{\ot}_S T_l).
\]
\end{corollary}

\begin{proof}
By Theorem \ref{thm:main_iso}, $\mathrm{HH}_\ast(T/S)\cong
\mathrm{HH}_\ast(B,T_S)$ and $\mathrm{HC}_\ast(T/S)\cong
\mathrm{HC}_\ast(B,T_S)$. Hence the claim follows by Theorem
\ref{thm:gr_main} and considerations following it, together with
Lemma \ref{lem:gr_Gal_SAYD}.
\end{proof}

A particular example $R\subseteq B$ of a Galois extension by a groupoid
${\mathcal G}$ is provided by the inclusion of the base algebra $R={\mathbb
K}{\mathcal G}^0$ in the groupoid algebra $B={\mathbb K}{\mathcal
G}^1$. Applying Corollary \ref{cor:rel_gal} to it, we obtain formulae for the
$R$-relative Hochschild and cyclic homologies of the groupoid ${\mathcal
G}$. Since the base algebra $R={\mathbb K}{\mathcal G}^0$ is separable over
${\mathbb K}$, relative homologies coincide with ordinary ones,
cf. \cite{Kad:CycHom}. Therefore we obtain the following corollary, extending
results in \cite{Burg} on cyclic homology of groups. Similar expressions
were derived also by Crainic in \cite{Cra} for {\em \'etale} groupoids.

\begin{corollary}\label{co:groupoid}
Let $B$ be the groupoid algebra of a small groupoid ${\mathcal G}$
with finite set of objects {over a field ${\mathbb K}$ of
characteristic zero}. Using the notations introduced in Section
\ref{sec:gr_SAYD} we have
\begin{align*}
\mathrm{HH}_{\ast }(B)& =\bigoplus\limits_{[l]\in \mathcal{T}(\mathcal{G})}\mathrm{H}
_{\ast }(\mathcal{G}_{l}^{1},\mathbb{K)}, \\
\mathrm{HC}_{\ast }(B)& =\left( \bigoplus\limits_{\substack{ \lbrack l]\in
  \mathcal{T}
(\mathcal{G}) \\ \mathrm{ord}(l)<\infty }}\bigoplus\limits_{i\geq 0}\mathrm{H
}_{\ast -2i}(\mathcal{G}_{l}^{1}/\left\langle l\right\rangle ,\mathbb{K)}
\right) \bigoplus \left( \bigoplus\limits_{\substack{ \lbrack l]\in \mathcal{
T}(\mathcal{G}) \\ \mathrm{ord}(l)=\infty }}\mathrm{H}_{\ast }(\mathcal{G}
_{l}^{1}/\left\langle l\right\rangle ,\mathbb{K})\right) .
\end{align*}
\end{corollary}

Recall that a groupoid $\mathcal{G}$ is connected if, for any
$x,y\in
\mathcal{G}^{0},$ there exists at least one $g\in \mathcal{G}^{1}$ such that
$s(g)=x$ and $t(g)=y.$ For a connected groupoid $\mathcal{G}$ we fix an object
$x\in \mathcal{G}^{0}$ and we denote by $G$ the group of loops $l\in
\mathcal{L}(\mathcal{G})$ such that $s(l)=t(l)=x.$ (Since ${\mathcal G}$ is
connected by assumption, different choices of $x$ lead to isomorphic subgroups
$G$ of ${\mathcal G}$.) Let $\mathcal{T}(G)$ denote the set of conjugacy
classes in $G$ and let $\{g_{\sigma }\mid $ $\sigma \in \mathcal{T}(G)\}$ be a
transversal of $\mathcal{T}(G).$

\begin{lemma} \label{le:T(G)}
Let $\mathcal{G}$ be a connected groupoid. Keeping the above
notation, there is an one-to-one correspondence between $\mathcal{T}(
\mathcal{G})\ $and $\mathcal{T}(G).$
\end{lemma}

\begin{proof}
We have to show that
{any orbit in ${\mathcal L}({\mathcal G})$ for the adjoint action contains
  precisely one element of the set $\{g_\sigma\mid $ $\sigma \in
  \mathcal{T}(G)\}$. That is,} any loop $l$ in $\mathcal{G}$ is equivalent
to a certain $g_{\sigma }$ and that two elements $g_{\sigma }$
and $g_{\tau }$ are equivalent if, and only if $\sigma =\tau $.
First, let us take $l\in \mathcal{L}(\mathcal{G}).$ Since $\mathcal{G}$ is
connected, there is a morphism $g$ such that $s(g)=x$ and $t(g)=s(l).$ Let
$l^{\prime }:=g^{-1}\circ l\circ g.$ By construction, $l^{\prime
}\in G,\ $so there are $h\in G$ and $\sigma \in \mathcal{T}(G)$
such that $h\circ l^{\prime }\circ h^{-1}=g_{\sigma }.$ Since
$l=(g\circ h^{-1})\circ g_{\sigma }\circ (g\circ h^{-1})^{-1}$ it
follows that $l$ and $g_{\sigma }$ are conjugated in
$\mathcal{L}(\mathcal{G}).$

Let us take $\sigma $ and $\tau $ in $\mathcal{T}(G)$ and assume that $
g_{\sigma }$ and $g_{\tau }$ define the same element in $\mathcal{T}(
\mathcal{G}).$ Then there is $g\in \mathcal{G}^{1}$ such that
$g_{\sigma }=g\circ g_{\tau }\circ g^{-1}.$ Since both loops
$g_{\sigma }$ and $g_{\tau }\ $have the same source, we deduce
that the source and the target of $g$
must be $x.$ Hence $g\in G$ and the conjugacy classes of $g_{\sigma }$ and $
g_{\tau }$ in $\mathcal{T}(G)$ are equal. In conclusion, $\sigma
=\tau .$
\end{proof}

\begin{remark}
{Lemma \ref{le:T(G)} tells us, in particular, that $\{g_{\sigma }\mid \sigma
\in \mathcal{T(}G)\}$ is also a transversal of $\mathcal{T}(\mathcal{G}).$
Explicitly, a given element $g_\sigma$ represents the following orbit
$\widetilde{\sigma } \in {\mathcal T}({\mathcal G})$.
Recall that we defined the group $G$ in terms of a fixed object $x$ in
$\mathcal{G}^{0}.$ For every $y\in \mathcal{G}^{0}$ we pick up a fixed
$g_{y}\in \mathrm{Hom}_{\mathcal{G}}(x,y)$. It defines a group isomorphism $G
\to \mathrm{Hom}_{\mathcal{G}}(y,y)$, $h \mapsto g_{y}\circ h\circ
g_{y}^{-1}$. It maps $\sigma \in
\mathcal{T}(G)$ to the conjugacy class $\sigma_{y}=\{g_{y}\circ h\circ
g_{y}^{-1}\mid h\in \sigma \}$ in the group
$\mathrm{Hom}_{\mathcal{G}}(y,y)$. The orbit of $g_\sigma$ for the adjoint
${\mathcal G}$-action is $\widetilde{\sigma }=\bigcup_{y\in
  \mathcal{G}^{0}}\sigma _{y}$. } 
\end{remark}

\begin{corollary}
\label{te:fc groupoid}Let $\mathcal{G}$ be a connected groupoid
with finitely many objects. Let $M$ be a stable anti Yetter-Drinfel'd module
over $B$,
where $B$ is the groupoid algebra of $\mathcal{G}$ over a field of
characteristic zero.

(i) $M_{G}:=\bigoplus_{g\in G}M_{g}$ is a stable anti Yetter-Drinfel'd module
over $\mathbb{K}G$ and the inclusions $\mathbb{K}G\subseteq B$ and $
M_{G}\subseteq M$ induce isomorphisms
\begin{equation*}
\mathrm{HH}_{\ast }(B,M)\cong \mathrm{HH}_{\ast
}(\mathbb{K}G,M_{G})\qquad
\text{and\qquad }\mathrm{HC}_{\ast }(B,M)\cong \mathrm{HC}_{\ast }(\mathbb{K}
G,M_{G}).
\end{equation*}

(ii) The inclusion $\mathbb{K}G\subseteq B$ induces isomorphisms
\begin{equation*}
\mathrm{HH}_{\ast }(B)\cong \mathrm{HH}_{\ast }(\mathbb{K}G)\qquad \text{
and\qquad }\mathrm{HC}_{\ast }(B)\cong \mathrm{HC}_{\ast
}(\mathbb{K}G).
\end{equation*}
\end{corollary}

\begin{proof}
(i) Obviously $M_{G}$ is a stable anti Yetter-Drinfel'd module over
  $\mathbb{K}G.$  Consider the following diagram.
\begin{equation*}
\begin{xy} \xymatrix{ Z_{\ast }(\mathbb{K}G,M_{G}) \ar[d] \ar[r]
&\bigoplus_{\sigma \in \mathcal{T} (G)}Z_{\ast
}(\mathbb{K}G,M_{\sigma })\ar[r]\ar[d]&\bigoplus_{\sigma \in
\mathcal{T}(G)}Z_{\ast }(\mathbb{K}B_{g_{\sigma }},M_{g_{\sigma
}})\ar[d]\\ Z_{\ast }(B,M) \ar[r] &\bigoplus_{{\sigma} \in
\mathcal{T}({G})}Z_{\ast }(B,M_{\widetilde{\sigma}
})\ar[r]&\bigoplus_{{\sigma}\in{T}({G})}Z_{\ast
}(B(x),M_{g_{\sigma }}) }\end{xy}
\end{equation*}
The leftmost morphism in the bottom row comes from the decomposition (\ref
{eq:gr_SAYD_decomp}), while the other one is the direct sum of the
arrows that were constructed in Proposition \ref{prop:red_to_grp}.
Note that both maps are induced by appropriate inclusions and they
are isomorphisms. The morphisms in the top row have the same
properties, as any group can be regarded as a groupoid with one
object. By definition, the vertical arrows are the canonical
morphisms induced by inclusions, so they make the squares
commutative. In view of the proof of Theorem \ref{thm:gr_main},
the rightmost vertical arrow gives isomorphisms both of
Hochschild and cyclic homologies. Then also the leftmost vertical
morphism does so.

(ii) {The subalgebra $R:={\mathbb K} X$ of $B$ is a separable ${\mathbb
    K}$-algebra. Hence $\mathrm{HH}_{\ast}(B)\cong \mathrm{HH}_{\ast }(B,R).$}
By Lemma \ref{lem:gr_Gal_SAYD} we know that $M:=R\widehat{\otimes }_{R}B$
is a stable anti Yetter-Drinfel'd module over $B$. Since {for any $g\in
  {\mathcal G}$,} either $[R,\mathbb{K}g]=\mathbb{K}g$ or $[R,
\mathbb{K}g]=0,$ depending on the fact that $g$ is a loop or not, we get $
M:=\bigoplus_{g\in \mathcal{L}(\mathcal{G)}}\mathbb{K}g.$ Hence $M_{G}:=
\mathbb{K}G.$ Obviously the action of $B$ on $\mathbb{K}G$ induced
from the Ulbrich-Miyashita action is the adjoint action of
$\mathbb{K}G$ on itself. We conclude by applying the
first part of the corollary and the isomorphisms $\mathrm{HH}_{\ast }(
\mathbb{K}G)\cong \mathrm{HH}_{\ast }(\mathbb{K}G,(\mathbb{K}G)_{ad})$ and $
\mathrm{\mathrm{HC}}_{\ast }(\mathbb{K}G)\cong \mathrm{\mathrm{HC}}_{\ast }(
\mathbb{K}G,(\mathbb{K}G)_{ad}),$ cf. \cite{JaSt:CycHom}.
\end{proof}

Our final aim is to compute \textrm{HC}$_{\ast }(B),$ the ordinary
cyclic homology of the groupoid algebra $B$ of a groupoid
$\mathcal{G}$ that may have infinite number of objects. For
such a groupoid, its groupoid algebra $B$ is not unital anymore.
Nevertheless, to define cyclic homology of $B$ one can proceed as for
unital algebras, cf. \cite[Chapter 2, \S 2.1]{Loday}. The point is
that Connes' complex $C^{\lambda }(B)$ still exists, although it
is now associated to a precyclic object, that is to a
presimplicial structure endowed with cyclic operators. Here, by
presimplicial object we mean a sequence of objects together only
with face maps. The defining properties of face maps and cyclic
operators are the same as in the definition of cyclic objects,
neglecting of course the relations that involve the degeneracy
maps. {The key ingredient of our computation is a description of any
groupoid as a direct limit of groupoids with finitely many objects. Since
cyclic homology is defined as homology of Connes' complex and the
homology functor commutes with direct limits, we obtain cyclic homology of an
arbitrary groupoid as a direct limit. Note however that, for non-unital
algebras, Hochschild homology is constructed in a different way. For the
definition, see for example \cite[Chapter 1, \S 1.2]{Loday}. Therefore, we can
not apply the same arguments to compute Hochschild homology of an arbitrary
groupoid.} 

\begin{theorem} \label{thm:inf_conn}
Let $\mathcal{G}$ be a connected groupoid. If $B$ is the groupoid
algebra over a field $ \mathbb{K}$ of characteristic zero, then
\begin{equation*}
\mathrm{HC}_{\ast }(B)\cong \mathrm{\mathrm{HC}}_{\ast
}(\mathbb{K}G).
\end{equation*}
\end{theorem}

\begin{proof}
Let $x$ be a given object in $\mathcal{G}^{0}$. Let $G$ denote the
group of loops $l\in \mathcal{L}(\mathcal{G})$ such that $s(l)=x$.
We order the set
\begin{equation*}
\mathfrak{X}:=\{X\subset \mathcal{G}^{0}\mid x\in X\text{ and
}X\text{ is finite}\}
\end{equation*}
with respect to inclusion. Trivially $\mathfrak{X}$ is a direct system. For $
X\in \mathfrak{X}$ we define $\mathcal{G}_{X}$ to be the full
subgroupoid of $\mathcal{G}$ such that $\mathcal{G}_{X}^{0}=X$ and
we denote its groupoid
algebra by $B_{X}.$ Note that $B=\bigcup_{X\in \mathfrak{X}}B_{X},$ so $
C^{\lambda }_\ast(B)=\bigcup_{X\in \mathfrak{X}}C_{\ast }^{\lambda
}(B_{X}).$ As the direct limit in the category of vector spaces is
exact, it follows that the homology functor and direct limit
commute. Thus
\begin{equation*}
\mathrm{HC}_{\ast }(B)\cong \mathrm{H}_{\ast }(C^{\lambda }_{\ast }(B))\cong
\mathrm{
H}_{\ast }\big(\lim_{\overrightarrow{X\in \mathfrak{X}}}C_{\ast
}^{\lambda }(B_{X})\big)\cong \lim_{\overrightarrow{X\in
\mathfrak{X}}}\mathrm{H}_{\ast }\big(C_{\ast }^{\lambda
}(B_{X})\big)\cong \lim_{\overrightarrow{X\in
\mathfrak{X}}}\mathrm{HC}_{\ast }(B_{X}),
\end{equation*}
where the latter direct system is defined by the canonical maps $\mathrm{HC}
_{\ast }(B_{X})\rightarrow \mathrm{HC}_{\ast }(B_{Y}),$ with $X,$ $Y$ in $
\mathfrak{X}$ such that $X\subset Y.$ We claim that these maps are
isomorphisms. Indeed, let us consider the following commutative diagram
\begin{equation*}
\xymatrix{ & Z_\ast(\mathbb{K}G) \ar[dr] \ar[dl] \\ Z_\ast(B_X)
\ar[rr] & & Z_\ast (B_Y) }
\end{equation*}
By the second part of Corollary \ref{te:fc groupoid}, the oblique
arrows induce isomorphisms in cyclic homology. Then, passing to
cyclic homology, also the horizontal map yields an isomorphism. We
can now conclude the proof of the theorem by remarking that, for
any $X\in \mathfrak{X},$
\begin{equation*}
\mathrm{HC}_{\ast }(B_{X})\cong \lim_{{\overrightarrow{X\in \mathfrak{X}}}}
\mathrm{HC}_{\ast }(B_{X}).
\end{equation*}
Thus, taking $X=\{x\}$ we get the required isomorphism.
\end{proof}

The computation performed in Theorem \ref{thm:inf_conn} can be extended to an
arbitrary (discrete) groupoid $\mathcal{G}.$ Let
$(\mathcal{G}_{i})_{i\in I}$ be the connected components of
$\mathcal{G}.$ For each $i$ we pick up $x_{i}\in
\mathcal{G}_{i}^{0}$ and we denote the set of loops $l$ with
$s(l)=x_{i}$ by $G_{i}.$ We have the following result.

\begin{corollary}
\label{co:general} Let $\mathcal{G}$ be a discrete groupoid. If
$B$ denotes the groupoid algebra of $\mathcal{G}$ over a field of
characteristic zero, then
\begin{equation}\label{iso:groupoid}
\mathrm{HC}_{\ast }(B)\cong \bigoplus_{i\in I}\mathrm{HC}_{\ast }(\mathbb{K}
G_{i}).
\end{equation}
{In addition, if $\mathcal{G}^0_i$ is a finite set for every $i\in I$,
  then a similar isomorphism holds in Hochschild homology}.
\end{corollary}

\begin{proof}
For $i\in I,$ let $B_{i}$ be the groupoid algebra of
$\mathcal{G}_{i}.$ As a vector space, $B$ is isomorphic to
$\bigoplus_{i\in I}B_{i}.$ Via this identification, the
multiplication of $B$ can be extended to the direct sum. It is
easy too see that, for two families $(b_{i}^{\prime })_{i\in I}$
and $(b_{i}^{\prime \prime})_{i\in I}$ in $\bigoplus_{i\in
I}B_{i}$, we have
\begin{equation*}
(b_{i}^{\prime })_{i\in I}\cdot (b_{i}^{\prime \prime })_{i\in
I}=(b_{i}^{\prime }b_{i}^{\prime \prime })_{i\in I}.
\end{equation*}
The corollary now follows by the isomorphism $\mathrm{HC}_{\ast
}(\bigoplus_{i\in I}B_{i})\cong \bigoplus_{i\in
I}\mathrm{HC}_{\ast }(B_{i}).$ If $I$ is a finite set,
then this isomorphism can be found in \cite[Exercise
2.2.1]{Loday}. Since cyclic homology and direct limits commute, the
isomorphism can be extended for an arbitrary set $I.$

{Since for non-unital algebras Hochschild homology is
constructed in a different way (cf. \cite[Chapter 1, \S 1.2]{Loday}), 
the above arguments can not be applied to deduce an isomorphism for Hochschild
homology, analogous to \eqref{iso:groupoid}.
Nevertheless, in the case when $\mathcal{G}^0$ is finite, $B$ is
an unital algebra. Thus $\mathrm{HH}_{\ast }(B)$ can be computed
as the Hochschild homology of Connes' cyclic object. The required
isomorphism now follows by \cite[Theorem 9.1.8]{We}, proceeding as
for cyclic homology.}
\end{proof}

As an application of Corollary \ref{co:general} we shall compute
$\mathrm{\mathrm{HH}}_{\ast }(B)$ and $\mathrm{HC} _{\ast }(B)$,
where $B$ is the groupoid algebra of the groupoid associated to a
$G$-set $X$.

Throughout the remaining part of the paper we fix an arbitrary
(discrete) group $G$ that acts to the left on an {arbitrary} set
$X$.
The action of $G$ on $X$ maps a pair $(x,g)\in $ $X\times G$ to
$^{g}x\in X.$ For a $G$-set $X$ as above, one construct a groupoid
$\mathcal{G}$ as follows. By definition, the set of objects in
$\mathcal{G}$ is $\mathcal{G}^{0}:=X$ while, for $x,y\in X, $ we
set
\begin{equation*}
\text{\textrm{Hom}}_{\mathcal{G}}(x,y)=\left\{ (x,g)\in X\times
G\mid {}^{g}x=y\right\} .
\end{equation*}
Note that the source of $(x,g)$ is $x$ and its target is $^{g}x.$
Thus the composition $(x,g)\circ (x^{\prime },g^{\prime })$ is
defined if, and only
if $x={}^{g^{\prime }}x^{\prime }$ and, in this case
\begin{equation*}
(x,g)\circ (x^{\prime },g^{\prime }):=(x',gg^{\prime }).
\end{equation*}
The set of morphisms in $\mathcal{G}$ is $\mathcal{G}^{1}=X\times
G.$ Therefore, the groupoid algebra $B$ of $\mathcal{G}$ has
$X\times G$ as a basis. To describe the multiplication on this
basis let us recall some well-known facts about twisted semigroup
algebras by a $2$-cocycle. Let $S$ be a semigroup. A function
$\omega :S\times S\rightarrow \mathbb{K}$ is
called a $2$-cocycle if, for any $p_{1},$ $p_{2},$ $p_{3}$ in $S$
\begin{equation*}
\omega (p_{1},p_{2})\omega (p_{1}p_{2},p_{3})=\omega
(p_{2},p_{3})\omega (p_{1},p_{2}p_{3}).
\end{equation*}
The semigroup algebra of $S$ is defined as in the group case: as a
vector space it has $S$ as a basis and the multiplication on this
basis is given by the multiplication in $S.$ We shall denote this
algebra by $\mathbb{K}S.$ The cocycle $\omega $ can be used to
deform the multiplication of $\mathbb{K} S$ such that we get
another associative algebra structure on the vector space
$\mathbb{K}S.$ Its multiplication is defined by
\begin{equation*}
p_{1}\cdot p_{2}=\omega (p_{1},p_{2})p_{1}p_{2}.
\end{equation*}
{The resulting algebra will be denoted by $\mathbb{K}_{\omega }S.$
Certainly, it is not unital in general. Still,} even if $S$ has no neutral
element, the algebra $\mathbb{K}_{\omega }S$ may have a unit $e:=\sum_{x\in X}x$,
where $X$ is an appropriate finite subset of $S$.
In fact, {it is easy to see that $e$ is a unit element in the algebra
  $\mathbb{K}_{\omega }$} if, and only if $\omega $ is $X$-normalized,
i.e. for any $p,q\in S$
\begin{equation*}
\sum\limits_{\{x\in X\mid xp=t\}}\omega (x,p)=\delta
_{p,q}=\sum\limits_{\{x\in X\mid px=q\}}\omega (p,x).
\end{equation*}
Let us turn back to the  groupoid algebra of $\mathcal{G},$ where
$X$ is a $G$-set. We define the semigroup $S:=X\times G$ with the
multiplication
\begin{equation*}
(x,g)(y,h)=(y,gh).
\end{equation*}
It is not difficult to see that $\omega_X :S\times S\rightarrow
\mathbb{K},$ given by
\begin{equation}\label{2-cocycle}
\omega_X \big((x,g),(y,h)\big):=\left\{
\begin{array}{cc}
0, & \text{if }^{h}y\neq x \\
1, & \text{if }^{h}y=x
\end{array}
\right.
\end{equation}
{is a $2$-cocycle, which is $X$-normalized if, and only if $X$
is finite}. Obviously, $B=\mathbb{K}_{\omega }S$ {as
non-unital algebras, in general}. For a finite $G$-set $X$, this
is an equality of unital algebras.

We denote the set of {${G}$}-orbits in $X$ by $\mathfrak{X}$.
Let us remark that there is an one-to-one correspondence between
$\mathfrak{X}$ and the set of connected components of ${\mathcal
G}$. This bijection maps an orbit $\mathfrak{0}\in \mathfrak{X}$
to the full subgroupoid $\mathcal{G}_{\mathfrak{o}}$ defined
uniquely such that $\mathcal{G}_{\mathfrak{o}}^{0}=\mathfrak{o}.$
We choose a transversal $\{x_{\mathfrak{o}}\in X\mid
\mathfrak{o}\in \mathfrak{X}\}$ for $\mathfrak{X}.$ Thus, any
$x\in X$ is in the orbit of a certain $x_{\mathfrak{o}}$ and
$x_{\mathfrak{o}^{\prime }}$ and $x_{\mathfrak{o}^{\prime \prime
}}$ are in the same orbit if, and only if $\mathfrak{o}^{\prime
}=\mathfrak{o}^{\prime \prime }.$ Moreover, the set
$G_{\mathfrak{o}}$ of loops $l$ such that $s(l)=x_{\mathfrak{o}}$
is the
stabilizer of  $x_{\mathfrak{o}}$
\begin{equation*}
G_{\mathfrak{o}}:=\{g\in G\mid
{}^{g}x_{\mathfrak{o}}=x_{\mathfrak{o}}\}.
\end{equation*}
Hence, a direct application of Corollary \ref{co:general} yields
the following.

\begin{theorem}
Let $X$ be a $G$-set. If $\omega_X$ denotes the 2-cocycle in
\eqref{2-cocycle} and $\mathbb{K}$ is a field of characteristic
zero, then
\begin{equation*}
\mathrm{HC}_{\ast }(\mathbb{K}_{\omega_X} G)=\bigoplus_{\mathfrak{o}\in
  \mathfrak{X}}\mathrm{HC} _{\ast }(\mathbb{K}G_{\mathfrak{o}}).
\end{equation*}
If in addition $X$ is finite, then
\begin{equation*}
\mathrm{HH}_{\ast }(\mathbb{K}_{\omega_X}
G)\mathbb{=}\bigoplus_{\mathfrak{o}\in \mathfrak{X}} \mathrm{HH}_{\ast
}(\mathbb{K}G_{\mathfrak{o}}).
\end{equation*}
\end{theorem}

\section*{Acknowledgments}
The first author was financially supported by the Hungarian Scientific
Research Fund OTKA K 68195 and the Bolyai J\'anos
Scholarship, while the second author was supported by Contract 2-CEx06-11-20 of
the Romanian Ministry of Education and Research. Both authors are grateful to
Bachuki Mesablishvili, Zoran \v Skoda and the referee for their helpful
comments.

\end{document}